\documentclass[11pt,twoside]{article}

\usepackage{multirow}
\usepackage{graphicx}
\usepackage{amsfonts}
\usepackage{amsmath}
\usepackage{amssymb}
\usepackage{colortbl}
\usepackage{indentfirst}
\usepackage{booktabs}
\usepackage{verbatim}
\usepackage{color}
\usepackage{amsthm}
\usepackage{booktabs}
\usepackage{caption}
\usepackage{pdflscape}
\usepackage{algorithm}
\usepackage[noend]{algpseudocode}
\usepackage{enumitem}

\usepackage[page,header]{appendix}

\newcommand{\bP}{\mathbf{P}}
\newcommand{\bE}{\mathbf{E}}

\newcommand{\BC}{{\mathbb{C}}}

\newcommand{\BE}{{\mathbb{E}}}

\newcommand{\BN}{{\mathbb{N}}}

\newcommand{\BR}{{\mathbb{R}}}

\newcommand{\CD}{{\mathcal{D}}}
\newcommand{\CE}{{\mathcal{E}}}

\newcommand{\CJ}{{\mathcal{J}}}

\newcommand{\CM}{{\mathcal{M}}}
\newcommand{\CN}{{\mathcal{N}}}
\newcommand{\CO}{{\mathcal{O}}}
\newcommand{\CP}{{\mathcal{P}}}
\newcommand{\CQ}{{\mathcal{Q}}}

\newcommand{\CY}{{\mathcal{Y}}}
\newcommand{\CZ}{{\mathcal{Z}}}

\newcommand{\ep}{\varepsilon}
\newcommand{\beq}{\begin{equation}}
\newcommand{\eeq}{\end{equation}}
\newcommand{\bde}{\begin{definition}}
\newcommand{\ede}{\end{definition}}
\newcommand{\brm}{\begin{remark}}
\newcommand{\erm}{\end{remark}}
\newcommand{\bex}{\begin{example}}
\newcommand{\eex}{\flushright{\(\qedsymbol\)}\end{example}}
\newcommand{\bthm}{\begin{theorem}}
\newcommand{\ethm}{\end{theorem}}
\newcommand{\bprp}{\begin{proposition}}
\newcommand{\eprp}{\end{proposition}}
\newcommand{\blm}{\begin{lemma}}
\newcommand{\elm}{\end{lemma}}
\newcommand{\ol}{\overline}
\newcommand{\bprf}{\begin{proof}}
\newcommand{\eprf}{\end{proof}}
\newcommand{\ben}{\begin{eqnarray}}
\newcommand{\een}{\end{eqnarray}}

\newcommand{\dd}{\mathrm{d}}

\newcommand{\kl}{\left(}
\newcommand{\kr}{\right)}
\newcommand{\bmr}{\left(\begin{array}}
\newcommand{\emr}{\end{array}\right)}

\newcommand{\bal}{\begin{aligned}}
\newcommand{\eal}{\end{aligned}}
\newcommand{\ti}{\mathrm{in}}
\newcommand{\fga}{\mathrm{FGA}}
\newcommand{\fgs}{\mathrm{FGS}}
\newcommand{\ufga}{u_{\mathrm{FGA}}}
\newcommand{\ufgs}{u_{\mathrm{FGS}}}
\newcommand{\ii}{\mathrm{i}}
\newcommand{\ee}{\mathrm{e}}
\newcommand{\rsp}{\mathrm{sp}}

\numberwithin{equation}{section}
\newtheorem{theorem}{Theorem}[section]
\newtheorem{lemma}[theorem]{Lemma}
\newtheorem{corollary}[theorem]{Corollary}
\newtheorem{proposition}[theorem]{Proposition}
\newtheorem{definition}[theorem]{Definition}
\newtheorem{remark}[theorem]{Remark}

\newtheorem{assumption}{Assumption}

\begin{document}

	\title{Frozen Gaussian Sampling: A Mesh-free Monte Carlo Method For Approximating Semiclassical Schr\"odinger Equations}
	
	\author{Yantong Xie\footnote{School of Mathematics Science, Peking University, Beijing, 100871, China (darkoxie@pku.edu.cn).} \   \
	Zhennan Zhou\footnote{Beijing International Center for Mathematical Research, Peking University, Beijing, 100871, China (zhennan@bicmr.pku.edu.cn).}}
	\maketitle

\begin{abstract}

In this paper, we develop a Monte Carlo algorithm named the Frozen Gaussian Sampling (FGS) to solve the semiclassical Schr\"odinger equation based on the frozen Gaussian approximation. Due to the highly oscillatory structure of the wave function, traditional mesh-based algorithms suffer from "the curse of dimensionality", which gives rise to more severe computational burden when the semiclassical parameter \(\ep\) is small. The Frozen Gaussian sampling outperforms the existing algorithms in that it is mesh-free in computing the physical observables and is suitable for high dimensional problems. In this work, we provide detailed procedures to implement the FGS for both Gaussian and WKB initial data cases, where the sampling strategies on the phase space balance the need of variance reduction and sampling convenience. Moreover, we rigorously prove that, to reach a certain accuracy, the number of samples needed for the FGS is independent of the scaling parameter \(\ep\). Furthermore, the complexity of the FGS algorithm is of a sublinear scaling with respect to the microscopic degrees of freedom and, in particular, is insensitive to the dimension number. The performance of the FGS is validated through several typical numerical experiments, including simulating scattering by the barrier potential, formation of the caustics and computing the high-dimensional physical observables without mesh.

\end{abstract}

{\small

}


\section{Introduction}

The Schr\"odinger equation, which depicts the evolution of the wave function of a quantum system, is given as follows:
\beq\label{physch}
\ii \hbar \frac{\partial}{\partial{t}}u(t,x)=-\frac{\hbar^{2}}{2 m_{\text{o}}}\Delta u(t,x)+E(x)u(t,x),
\eeq
where \(m_{\text{o}}\) denotes the mass of the particles and \(\hbar\) denotes the reduced Planck constant. By adopting the nondimensionalisation techniques (see Appendix A), one obtains the following semiclassical Schr\"odinger equation in the dimensionless form: 
\begin{equation}
\label{SemiSchrodinger}
\begin{cases}
\ii\varepsilon\frac{\partial}{\partial{t}}u(t,x)=-\frac{\varepsilon^2}{2}\Delta u(t,x)+E(x)u(t,x),&\quad(t,x)\in\mathbb{R}\times\mathbb{R}^m,\\
u(0,x)=u_{\ti}(x),&\quad x\in\BR^m,
\end{cases}
\end{equation}
where the parameter \(\ep\), named the semiclassical parameter, describes the microscopic macroscopic scale ratio and is small  within the semiclassical regime (i.e. \(0<\ep\ll 1\)).  The quantum mechanical wave function \(u(t,x)\), valued in \(\BC\), often serves as an auxiliary quantity to compute the physical quantities of interest. The physical quantities, including position, momentum and energy, etc., are referred to as "observables", which are of the form of quantum averages with respect to the wave functions of certain self-adjoint operators:
\beq
O=\langle u|\hat{O}|u\rangle= \int_{\mathbb{R}^{m}} \overline{u(t,x)}\,\hat{O}\, u(t,x)
\dd x.
\eeq
Here, \(\hat{O}\) denotes the operator related to a physical observable. For examples, \(\hat{O}^q=x\) and \(\hat{O}^p=-\ii\ep\nabla\) are the operators corresponding to the position and momentum observables, respectively.

In the semiclassical regime, namely \(\ep\ll 1\), the solution to equation \eqref{SemiSchrodinger} is highly
oscillatory both in space and time on the scale \(\CO(\ep)\) and doesn't converge in strong sense as \(\ep\to 0\), which gives rise to significant computational burdens \cite{JMS2011,LL2020}. Nevertheless, various effective numerical algorithms have been proposed so far, succeeding in capturing the highly oscillatory wave functions accurately with affordable computational complexity. The existing algorithms can be divided into two classes: direct discretization methods and asymptotic methods. Among the direct discretization methods, one of the most prevailing choices is the time splitting spectral method proposed in \cite{BJM2003}, where the meshing strategy \(\Delta t=\CO(\ep)\) and \(\Delta x=\CO(\ep)\) is practically sufficient for accurate approximation of the wave function. The time step size of the time splitting spectral method can be even relaxed to \(\CO(1)\), in other words, independently of \(\ep\), if one is only interested in capturing the correct physical observables, whereas one still needs to resolve the spatial oscillations. The asymptotic methods, including the WKB methods \cite{ER1996,ER2003}  and  various wave packets based methods \cite{TQR2007,JWY2008,LY2011-1,LY2011-2,LY2011-3,RZ2019,JLRZ2019,MRZ2020}, provide accurate numerical approximation to the wave function with a relatively less computational cost with respect to the parameter \(\ep\) compared to the direct discretization methods. The wave packets based methods have been extended to more complicated quantum dynamics, such as the vector potential cases and systems with random inputs \cite{Z2014,RZ2019,JLRZ2019,MRZ2020}, providing broad application prospects in many applications.

Although the existing algorithms have made great progress in efficiently simulating the wave functions and the physical observables in the semiclassical regime, they still suffer from the "the curse of dimensionality" when considering high dimensional cases, since those methods necessarily reply on an Eulerian mesh either in the initial decomposition or in the time evolution. When the dimension number \(m\gg 1\), the computation and memory amounts of the both classes of methods increase exponentially as $m$ increases, which often leads to more infeasible computational cost when $\ep \ll 1$. For example, for the time splitting spectral method, the degree of freedom to approximate the wave function on an \(m-\)dimensional mesh is \(\CO(\ep^{-m})\), not to mention the increased computational complexity to evolve the wave function on this high dimensional mesh. Therefore, a novel numerical method for the Schr\"odinger equation in semiclassical regime, whose cost is insensitive to \(\ep\) and the dimension number $m$, is highly desirable.

To the best of our knowledge, there has been very few existing methods, which are feasible or have the potential to be suitable for high dimensional simulations. From this perspective, the path integral sampling method proposed recently in \cite{LZ2016,LZ2017} can be viewed as a pioneering work, which is constructed on the so-called Frozen Gaussian approximation (FGA). The FGA represents the wave function with a superposition of Gaussian functions in the phase space with fixed widths and is often wrote in the form of an integral on the phase space:
\begin{equation}\label{absans}
\ufga(t,x)=\frac{1}{(2\pi\varepsilon)^{\frac{3m}{2}}}\int_{\mathbb{R}^{2m}}
A(t,q,p)e^{\frac{i\Theta (t,x,q,p)} {\varepsilon}}\dd q\dd p,
\end{equation}
where the variables \(A\) and \(\Theta\) are can be approximately obtained through evolving a system of Lagrangian ODEs, thus the full quantum dynamics is reduced to an ensemble of wave packets dynamics. There have been quite a few works related to the FGA, including the theoretical studies in the chemistry community \cite{H1981,HK1984,K1994,K2006}, the rigorous convergence analysis  \cite{SR2009} and the FGA based algorithms for solving hyperbolic systems  which reply on the direct phase space discretization \cite{LY2011-1,LY2011-2,LY2011-3}. However, it is not until \cite{LZ2016,LZ2017} that people find another distinct prospective, giving the FGA ansatz \eqref{absans} a probabilistic interpretation, and show that the Monte-Carlo sampling method based on the FGA gives a feasible algorithm which is mesh-free in computing the physical observables. The main spirit of \cite{LZ2016,LZ2017} can be summarized as follows: firstly, as an integral on the phase space, the FGA ansatz can be viewed as an expectation of random variables on the phase space, which can be calculated by solving the parameter ODEs; in addition, when the initial data is Gaussian, the module function \(|A(0,q,p)|\) is also Gaussian and naturally serves as an ideal choice of a sampling density function on the phase space. However, what remains to be done is that there have been no rigorous numerical analysis on the sampling error nor discussion on how to extend the sampling method to general cases when the initial data is not Gaussian. Our current work aims to fill the vacancy on these issues.

In this work, we develop a Monte Carlo method which relies on the probabilistic interpretation on the FGA ansatz, and we shall name it the Frozen Gaussian Sampling (FGS). The essence of the FGS is to properly choose a reference density function on the phase space, and reformulate the FGA representation \eqref{absans} as an expectation of a complex-valued functional with respect to the Frozen Gaussian trajectories. In this view, the FGS is a path-integral approximation of the wave function, and in the dynamical picture, the evolution of the quantum state is replaced by propagating an ensemble of Frozen Gaussian wave packets, where the initial phase space coordinates are sampled from the reference density function.  The design of an efficient FGS algorithm boils down to proposing a suitable reference sampling measure, which is closely related the form of the initial condition of the wave function, and the quality of an FGS algorithm is measured by the variance of the functional to be sampled. Hence, the FGS provides a subtle shift in the computation paradigm and the analysis framework, which, as we shall elaborate, turns to show great promises as a simulation tool for high dimensional problems.

We focus on two of  the most important types of initial quantum states: the Gaussian wave packet and the WKB  data, and study the respective strategies in designing the density functions that balance the need of variance reduction and sampling convenience. 
For the Gaussian initial condition scenario, we rigorously prove that, with the feasible choice of the sampling density, the variance of the functional to be sampled is independent of the scaling parameter \(\ep\) . Furthermore, the number of samples to reach a certain sampling error $\delta$ is in the order $\CO \left(c^{m} \delta^{-2}  \right)$, where $c$ is a $\CO(1)$ constant and $m$ is the dimension number.  In particular, if one aims to halve the sampling error, only four times of the samples are needed, which is independent of the dimension $m$ or the semiclassical parameter $\varepsilon$. Different from all the existing methods, this result shows that the number of wave packets does not to "resolve" any spatial oscillations, and is of a sublinear scaling of the microscopic degrees of freedom. Besides, the computation cost grows rather gently as $m$ increases, and the whole FGS algorithm is parallelizable, which makes high dimensional simulations rather affordable.   
For WKB initial data cases, we propose a sequential sampling strategy leveraging the mono-kinetic structure of the highly oscillatory wave function. 
We numerically verify that the sample size to reach a certain error is also insensitive to \(\ep\), which still outperforms the existing mesh-based algorithms whose computational cost grows rapidly as \(\ep\to0\). 

It is worthy emphasizing that, we only require that the initial condition to be in certain forms to carry out efficient initial sampling, but we do not require that the solutions remain in such ansatz forms. This is a natural consequence of the construction of the FGS algorithm, while most asymptotic methods break down when the solution ansatz is no longer valid. And we have carried a set of 2D numerical tests to illustrate such a property, where the initial Gaussian wave packet is scattered into various profiles by different potential barriers, and the FGS method is able to capture the correct evolution behavior with a reasonably small sample size.  

Another main highlight of the FGS lies in that it is mesh-free in computing the physical observables. By taking advantage of path integral representation, the FGS succeeds in computing a wide class of physical observables without reconstructing the wave functions. 
Therefore, if one is only interested in physical observables rather than the detained behavior of the wave function, the FGS serves as an efficient algorithm to simulate high dimensional Schr\"odinger equation, since it is a mesh-free algorithm of a sublinear scaling complexity and is naturally compatible with parallel computing. To demonstrate the performance of the FGS in computing observables, we implement the FGS to compute the  physical observables in high  dimensional quantum systems with a harmonic potential and it is shown numerically that the computational cost to reach a certain error is rather affordable in various dimensions. With the sample size quadrupling when the dimension number increases by one, we observe that the errors in the observables is uniformly bounded in all dimensions tested.

The rest of this paper is organized as follows: in Section 2, we review the basic setup the FGA and give the detailed procedures to compute the wave functions and physical observables with the FGS. We also present the error analysis framework of the FGS and elaborate the difficulties to be overcome when choosing the density function for the FGS. In Section 3, we focus on Equation \eqref{SemiSchrodinger} with Gaussian initial data and prove that the quantitative estimates of the FGS sampling error . In Section 4, we propose a feasible sampling strategy to apply the FGS on Equation \eqref{SemiSchrodinger} with WKB initial data. In Section 5, the numerical performance of the FGS is validated through five carefully designed numerical experiments, including simulating the Gaussian barrier potential and the caustics and computing the high-dimensional position and momentum observables without mesh.

\section{The frozen Gaussian approximation and the frozen Gaussian sampling}

In this section, we briefly review the forzen Gaussian approximation (FGA) to the semiclassical Schr\"odinger equation. The frozen Gaussian approximation is in the form of an integral on the phase space, which can be given a probabilistic interpretation. Then we present the frozen Gaussian sampling (FGS) algorithm as a general computation methodology for semiclassical quantum dynamics and discuss the practical implementation of computing physical observables based on the FGS algorithm.

\subsection{The frozen Gaussian representation for the semiclassical Schr\"odinger equation}

The frozen Gaussian approximation is an asymptotic approximation to the solution of Equation \eqref{SemiSchrodinger} with \(\CO(\ep)\) error \cite{SR2009}. It is based on an integral representation on the phase space:
\begin{equation}\label{FGAansatz}
\ufga(t,x)=\frac{1}{(2\pi\varepsilon)^{\frac{3m}{2}}}\int_{\mathbb{R}^{2m}}
A(t,q,p)\ee^{\frac{\ii\Theta (t,x,q,p)} {\varepsilon}}\dd q\dd p.
\end{equation}
Here the phase function \(\Theta\) is given by
\begin{eqnarray}
&&\label{ThetaFGA}\Theta(t,x,q,p)=S(t,q,p)+P(t,q,p)\cdot(x-Q(t,q,p))+\frac{\ii}{2}|x-Q(t,q,p)|^{2}.
\end{eqnarray}

For a given initial condition \(q\) and \(p\), the evolution of \(Q\) and \(P\) is governed by the classical Hamitonian \(H(q,p)=E(q)+\frac{1}{2}|p|^2\) and defines a trajectory on phase space \(\BR^{2m}\) called an \textbf{FGA tracjectory}, which also determines the value of the amplitude \(A\) and action \(S\). We name the variables \(Q,P,S,A\) as \textbf{FGA variables}. The evolution of the FGA variables is determined by the ODE system as follows:
\begin{eqnarray}
\label{PhaseQ}\frac{\dd}{\dd t} Q(t,q,p)&=&P,\\
\label{PhaseP}\frac{\dd}{\dd t} P(t,q,p)&=&-\nabla E(Q),\\
\label{PhaseS}\frac{\dd}{\dd t}S(t,q,p) &=&\frac{1}{2}|P|^{2}-E(Q),\\
\label{PhaseA}\frac{\dd}{\dd t}A(t,q,p)&=&\frac{A}{2} \mathrm{tr}\left(Z^{-1}\left(\partial_zP-\ii\partial_zQ\nabla_Q^2E(Q)\right)\right),
\end{eqnarray}
where we use the short hand notations
\[
\partial_z=\partial_q-\ii\partial_p\quad\text{and}\quad Z=\partial_z(Q+\ii P).
\]
Moreover, the initial condition is given by
\begin{eqnarray}
\label{PhaseQ0}Q(0,q,p)&=&q,\\
\label{PhaseP0}P(0,q,p)&=&p,\\
\label{PhaseS0}S(0,q,p)&=&0,\\
\label{A0Compute}A(0,q,p)&=&2^{\frac{m}{2}}\int_{\mathbb{R}^m}u_{\ti}(y)e^{\frac{\ii}{\varepsilon}(-p\cdot(y-q)+\frac{\ii}{2}|y-q|^{2})} \dd y.
\end{eqnarray}
The derivation of the ODE system \eqref{PhaseQ}-\eqref{PhaseA} is obtained through asymptotic matching. See \cite{SR2009} for detailed discussion on this issue. To summarize, with the help of the FGA, one can calculate the asymptotic solution to semiclassical Schr\"odinger equation \eqref{SemiSchrodinger} through solving an ensemble of ODE systems of the FGA variables instead of a high dimensional PDE. 

In addition, we remark that the phase space integral \eqref{FGAansatz} can also be viewed as an integral of the FGA trajectory. Denote the FGA tracjectory by \(z_t=(Q(t),P(t))\), where we hide the dependence of the FGA variables on \((q,p)\) for simplicity (we use the similar short hand notation in the rest of the paper). Then \(z_t\) satisfies the Hamiltonian flow and defines a one-to-one map \(z_0\mapsto z_t\) on the phase space. In other words, the FGA trajectory \(z_t\) is specified once the initial condition \((q,p)\) is fixed and the FGA variables \(S\) and \(A\) are also specified with a given trajectory \(z_t\). The idea of viewing the FGA ansatz \eqref{FGAansatz} through FGA trajectory is important in the introduction of the frozen Gaussian sampling.

To conclude this subsection, we directly quote two theorems  from \cite{LY2011-1,SR2009}. The first theorem shows the FGA ansatz integration \eqref{FGAansatz} reproduces the initial condition:
\bthm\label{Thmini}
For \(u_\ti\in L^2(\BR^{m})\), we have
\beq
u_\ti(x)= \frac{1}{(2\pi\varepsilon)^{\frac{3m}{2}}}\int_{\mathbb{R}^{2m}}
A(0,q,p)\ee^{\frac{\ii\Theta (0,x,q,p)} {\varepsilon}}\dd q\dd p,
\eeq
where the initial amplitude \(A(0,q,p)\) given in Equation \eqref{A0Compute} is determined by \(u_\ti\).
\ethm

The second theorem shows that the \(L^2\) error of \(\ufga\) with respect to the exact solution is \(\CO(\ep)\):  
\bthm\label{ThmOep}
Consider the semiclassical Schr\"odinger equation \eqref{SemiSchrodinger} with the exact solution \(u(t,x)\). Assume that the potential function \(E(x)\in C^\infty(\BR^m)\). For any \(t>0\), we have the following \(L^2\) error estimate result for the FGA ansatz \eqref{FGAansatz}:
\begin{equation}
\left\|\ufga(t,\cdot)-u(t,\cdot)\right\|_{L^2}\leq C_A(t,m)\ep\left\|u_{\ti}(\cdot)\right\|_{L^2},
\end{equation}
where \(C_A(t,m)\) is a parameter independent of the scaling parameter \(\ep\).
\ethm

\subsection{The frozen Gaussian sampling algorithm}

We have seen in Section 2.1 that the FGA ansatz \eqref{FGAansatz} is an integral on the phase space. With the FGA, the task of solving a PDE is converted into solving an ensemble of ODEs \eqref{PhaseQ}-\eqref{PhaseA}. However, difficulties still remain in computing the integral \eqref{FGAansatz}, which can be of high dimensions. Direct mesh method of numerical integration may suffer from "the curse of dimensionality" and does not lead to a practical algorithm. To deal with the high-dimensional integral \eqref{FGAansatz}, we seek a stochastic method for calculation. To this end, we write the FGA ansatz in the form of an expectation and approximate the integral \eqref{FGAansatz} numerically with a Monte Carlo method. The resulting stochastic algorithm is called the \textbf{frozen Gaussian sampling} (FGS), which also gives a mesh-free implement to compute physical observables. In this subsection, we elaborate the FGS algorithm and discuss its advantages and difficulties to be further investigated.

To write the FGA ansatz \eqref{FGAansatz} into an expectation, we introduce a real-valued probability density function \(\pi(\cdot)\) on the phase space \(\BR^{2m}\). The density \(\pi\) then defines a probability measure \(\bP_0\) on \(\BR^{2m}\):
\begin{equation}
\mathbf{P}_0(\Omega)=\int_{\Omega}\pi(z_0)\dd z_0,\quad\forall\Omega\subset\BR^{2m},
\label{measure}
\end{equation}
where \(\bP_0(\BR^{2m})=\int_{\BR^{2m}}\pi(z_0)\dd z_0=1\). With the probability measure, we can rewrite the FGA ansatz \eqref{FGAansatz} as 
\begin{eqnarray}
\ufga(t,x)&=&\frac{1}{(2\pi\varepsilon)^{\frac{3m}{2}}}\int_{\mathbb{R}^{2m}}\pi(z_0)\frac{A(t,z_0)}{\pi(z_0)} \ee^{\frac{\ii\Theta (t,x,z_0)}{\varepsilon}}\dd z_0\nonumber\\
&=&\bE_{z_0\sim\pi}\left[\frac{1}{(2\pi\ep)^{\frac{3m}{2}}}\frac{A(t,z_0)}{\pi(z_0)}\ee^{\frac{\ii\Theta(t,x,z_0)}{\varepsilon}}\right].
\label{ProbInterrupt}
\end{eqnarray}

Thus, we may use a Monte Carlo sampling for \(\ufga(t,x)\):
\begin{equation}\label{Probz0j}
\ufga(t,x)\approx \ufgs\left(t,x;M,\left\{z_0^{(j)}\right\}_{j=1}^M\right)=\frac{1}{M}\sum_{j=1}^M \Lambda\left(t,x,z_0^{(j)}\right),
\end{equation}
where
\beq\label{Lambda}
\Lambda\left(t,x,z_0\right)=\frac{1}{(2\pi\ep)^{\frac{3m}{2}}}\frac{A\kl t,z_0\kr}{\pi\kl z_0\kr} \exp\kl\frac{\ii\Theta\kl t,x,z_0\kr}{\varepsilon}\kr.
\eeq
Here \(\left\{z_0^{(j)}\right\}_{j=1}^M\) are \(M\) independent identically distributed samples with the probability density function \(\pi\). We refer to \(\ufgs\) as the frozen Gaussian sampling (FGS) wave function. We also note that \(\Lambda\left(t,x,z_0\right)\) can be treated as a random variable related to \(z_0\) with the expectation
\beq
\bE_{z_0\sim\pi}\Lambda\kl t,x,z_0\kr=\ufga(t,x).
\eeq

Algorithmically, once we have sampled \(z_0\), the FGA trajectory and FGA variables are available through evolving the ODE system Equation \eqref{PhaseQ}-\eqref{PhaseA} up to time \(t\). With enough samples of the FGA trajectories, we approximates the FGS wave function \(u_\fgs\) numerically by Equation \eqref{Probz0j}. From the discussion above, we give the pseudo-code of the frozen Gaussian sampling (FGS) in Algorithm \ref{alg:FGA}. For simplicity, in the rest of this paper, we use the short hand notations to the FGA trajectories and FGA variables:
\begin{eqnarray}
Q_j(t)=Q\kl t,z_0^{(j)}\kr,\quad P_j(t)=P\kl t,z_0^{(j)}\kr,&&\nonumber\\ S_j(t)=S\kl t,z_0^{(j)}\kr,\quad A_j(t)=A\kl t,z_0^{(j)}\kr,&&\quad\forall j=1,\cdots,M.\nonumber
\end{eqnarray}
We note that analogous short hand notation may be also applied to other variables without confusion.

\begin{algorithm}
\caption{The frozen Gaussian sampling} \label{alg:FGA}
\begin{algorithmic}[1]
\State \textbf{Initial sampling:} sample \(M\)  independent identically distributed initial conditions \(\left\{z_0^{(j)}\right\}_{j=1}^M\) that obey the density function \(\pi\) on the phase space \(\BR^{2m}\).
\State \textbf{Initial amplitude computing:} for \(j=1:M\), compute the initial amplitude \(A_j(0)\) by Equation \eqref{A0Compute}.
\State \textbf{Trajectories evolution:} for \(j=1:M\), evolve the ODE system Equation \eqref{PhaseQ}-\eqref{PhaseA} up to time \(t\) with the initial condition given by Equation \eqref{PhaseQ0}-\eqref{A0Compute} to obtain the FGA variables \(Q_j(t),P_j(t),S_j(t),A_j(t)\).
\State \textbf{Reconstruction:} compute the wave function or physical observable based on Equation \eqref{Probz0j} or Equation \eqref{EqDoubleS} respectively.
\end{algorithmic}
\end{algorithm}

However, there are certain difficulties when actually executing the FGS algorithm, especially in high dimensional cases. First, to compute the initial amplitude, one needs to calculate the integral \eqref{A0Compute}. For a Gaussian wave packet initial condition \(u_\ti\), the initial amplitude \(A(0,q,p)\) can be calculated explicitly in Equation \eqref{modA}, which has been pointed out in \cite{LZ2016}. However, for a more complex initial condition \(u_\ti\) such as a WKB initial condition, the explicit expression of \(A(0,q,p)\) is unavailable. To compute the initial amplitude \eqref{A0Compute}, the computational cost of numerical integration is huge due to the high dimension and high frequency oscillation structure of the integrand. 

Next, let us discuss the error analysis framework of the FGS wave function \(\ufgs\). Obviously, the \(L^2\) error function \(E_0=\left\|\ufgs-u\right\|_{L^2}\) of the FGS wave function with respect to the exact solution can be divided it into two parts:
\beq\label{Frame}
\left\|\ufgs-u\right\|_{L^2}\leq\left\|\ufga-u\right\|_{L^2}+\left\|\ufgs-\ufga\right\|_{L^2}=E_{\fga}+E_{S},
\eeq
Here \(E_\fga=\left\|\ufga-u\right\|_{L^2}\) denotes the asymptotic error of the FGA ansatz and \(E_S=\left\|\ufgs-\ufga\right\|_{L^2}\) denotes the Monte Carlo sampling error. Thanks to Throrem \ref{ThmOep}, we arrive at 
\beq\label{EFGA}
E_\fga=\left\|\ufga-u\right\|_{L^2}\leq C_A(t,m)\ep,
\eeq
where \(C_A(t,m)\) is a constant depending on time \(t\) and dimension \(m\) and independent of the rescaled parameter \(\ep\). However, with unreasonable choice of density function \(\pi\), the Monte Carlo sampling error \(E_S\) may increase dramatically as the dimension number \(m\) increases, thus formidably many samples are needed. Therefore, to apply the FGS in high dimensional cases, variance reduction techniques such as importance sampling are necessary to reduce the number of samples. Ideally speaking, we aim to choose a reference density \(\pi\) which not only leads to a minimized sampling variance, but also easy to sample numerically. If we solely consider the variance reduction prospective, an obvious choice of sampling density can be given as follows:
\beq\label{piabsA0}
\pi(z_0)=\frac{1}{\int_{\BR^{2m}}|A(0,z)|\dd z}|A(0,z_0)|.
\eeq
For a Gaussian wave packet initial condition, the density function given by Equation \eqref{piabsA0} actually defines a multivariate normal distribution and the resulting FGS algorithm is proven numerically practical and effective \cite{LZ2016,LZ2017}. For a more complex WKB initial condition when the explicit expression of \(|A(0,z_0)|\) in Equation \eqref{piabsA0} is unavailable, how to choose \(\pi(z_0)\) is still unclear. As far as we know, no rigorous numerical analysis on sampling error \(E_S\) has been carried out prior to this work.

The aim of the rest of this paper is to investigate and overcome the above difficulties of the FGS. In Section 3, we prove rigorously that the sampling error \(E_S\) of FGS with sampling density \eqref{piabsA0} is independent of \(\ep\) with a Gaussian initial condition. In Section 4, we introduce an asymptotic expression approximating the initial amplitude \eqref{A0Compute} when the initial data is of a WKB form and present a new initial sampling method of the FGS algorithm based on this approximation. In Section 5, we show various numerical examples to test the practical performance of the FGS with these new techniques, especially in high dimensional cases.

\subsection{Computing observables without mesh}

For many applications, the goal is not to approximate
the wave function in semiclassical Schr\"odinger equation \eqref{SemiSchrodinger} itself, but rather computing certain observables. To calculate an observable, reconstructing the wave function on a mesh is often neither feasible nor necessary, especially for high dimensional cases. In this subsection, we introduce a mesh-free method to compute a wide class of observables, which only use the FGA trajectories sampled by the FGS and doesn't require the reconstruction of the FGS wave function. The method exhibits great advantages on computing high dimensional problems.

Let \(O=\langle u|\hat{O}|u\rangle\) where \(\hat{O}\) denotes a self-adjoint operator asscoiated with a physical observable. Then we compute \(O\) based on the FGS wave function \(\ufgs\) as follows:
\begin{eqnarray}\label{EqDoubleS}
O_\fga&=&\left\langle\ufgs\kl t;M,\left\{z_0^{(j)}\right\}_{j=1}^M \kr\right|\hat{O}\left|\ufgs\kl t;M,\left\{z_0^{(j)}\right\}_{j=1}^M\kr\right\rangle\nonumber\\
&=&\frac{(2\pi\varepsilon)^{-3m}}{M^2}\sum_{j,k=1}^M\frac{\ol{A_j(t)}}{\kl\pi\kl z_0^{(j)}\kr\kr}\cdot \frac{A_k(t)}{\kl\pi\kl z_0^{(k)}\kr\kr}g_{j,k}\kl t,x\kr,
\end{eqnarray}
where
\beq
g_{j,k}\kl t,x\kr=\left\langle \ee^{\frac{\ii\Theta_j(t,x)}{\varepsilon}}\right|\hat{O}\left|\ee^{\frac{\ii\Theta_k(t,x)}{\varepsilon}}\right\rangle.
\eeq

We observe that the only dependence on \(x\) in Equation \eqref{EqDoubleS} is in the term \(g_{j,k}\), which is an integration over the Gaussian wave packets. For a wide class of observables \(\hat{O}\), the term \(g_{j,k}\kl t,x\kr\) can be carried out analytically through Gaussian integration formula, hence the observables \(O=\langle u|\hat{O}|u\rangle\) can be computed without mesh through computing the double summation in Equation \eqref{EqDoubleS}. This class of observables include functions of the position variable and the momentum variable, which is in a polynomial form or can be approximated by  a polynomial.

Now we take the position observable and the momentum observable as examples to show the analytical expression of \(g_{j,k}\). For the position observable \(\hat{O}=x\), we obtain
\begin{eqnarray}\label{ObsP}
g_{j,k}\kl t,x\kr=(\pi\ep)^{\frac{m}{2}}\left[\frac{1}{2}\kl Q_k(t)+Q_j(t)\kr+\frac{\ii}{2}\kl P_k(t)-P_j(t)\kr\right]\times \ee^{\frac{i}{\ep}\kl S_k(t)-S_j(t)\kr}\nonumber\\
\times \ee^{-\frac{i}{2\ep}\kl Q_k(t)-Q_j(t)\kr\cdot\kl P_k(t)+P_j(t)\kr}\times \ee^{-\frac{1}{4\ep}\left| Q_k(t)-Q_j(t)\right|^2-\frac{1}{4\ep}\left|P_k(t)-P_j(t)\right|^2}.
\end{eqnarray}
For the momentum observable \(\hat{O}=-\ii\ep\nabla_x\), we obtain
\begin{eqnarray}\label{ObsM}
g_{j,k}\kl t,x\kr=(\pi\ep)^{\frac{m}{2}}\left[\frac{\ii}{2}\kl Q_j(t)-Q_k(t)\kr+\frac{1}{2}\kl P_k(t)+P_j(t)\kr\right]\times \ee^{\frac{\ii}{\ep}\kl S_k(t)-S_j(t)\kr}\nonumber\\
\times \ee^{-\frac{\ii}{2\ep}\kl Q_k(t)-Q_j(t)\kr\cdot\kl P_k(t)+P_j(t)\kr}\times \ee^{-\frac{1}{4\ep}\left| Q_k(t)-Q_j(t)\right|^2-\frac{1}{4\ep}\left|P_k(t)-P_j(t)\right|^2}.
\end{eqnarray}

To sum up, given the FGA trajectories and FGA variables, we can reconstruct the observables in a mesh-free fashion without computing the FGS wave function. Clearly, the complexity for the physical observable evaluation \eqref{EqDoubleS} is \(\CO(M^2)\). In the following sections, we shall further explore the dependence of the ensemble size \(M\) on the scaling parameter \(\ep\) and the dimension number \(m\) and show great potential of the FGS in computing physical observables in high-dimensional cases.

\section{Frozen Gaussian sampling with Gaussian initial conditions}

In Section 3, we consider a \(m-\)dimension normalized semiclassical Gaussian wave packet centering at \(\tilde{q}\) with momentum \(\tilde{p}\):
\begin{equation}
\label{InitialGaussian}
u_{\ti}(x)=\left(\prod_{j=1}^m a_j\right)^\frac{1}{4}(\pi\varepsilon)^{-\frac{m}{4}}\exp \left(\frac{\ii}{\varepsilon} \tilde{p} \cdot(x-\tilde{q})\right) \exp \left(-\sum_{j=1}^{m} \frac{a_{j}}{2 \varepsilon}\left(x_{j}-\tilde{q}_{j}\right)^{2}\right),
\end{equation} 
where \(\tilde{p},\tilde{q}\in\mathbb{R}^m\) and \(a_j(1\leq j\leq m)\) are positive.

The semiclassical Gaussian wave packet is one of the simplest yet most important wave function in quantum dynamics, whose semiclassical limit arrives at a delta function on the phase space. With the Gaussian integration formula, we are able to calculate the initial amplitude \(A(0,q,p)\) (see Equation \eqref{A0Compute}) explicitly and the resulting density function \(\pi\) as shown in Equation \eqref{piabsA0} is Gaussian. Finally, we rigorously prove the sampling error \(E_S\) of the FGS with density function is independent of \(\ep\) and discuss its dependence on the dimension number \(m\). 

\subsection{Initial sampling and initial amplitude for Gaussian initial condition}

Plug the Gaussian initial condition \eqref{InitialGaussian} into Equation \eqref{A0Compute}, we calculate \(A(0,q,p)\) analytically with the Gaussian integration formula:
\begin{eqnarray}
A(0,q,p)=2^{m}(\pi\varepsilon)^{\frac{m}{4}}\prod_{j=1}^{m}\left[\left(\frac{\sqrt{a_j}}{1+a_{j}}\right)^{\frac{1}{2}} \exp \left(-\frac{\left(\tilde{p}_{j}-p_{j}\right)^{2}+a_{j}\left(\tilde{q}_{j}-q_{j}\right)^{2}}{2\left(1+a_{j}\right) \varepsilon}\right)\right.\nonumber\\ \left. \times\exp \left(\frac{\ii\left(a_{j} \tilde{q}_{j}+q_{j}\right)\left(\tilde{p}_{j}-p_{j}\right)}{\left(1+a_{j}\right) \varepsilon}+\frac{\ii\left(p_{j} q_{j}-\tilde{p}_{j} \tilde{q}_{j}\right)}{\varepsilon}\right)\right],
\label{modA}
\end{eqnarray}
and
\begin{equation}
\left|A(0,q,p)\right| =2^{m} (\pi\varepsilon)^{\frac{m}{4}}\prod_{j=1}^{m}\left(\frac{\sqrt{a_j}}{1+a_{j}}\right)^{\frac{1}{2}} \exp \left(-\sum_{j=1}^m\frac{\left(\tilde{p}_{j}-p_{j}\right)^{2}+a_{j}\left(\tilde{q}_{j}-q_{j}\right)^{2}}{2\left(1+a_{j}\right) \varepsilon}\right).
\label{modA0}
\end{equation}

As we have stated, the probability density function \(\pi(z_0)\) given by Equation \eqref{piabsA0} is a feasible choice of initial sampling. For Gaussian initial condition \eqref{InitialGaussian}, the resulting probability distribution for the initial position and momentum \(z_0=(q,p)\) is a multivariate normal distribution centered at \(q=\tilde{q}\) and \(p=\tilde{p}\) with probability density function:
\beq\label{GauGau}
\pi(q,p)= (2\pi\varepsilon)^{-m}\prod_{j=1}^{m}\left(\frac{\sqrt{a_j}}{1+a_{j}}\right) \exp \left(-\sum_{j=1}^m\frac{\left(\tilde{p}_{j}-p_{j}\right)^{2}+a_{j}\left(\tilde{q}_{j}-q_{j}\right)^{2}}{2\left(1+a_{j}\right) \varepsilon}\right).
\eeq
Algorithmically, we can directly sample \(z_0=(q,p)\) from Gaussian density function in Equation \eqref{GauGau}, then the initial amplitude \(A(0,q,p)\) can be computed analytically.

\subsection{Main result of sampling error}

In this subsection, we aim to establish a rigorous error estimate for the FGS algorithm. We justify that the Monte Carlo sampling error \(E_S=\left\|\ufgs-\ufga\right\|_{L^2}\) is independent of the scaling parameter \(\ep\) and study the relationship between \(E_S\) and the dimension number \(m\). With the error estimate, we are able to discuss the proper sample size for the FGS for given \(\ep\) and \(m\) and illustrates the advantages of the FGS algorithm, especially in high-dimensional cases.

Recall that \(\ufga\) denotes the FGA ansatz defined in Equation \eqref{FGAansatz} and \(\ufgs\) denotes the FGS wave function given in Equation \eqref{Probz0j}. Now that \(\ufgs\) is determined by a series of random variables \(\left\{z_0^{(j)}\right\}\) on the phase space, we define the quadratic expectation of sampling error \(E_S\) and total error \(E_0\) as follows:
\begin{eqnarray}
\label{CES}\CE_S(t,M)&=&\bE_{z_0^{(j)}\sim\pi}\left\|\ufga(t,\cdot)-\ufgs\kl t,\cdot;M,\left\{z_0^{(j)}\right\}_{j=1}^M\kr\right\|_{L^2}^2,\\
\label{CE0}\CE_0(t,M)&=&\bE_{z_0^{(j)}\sim\pi}\left\|\ufgs\kl t,\cdot;M,\left\{z_0^{(j)}\right\}_{j=1}^M\kr-u(t,\cdot)\right\|_{L^2}^2.
\end{eqnarray}
Then plug Equation \eqref{Probz0j} into the expression of \(\CE_S\):
\begin{eqnarray}\label{eq:norm2est}
\CE_S(t,M)&=&\int_{\mathbb{R}^m}\left[\bE_{z_0^{(j)}\sim\pi}\left|\ufgs(t,x)-\ufga(t,x)\right|^2\right]\dd x\nonumber\\
&=&\int_{\mathbb{R}^m}\left[\bE_{z_0^{(j)}\sim\pi}\left|\frac{1}{M}\sum_{j=1}^M\left(\Lambda\left(t,x,z_0^{(j)}\right)-\ufga(t,x)\right)\right|^2\right]\dd x\nonumber\\
&=&\frac{1}{M}\int_{\mathbb{R}^m}\left(\bE_{z_0\sim\pi}\left|\Lambda\kl t,x,z_0\kr-\ufga(t,x)\right|^2\right)\dd x=:\frac{I_\fga}{M},
\end{eqnarray}
where
\beq\label{IFGA}
I_\fga=\int_{\mathbb{R}^m}\left(\bE_{z_0\sim\pi}\left|\Lambda\kl t,x,z_0\kr-\ufga(t,x)\right|^2\right)\dd x,
\eeq
is regarded as the sampling variance of the FGS wave function. The derivation of Equation \eqref{eq:norm2est} uses the fact that \(z_0^{(j)}\) are independent identically distributed with density \(\pi\). We aim to estimate an upper bound for the integral \(I_\fga\).

Before we present the main result of error estimate, we introduce the definition of the subquadratic condition. We say the potential function \(E(x)\) satisfies the subquadratic condition if there exists a constant \(C_E\) such that
\beq\label{SubQ}
\sup _{q \in \mathbb{R}^{m}}\left|\partial_{\alpha} E_{k}(q)\right| \leq C_{E}, \quad \forall|\alpha|=2.
\eeq
The subquadratic condition is commonly used in analysis of hyperbolic systems \cite{SR2009,LY2011-2,LY2011-3,LZ2017} and we also assume the subquadratic condition in our theorem. With all of the preparation, let us show the main result of this section:
\bthm\label{ThmIFGA}
Consider the semiclassical Schr\"odinger equation \eqref{SemiSchrodinger} up to time \(t\) with the normalized Gaussian initial condition \eqref{InitialGaussian}. Suppose that the potential function \(E(x)\in C^\infty(\BR^m)\) and satisfies the subquadratic condition \eqref{SubQ}. When applying the frozen Gaussian sampling (FGS) as in Equation \eqref{Probz0j} with the initial sampling density \(\pi\) as in Equation \eqref{piabsA0}, there exists a constant \(C_0(m,t)\) independent of the scaling parameter \(\ep\) and the initial data \(u_\ti\) such that the sampling variance given in Equation \eqref{IFGA} satisfies 
\beq\label{EqIFGA0}
I_\fga\leq C_0(t,m)\prod_{j=1}^m\left(\frac{1+a_j}{\sqrt{a_j}}\right)^{\frac{1}{2}},\quad\forall\ep\in(0,1].
\eeq
Moreover, if the assumptions in Lemma \ref{LmA} are valid, the constant \(C_0(t,m)\) in Equation \eqref{EqIFGA0} follows
\beq\label{EqIFGA2}
C_0(t,m)\leq \kl\gamma_0(t)\kr^m,
\eeq
where \(\gamma_0(t)>1\) is a constant independent of \(\ep\) and \(m\).
\ethm

We leave the proof of Theorem \ref{ThmIFGA} to Section 3.3 and discuss on the theorem first. As a corollary, the following estimate on the quadratic expectation \(\CE_S\) and \(\CE_0\) defined in Equation \eqref{CES} and \eqref{CE0} is obvious:

\begin{corollary}\label{CorFGA}
Under the same assumptions in Theorem \ref{ThmIFGA}, we have the following estimates to \(\CE_S\) and \(\CE_0\):
\beq\label{CES}
\CE_S\leq\frac{\kl(\gamma_0(t)\kr^m}{M}\prod_{j=1}^m\left(\frac{1+a_j}{\sqrt{a_j}}\right)^{\frac{1}{2}}.
\eeq
and
\beq
\CE_0\leq\kl C_A(t,m)\ep\right)^2+\frac{\kl(\gamma_0(t)\kr^m}{M}\prod_{j=1}^m\left(\frac{1+a_j}{\sqrt{a_j}}\right)^{\frac{1}{2}}.
\eeq
\end{corollary}

To conclude this subsection, let us demonstrate the advantages of the FGS over the time-splitting spectral methods (see \cite{BJM2003}) through the error estimate in Corollary \ref{CorFGA}. For simplicity, we assume \(a_j=1\) for \(j=1,2,\cdots, m\) and ignore the asymptotic error \(E_\fga=\left\|\ufga-\ufgs\right\|_{L^2}\) of the FGA ansatz, then the \(L^2\) error \(E_0=\left\|u-\ufgs\right\|_{L^2}\)  can be approximated by the quadratic expectation \(\CE_S\), i.e.
\beq
E_0=\left\|u-\ufgs\right\|\approx \kl\CE_S\kr^{\frac{1}{2}}\leq \sqrt{\frac{2^{\frac{m}{2}}\kl(\gamma_0(t)\kr^m}{M}}=\frac{\kl\gamma_1\kr^m}{\sqrt{M}}.
\eeq
where \(\gamma_1=2^{\frac{1}{4}}\kl\gamma_0(t)\kr^{\frac{1}{2}}\) is another constant independent of \(\ep\) and \(m\). In order to guarantee that \(E_0\leq\delta\), the sample size of the FGS should be more than
\beq
\CN_\fgs=\frac{\kl\gamma_1\kr^{2m}}{\delta^2},\quad \delta>0.
\eeq
In contrast, if we use a mesh-based algorithm such as time-splitting spectral methods in \cite{BJM2003} to numerically approximate the semiclassical Schr\"odinger equation \eqref{SemiSchrodinger}, the spatial step length of the mesh should be at least of order \(O(\ep)\) to capture the highly oscillatory wave function. For an \(m-\)dimensional semiclassical Schr\"odinger equation, the spatial degree of freedom of the mesh-based method should be at least
\beq
\CN_{\textrm{TSP}}=\CO\kl \ep^{-m}\kr.
\eeq

Recall that the sample size \(\CN_\fga\) is the \(m-\)th power of an \(\CO(1)\) parameter. Therefore, when scaling parameter \(\ep\ll 1\) or the dimension number \(m\) is large, the computational cost to solve a given semiclassical Schr\"odinger equation with the FGS is much smaller than those with the grid based methods, i.e. \(\CN_\fga\ll \CN_{\textrm{TSP}}\).

\brm
We remark that though the sampling size \(\CN_\fgs\) increases exponentially as \(m\) increases, we examine from numerical examples that the increasing rate of \(\CN_\fgs\) is much smaller compared to \(\ep^{-1}\). For more information of this issue, see Example 5 in Section 5.3.
\erm

\subsection{Proof of Theorem \ref{ThmIFGA}}

Before we prove Theorem \ref{ThmIFGA}, let us first collect some properties of the Hamiltonian flow on the phase space. In this subsection, we assume the potential function \(E(x)\in C^\infty(\BR^m)\) and satisfies the subquadratic condition  \eqref{SubQ} without special instructions. 

For simplicity, we denote the map on the phase space from initial time 0 to time \(t\) as:
\beq
\begin{aligned}
\kappa_{t}: & \mathbb{R}^{2 m} \rightarrow \mathbb{R}^{2 m} \\
&z_0=(q, p) \longmapsto z_t=\left(Q(t,q, p), P(t,q, p)\right)=:\left(Q^{\kappa_t}(q, p), P^{\kappa_t}(q, p)\right).
\end{aligned}
\eeq
For the transformation \(\kappa\) on the phase space, we denote its Jacobian matrix as
\beq\label{Jac}
J^{\kappa_t}(q, p)=\left(\begin{array}{ll}
\left(\partial_{q} Q^{\kappa_t}\right)^{T}(q, p) & \left(\partial_{p} Q^{\kappa_t}\right)^{T}(q, p) \\
\left(\partial_{q} P^{\kappa_t}\right)^{T}(q, p) & \left(\partial_{p} P^{\kappa_t}\right)^{T}(q, p)
\end{array}\right).
\eeq
Note that the transformation \(\kappa\) is canonical, i.e. its Jacobian matrix \(J^{\kappa_t}\) satisfies
\beq
\left(J^{\kappa_t}\right)^{T}\left(\begin{array}{cc}
0 & I_{m} \\
-I_{m} & 0
\end{array}\right) J^{\kappa_t}=\left(\begin{array}{cc}
0 & I_{m} \\
-I_{m} & 0
\end{array}\right),
\eeq
for any \((q,p)\in\BR^{2m}\). Here \(I_m\) denotes \(m\times m\) identity matrix.

The following propositions give uniform bounds of the Jacobian matrix \(J^{\kappa_t}(q, p)\) and its derivatives.
\begin{proposition}\label{ProQqPp}
For given time \(t>0\) and \(k\in\BN\), there exists a constant \(K_0(t,k)\) such that
\beq
\sup_{(q,p)\in\BR^{2m}}\max_{\left|\alpha_{p}\right|+\left|\alpha_{q}\right| \leqslant k}\left|\partial_{q}^{\alpha_{q}} \partial_{p}^{\alpha_{p}}\left[J^{\kappa_{t}}(q, p)\right]\right| \leqslant K_0(t,k),
\eeq
where \(|\cdot|\) denotes the matrix norm.
\end{proposition}

Note that Proposition \ref{ProQqPp} was proven in Proposition 1 of \cite{SR2009}.

For the canonical transformation \(\kappa\), we define
\beq
Z^{\kappa_t}(q, p)=\partial_{z}\left(Q^{\kappa_t}(q, p)+\ii P^{\kappa_t}(q, p)\right),
\eeq
where \(\partial_z=\partial_q-\ii\partial_p\). The invertibility and boundness estimate of \(Z^{\kappa_t}\) can be proven by mimicking the proof of Lemma 5.1 in \cite{LY2011-3}. We summarize them into the following proposition.
\begin{proposition}\label{ProZ}
For given time \(t>0\) and \(k\in\BN\), there exists a constant \(K_0(t,k)\) such that
\beq
\sup_{(q,p)\in\BR^{2m}}\max _{\left|\alpha_{p}\right| +\left|\alpha_{q}\right| \leqslant k}\left|\partial_{q}^{\alpha_{q}} \partial_{p}^{\alpha_{p}}\left[\left(Z^{\kappa_{t}}(q, p)\right)^{-1}\right]\right| \leqslant K_0(t,k),
\eeq
where \(|\cdot|\) denotes the matrix norm.
\end{proposition}

With the properties of the matrix \(J^{\kappa_t}\) and  \(\left(Z^{\kappa_{t}}\right)^{-1}\), we finally establish the estimate of the FGA variable \(A(t,q,p)\):
\begin{lemma}\label{LmA}
Consider the FGA variable \(A(t,q,p)\) that solves Equation \eqref{PhaseA} on the phase space. For given time \(t>0\), there exists a constant \(K_1(t,m)\) such that
\beq\label{FGAAbound}
\sup_{(q,p)\in\BR^{2m}}\left|\frac{A(t,q,p)}{A(0,q,p)}\right|\leq K_1(t,m).
\eeq
Moreover, if we assume the constant \(K_0(t,1)\) in Proposition \ref{ProQqPp} and \ref{ProZ} is independent of the dimension number \(m\), there exists a parameter \(\gamma^*(t)>1\) such that 
\beq\label{LemBound}
K_1(t,m)\leq \kl \gamma^*(t)\kr^m.
\eeq
\end{lemma}

The proof is straightforward through Gronwall's inequality.

\bprf
Based on Equation \eqref{PhaseA}, we arrive at
\beq\label{PfLm1}
\frac{\dd}{\dd t}\left|A(t)\right|\leq \frac{\left|A(t)\right|}{2}\left|\text{tr}\left[Z^{-1}\left(\partial_zP-\ii\partial_zQ\nabla_Q^2E(Q)\right)\right]\right|.
\eeq
Note that the constants in Proposition \ref{ProQqPp} and \ref{ProZ} are independent of the dimension number \(m\), we use the constant \(K^*(t)=K_0(t,1)\) to estimate Equation \eqref{PfLm1}:
\beq
\frac{1}{2}\left|\text{tr}\left[Z^{-1}\left(\partial_zP-\ii\partial_zQ\nabla_Q^2E(Q)\right)\right]\right|\leq 2m \left(K^*(t)\right)^2.
\eeq
Then apply Gronwall's inequality to Equation \eqref{PfLm1}:
\beq
|A(t)|\leq \ee^{2m \left(K^*(t)\right)^2}=\kl\gamma(t)\kr^m,
\eeq
where \(\gamma^*(t)=\ee^{2\left(K^*(t)\right)^2}\) is independent of \(\ep\) and \(m\).
\eprf

\brm
We remark that the importance of the three propositions and lemmas lie in the fact that the estimate to the FGA trajectories and FGA variables is uniform with respect to the initial phase \(z_0=(q,p)\). Hence we can make global estimate on the integration \(I_\fga\). Validity of the propositions and lemmas rely on the subquadratic condition of the potential function \(E(x)\). 
\erm

After the preparation above, we finally move to the proof of Theorem \ref{ThmIFGA}:

\bprf
According to Equation \eqref{IFGA}, we have
\beq
I_\fga=\int_{\mathbb{R}^m}\left(\bE_{z_0\sim\pi}\left|\Lambda\kl t,x,z_0\kr-\ufga(t,x)\right|^2\right)\dd x\leq \int_{\mathbb{R}^m}\left(\bE_{z_0\sim\pi}\left|\Lambda\kl t,x,z_0\kr\right|^2\right)\dd x,
\eeq
and it suffices that we focus on the following secondary moment:
\beq\label{M2E}
M_{(2)}(x)=\bE_{z_0\sim \pi}\left|\Lambda\kl t,x,z_0\kr\right|^2.
\eeq
Then we plug Equation \eqref{Lambda} into Equation \eqref{M2E}
\begin{eqnarray}\label{varstart}
M_{(2)}(x)
&=&\frac{1}{(2\pi\ep)^{3m}}\int_{\BR^{2m}}\pi(z_0)\left|\frac{A(t,z_0)}{\pi(z_0)}\right|^2\exp\left(-\frac{\left|x-Q(t,z_0)\right|^2}{\varepsilon}\right)\dd z_0\nonumber\\
&=&\frac{1}{(2 \pi \varepsilon)^{3 m}}\int_{\mathbb{R}^{2m}}\frac{\left|A(t,z_0)\right|^2}{\pi\left(z_0\right)}\exp\kl-\frac{\left|x-Q(t,z_0)\right|^2}{\varepsilon}\kr\dd z_0.
\end{eqnarray}
Now recall that the sampling density \(\pi\) is given by Equation \eqref{piabsA0}, then we can simplify Equation \eqref{varstart}:
\beq \label{varest1}
M_{(2)}(x)=\frac{\int_{\BR^{2m}}|A(0,z_0)|\dd z_0}{(2 \pi \varepsilon)^{3 m}}\int_{\mathbb{R}^{2m}}\frac{\left|A(t,z_0)\right|^2}{|A(0,z_0)|}\exp\kl-\frac{\left|x-Q(t,z_0)\right|^2}{\varepsilon}\kr\dd z_0.
\eeq

Next we integrate \(M_{(2)}(x)\) in Equation \eqref{varest1} with respect to the position variable \(x\):
\begin{eqnarray}\label{varL1est}
&&\int_{\mathbb{R}^m}M_{(2)}(x)\dd x\nonumber\\
&=&\frac{\int_{\BR^{2m}}|A(0,z_0)|\dd z_0}{(2 \pi \varepsilon)^{3 m}}\int_{\mathbb{R}^m}\left[\int_{\mathbb{R}^{2m}}\frac{\left|A(t,z_0)\right|^2}{|A(0,z_0)|}\exp\kl-\frac{\left|x-Q(t,z_0)\right|^2}{\varepsilon}\kr \dd z_0\right]\dd x\nonumber\\
&=&\frac{\int_{\BR^{2m}}|A(0,z_0)|\dd z_0}{(2 \pi \varepsilon)^{3 m}}\int_{\mathbb{R}^{2m}}\frac{\left|A(t,z_0)\right|^2}{|A(0,z_0)|}\left(\int_{\mathbb{R}^m}\exp\kl-\frac{\left|x-Q(t,z_0)\right|^2}{\varepsilon}\kr\dd x\right)\dd z_0\nonumber\\
&=&\frac{\int_{\BR^{2m}}|A(0,z_0)|\dd z_0}{(2 \pi \varepsilon)^{\frac{5m}{2}}}\int_{\mathbb{R}^{2m}}\frac{\left|A(t,z_0)\right|^2}{|A(0,z_0)|}\dd z_0.
\end{eqnarray}

Now we apply Lemma \ref{LmA} in Equation \eqref{varL1est}:
\begin{eqnarray}\label{EqIFGA}
I_\fga&\leq &\int_{\mathbb{R}^m}M_{(2)}(x)\dd x\nonumber\\
&\leq&\frac{\int_{\BR^{2m}}|A(0,z_0)|\dd z_0}{(2 \pi \varepsilon)^{\frac{5m}{2}}}\int_{\mathbb{R}^{2m}}\left(K(t,m)\right)^2 |A(0,z_0)|\dd z_0\nonumber\\
&=&\left(K_1(t,m)\right)^2\frac{\left(\int_{\BR^{2m}}|A(0,z_0)|\dd z_0\right)^2}{(2 \pi \varepsilon)^{\frac{5m}{2}}}.
\end{eqnarray}
Here \(K_1(t,m)\) is a constant independent to \(\ep\).

According to Equation \eqref{modA0}, we can calculate the integration of \(|A(0,z_0)|\) on the phase space:
\beq\label{IA0mod}
\int_{\BR^{2m}}|A(0,z_0)|\dd z_0=2^{2m}(\pi\ep)^{\frac{5m}{4}}\prod_{j=1}^{m}\left(\frac{1+a_{j}}{\sqrt{a_{j}}}\right)^{\frac{1}{2}}.
\eeq
Finally we use Equation \eqref{IA0mod} in Equation \eqref{EqIFGA}:
\beq\label{IFGAConstant}
I_\fga\leq 2^{\frac{3m}{2}}(K_1(t,m))^2\prod_{j=1}^m\left(\frac{1+a_j}{\sqrt{a_j}}\right)^{\frac{1}{2}}.
\eeq
Hence the right hand side of Equation \eqref{IFGAConstant} is a constant independent to \(\ep\). Let \(C_0(t,m)=2^{\frac{3m}{2}}\kl K_1(t,m) \kr^2\), we finally arrive at Equation \eqref{EqIFGA0}.

Through the expression of \(C_0(t,m)\), it is obvious that Equation \eqref{EqIFGA2} is valid provided that Equation \eqref{LemBound} in Lemma \ref{LmA} is satisfied.
\eprf

\section{Frozen Gaussian sampling with WKB initial conditions}

In this section, we study the semiclassical Schr\"odinger equation \eqref{SemiSchrodinger} with the following normalized WKB initial condition:
\beq\label{WKBInitial}
u_\ti(x)=a_\ti(x)\exp\left(\frac{\ii}{\varepsilon}S_\ti(x)\right),
\eeq
where \(a_\ti\) and \(S_\ti\) are two real-valued functions on \(\BR^m\), which are both independent of \(\ep\). The initial amplitude \(a_\ti\) follows
\beq
\int_{\BR^m}a_\ti^2(x)\dd x=1.
\eeq
With the celebrated WKB approximation (see Appendix B for a brief introduction), the exact solution to Equation \eqref{SemiSchrodinger} has an asymptotic description of WKB type before caustics onset, provided that \(a_\ti\) is smooth and decays rapidly \cite{JMS2011}. Practically, we mainly consider the cases with the following Gaussian initial amplitude:
\begin{equation}\label{a0WKB}
a_{\text{in}}(x)=\left(\prod_{j=1}^m a_j\right)^{\frac{1}{4}}\pi^{-\frac{m}{4}}\exp\left(-\frac{1}{2}\sum_{j=1}^m a_j(x_j-\tilde{x}_j)^2\right),
\end{equation}
where \(a_j>0\) are sufficiently big and \(\tilde{x}_j\in\BR\) (\(j=1,2,\cdot,m\)) are given parameters. We denote the mean of \(a_\ti\) as 
\beq\label{deftx}
\tilde{x}=\kl x_1,x_2,\cdots,x_n\kr^T.
\eeq
The cases with general initial amplitude \(a_\ti\) may be treated in a similar fashion and we refer to Remark \ref{no4.2} for discussion on this issue.

Compared to Gaussian wave packet initial data cases where the explicit expression of \(A(0,q,p)\) is available, the difficulty of applying the FGS in WKB initial data cases lies in that there is no natural choice for the reference density function \(\pi(z_0)\) which achieves the balance between variance reduction and convenience of sampling. Moreover, the WKB approximation to Equation \eqref{SemiSchrodinger} may break down after caustics onset, which brings additional difficulties to the numerical approximation of Equation \eqref{SemiSchrodinger}. In this section, we put forward a new initial sampling algorithm based on the stationary phase approximation to the initial amplitude \eqref{A0Compute}. Then we validate the numerical performance of the FGS associating with the new initial sampling strategy through various typical examples in Section 5.

\subsection{Description of the algorithm}

As we have stated, one of the difficulties of the FGS algorithm with a WKB initial condition is that \(A(0,q,p)\) (given in Equation \eqref{A0Compute}) can not be calculated explicitly, hence there is no obvious way to do initial sampling. In this section, we apply the stationary phase method to approximate the initial amplitude \(A(0,q,p)\) in Equation \eqref{A0Compute} and give a practical sampling method based on the analytical approximation to \(A(0,q,p)\). Before we start, let us summarize the stationary phase method into the following lemma:
\begin{lemma}\label{ThmStatPhase}
Let \(\alpha(x)\) and \(\Psi(x)\) be two smooth functions on \(\BR^m\). Suppose that \(\alpha(x)\in C_0^\infty(\BR^m)\) and \(\Psi(x)\in C^\infty(\BR^m)\). Assume that \(\nabla\Psi(x)\) has a finite number of zero points \(\{y_k\}_{k=1}^N\) on the support of \(\alpha\) and \(\det\kl\nabla^2\Psi(y_k)\kr\neq 0\) for each zero point \(y_k\). We consider the highly-oscillatory integration \beq
I_\ep=\int_{\BR^m}\alpha(x)\exp\kl\frac{\ii\Psi(x)}{\ep}\kr\dd x.
\eeq
As \(\ep\to0\), \(I_\ep\) has the following asymptotic behavior:
\beq\label{ThmEqSP}
I_\varepsilon=(2 \pi \varepsilon)^{\frac{m}{2}} \sum_{k=1}^N \frac{\exp\kl\frac{\ii\Psi(x)}{\ep}+\frac{\ii\pi}{4}\mathrm{sgn}\left(\nabla^2 \Psi\left(y_{k}\right)\right)\kr}{\left|\det \nabla^{2} \Psi\left(y_{k}\right)\right|^{\frac{1}{2}}}\left(\alpha\left(y_{k}\right)+\CO(\varepsilon)\right),
\eeq
where the term \(\mathrm{sgn}\left(\nabla^2 \Psi\left(y_{k}\right)\right)\) denotes the number of positive eigenvalues of the matrix \(\nabla^2 \Psi\left(y_{k}\right)\) minus the number of negative eigenvalues.
\end{lemma}

A proof to Lemma \ref{ThmStatPhase} can be found in \cite{PDE}.

According to Equation \eqref{ThmEqSP}, the initial amplitude \(A(0,q,p)\) in Equation \eqref{A0Compute} can be calculated asymptotically as follows:
\begin{eqnarray}
\label{SPhaseA0}
A(0,q,p)&=&2^{\frac{m}{2}}\int_{\mathbb{R}^m}a_{\text{in}}(y)\ee^{\frac{i}{\varepsilon}S_{\text{in}}(y)}e^{\frac{\ii}{\varepsilon}(-p\cdot(y-q)+\frac{\ii}{2}|y-q|^{2})} \dd y\nonumber\\
&=&2^{\frac{m}{2}}\int_{\mathbb{R}^m}\left(a_{\text{in}}(y)e^{-\frac{|y-q|^2}{2\varepsilon}}\right)\ee^{\frac{\ii}{\varepsilon}(S_{\text{in}}(y)-p\cdot(y-q))} \dd y\nonumber\\
&\approx&2^m(\pi\varepsilon)^{\frac{m}{2}}\sum_{\nabla S_{\text{in}}(y_k)=p}\frac{\exp\left(\frac{\ii}{\varepsilon}(S_{\ti}(y_k)-p\cdot(y_k-q))+\frac{\ii \pi}{4}\operatorname{sgn}\left(\nabla^{2} \phi\left(y_{k}\right)\right)\right)}{\left|\mathrm{det} \nabla^{2}S_{\ti}\left(y_{k}\right)\right|^{\frac{1}{2}}}\nonumber\\&&\qquad\qquad\qquad\qquad\qquad\qquad\qquad\qquad\times\left(a_{\ti}(y_k)\ee^{-\frac{|y_k-q|^2}{2\varepsilon}}+\CO(\ep)\right).
\end{eqnarray}

In order to illustrate the spirit of sampling strategy, we make some reasonable assumptions which simplify the asymptotic description of \(A(0,q,p)\) in Equation \eqref{SPhaseA0}. The cases when the assumptions break down will be discussed in Remark 4.3. We will show the 

\begin{assumption}
For any \(p\in\BR^m\), there exists at most one \(y\in\BR^m\) such that \(\nabla S_{\ti}(y)=p\), i.e. \(\nabla S_\ti\) is an injective mapping. The inverse mapping of \(\nabla S_\ti\) is given by \(T:p\mapsto y\). The definition domain of \(T\) may not be \(\BR^n\) and we denote it by \(\CD(T)\), i.e.
\beq\label{deft}
T:\CD(T)\to\BR^m,\quad p\mapsto y.
\eeq
\end{assumption}

Recall that we assume the amplitude of the WKB initial data \(a_\ti(x)\) decays exponentially as \(x\to\infty\), the support of \(a_\ti\) can be viewed as a small interval. Hence it is reasonable to assume that for any \(p\in\BR^m\), there exists at most one \(y\in\BR^m\) such that \(\nabla S_{\text{in}}(y)=p\) and \(a_\ti(y)\neq 0\). With this assumption, the summation in Equation \eqref{SPhaseA0} can be reduced to only one term. Besides, we will show later in Remark 4.3 that the difficulties caused by multivalueness of the mapping \(T\) are not essential. 

For simplicity, we use the variable \(y\) to denote the unique image of the mapping \(T\) defined in Equation \eqref{deft} at \(p\in\CD(T)\) in the rest of this section without special instructions. 

The next assumption is a technique assumption without lose of generality:
\begin{assumption}
The inverse mapping \(T\) defined in Equation \eqref{deft} is smooth and satisfies
\beq
\mathrm{det}\kl\nabla^{2}S_{\ti}(y)\kr\neq0,\quad\forall p\in\CD(T).
\eeq
\end{assumption}

Based on the two assumptions, the asymptotic approximation \eqref{SPhaseA0} of the initial amplitude \(A(0,q,p)\) can be simplified as follows:
\begin{eqnarray}\label{WKBA0sp}
A_{\rsp}(q,p)&=&2^m(\pi\varepsilon)^{\frac{m}{2}}\chi_{\CD(T)}(p)a_\text{in}(y)\exp\left(-\frac{|y-q|^2}{2\varepsilon}\right)\nonumber\\&&\times\frac{\exp\left(\frac{\ii}{\varepsilon}(S_\text{in}(y)-p\cdot(y-q))\right)+\frac{\ii \pi}{4}\operatorname{sgn}\left(\nabla^{2} \phi\left(y\right)\right)}{\left|\operatorname{det} \nabla^{2}S_\text{in}\left(y\right)\right|^{\frac{1}{2}}},
\end{eqnarray}
where \(\chi_{\CD(T)}(p)\) stands for the indicative function on \(\CD(T)\). Then the modulus of Equation \eqref{WKBA0sp} is given by
\begin{equation}\label{WKBansA0}
|A_{\rsp}(q,p)|=2^m(\pi\varepsilon)^{\frac{m}{2}}\chi_{\CD(T)}(p)a_\text{in}(y)\frac{\exp\left(-\frac{|y-q|^2}{2\varepsilon}\right)}{\left|\operatorname{det} \nabla^{2}S_\text{in}\left(y\right)\right|^{\frac{1}{2}}}.
\end{equation}

For Gaussian initial condition cases, we directly use the real-valued Gaussian function \(|A(0,q,p)|\) on the phase space as the sampling density \(\pi\) in Equation \eqref{ProbInterrupt}, which results in the minimum sampling variance. However, for the WKB initial condition cases, it is impossible to sample towards \(|A_{\rsp}(q,p)|\) given by Equation \eqref{WKBansA0} directly since the non-Gaussian function \(|A_{\rsp}(q,p)|\) is hard to sample. To balance between variance reduction and convenience of sampling, we choose the following sampling density function:
\begin{equation}\label{densityWKBsp}
\pi(q,p)=\frac{1}{\CZ}2^m(\pi\varepsilon)^{\frac{m}{2}}\chi_{\CD(T)}(p)a_\text{in}(y)\frac{\exp\left(-\frac{|y-q|^2}{2\varepsilon}\right)}{\left|\operatorname{det} \nabla^{2}S_\text{in}\left(y\right)\right|},
\end{equation}
where \(\CZ\) is a normalization parameter to ensure \(\pi\) a probability density, i.e. 
\[
\int_{\BR^{2m}}\pi(q,p)\dd q\dd p=1.
\]

We note that the only difference between \(|A_{\rsp}(q,p)|\) in Equation \eqref{WKBansA0} and the probability density \(\pi\) in Equation \eqref{densityWKBsp} lies in the exponent of the term \(\left|\det \nabla^{2}S_\ti\left(y\right)\right|\) in the denominator. Hence the choice of \(\pi\) in Equation \eqref{densityWKBsp} is rather acceptable from the prospective of variance reduction. The following proposition shows that the probability density \(\pi\) in Equation \eqref{densityWKBsp} is also easy to sample, provided that the the amplitude \(a_\ti\) is given by Equation \eqref{a0WKB}:

\begin{proposition}\label{WKBpiPro}
Consider the WKB initial data \eqref{WKBInitial} where the amplitude function \(a_\ti\) of is given by Equation \eqref{a0WKB}. Assume that Assumption 1 and Assumption 2 are valid. Suppose that the random variables \(\CQ\) and \(\CP\) are sampled with respect to the probability density \(\pi(q,p)\) defined in Equation \eqref{densityWKBsp}, then \(\CQ\) and \(\CP\) follows 
\begin{enumerate}[itemindent=0em]
    \item The marginal distribution of the random variable \(\CQ\) is a normal distribution with expectation \(\BE_\pi q=\tilde{x}\), where \(\tilde{x}\) as in Equation \eqref{deftx} denotes the mean of \(a_\ti\).
    \item Denote \(\CY=T(\CP)\), then the conditional distribution of the random variable \(\CY\) with respect to \(\CQ\) is a normal distribution.
    \item The normalization parameter \(\CZ\) in Equation \eqref{densityWKBsp} is given as follows:
    \begin{equation}\label{ZWKB}
    \CZ=2^{2m}\pi^{\frac{5m}{4}}\varepsilon^{\frac{m}{2}}\left(\prod_{j=1}^m \frac{\sqrt{a_j}\left(\frac{1}{a_j}+\varepsilon\right)}{a_j+\frac{1}{\varepsilon}}\right)^{\frac{1}{2}}.
    \end{equation}
\end{enumerate}
\end{proposition}

\bprf
Recall that the mapping \(T\) as in Equation \eqref{deft} is the inverse mapping of the gradient function \(\nabla S_\ti\). We refer to \(\CJ_T(p)\) as the Jacobian matrix of the mapping \(T\) at \(p\). Based on Assumption 2, the determinant of \(\CJ_T(p)\) follows:
\begin{equation}\label{det-1}
\left|\mathrm{det}\CJ_T(p)\right|=\frac{1}{\left|\mathrm{det} \nabla^{2}S_\text{in}\left(y\right)\right|},\quad\forall p\in\CD(T).
\end{equation}
Then we apply Equation \eqref{det-1} and the integral substitution rule to compute the marginal density \(\pi_\CQ(q)\) of the random variable \(\CQ\):
\begin{eqnarray}\label{piqq}
\pi_\CQ(q)
&=&\frac{1}{\CZ}\left(\prod_{j=1}^m a_j\right)^{\frac{1}{4}}2^m\pi^{\frac{m}{4}}\varepsilon^{\frac{m}{2}}\int_{\CD(T)}\frac{\exp\left[-\sum_{j=1}^m\frac{a_j}{2}(y_j-\tilde{x}_j)^2-\frac{1}{2\varepsilon}|y-q|^2\right]}{\left|\mathrm{det} \nabla^{2}S_\ti\left(y\right)\right|}\dd p\nonumber\\
&=&\frac{1}{\CZ}\left(\prod_{j=1}^m a_j\right)^{\frac{1}{4}}2^m\pi^{\frac{m}{4}}\varepsilon^{\frac{m}{2}}\int_{\mathbb{R}^m}\exp\left[-\sum_{j=1}^m\frac{a_j}{2}(y_j-\tilde{x}_j)^2-\frac{1}{2\varepsilon}|y-q|^2\right]\dd y\nonumber\\
&=&\frac{1}{\CZ}\left(\prod_{j=1}^m \frac{\sqrt{a_j}}{a_j+\frac{1}{\varepsilon}}\right)^{\frac{1}{2}}2^{\frac{3m}{2}}\pi^{\frac{3m}{4}}\varepsilon^{\frac{m}{2}}\exp\left[-\sum_{j=1}^m\frac{(q_j-\tilde{x}_j)^2}{2\left(\varepsilon+\frac{1}{a_j}\right)}\right],
\end{eqnarray}
Then the marginal distribution of the random variable \(\CQ\) is a normal distribution \(\CN(\tilde{x},\Sigma_1)\), where the covariance matrix \(\Sigma\) is a diagonal matrix as follows:
\beq\label{SM1}
\Sigma_1=\mathrm{diag}\kl \ep+\frac{1}{a_1},\ep+\frac{1}{a_2},\cdots,\ep+\frac{1}{a_n}\kr.
\eeq
To calculate the expression of the normalization parameter \(\CZ\) given in Equation \eqref{ZWKB}, we just need to integrate \(\pi_\CQ(q)\) by \(q\) to ensure that \(\int\pi_\CQ(q)\dd q=1\), i.e.
\beq
\CZ=\left(\prod_{j=1}^m \frac{\sqrt{a_j}}{a_j+\frac{1}{\varepsilon}}\right)^{\frac{1}{2}}2^{\frac{3m}{2}}\pi^{\frac{3m}{4}}\varepsilon^{\frac{m}{2}}\int_{\BR^m}\exp\left[-\sum_{j=1}^m\frac{(q_j-\tilde{x}_j)^2}{2\left(\varepsilon+\frac{1}{a_j}\right)}\right]\dd q.
\eeq

For fixed \(\CQ=q\), we denote the conditional density of the random variables \(\CP\) and \(\CY\) with respect to \(\CQ\) as \(\pi_{\CP|\CQ}(p)\) and \(\pi_{\CY|\CQ}(y)\). According to Equation \eqref{piqq}, we compute \(\pi_{\CP|\CQ}(p)\) as follows:
\begin{equation}
\pi_{\CP|\CQ}(p)=\frac{\pi(q,p)}{\pi_\CQ(q)}=\frac{2^m(\pi\varepsilon)^{\frac{m}{2}}}{\CZ\pi_\CQ(q)}a_\ti(y)\frac{\exp\left(-\frac{|y-q|^2}{2\varepsilon}\right)}{\left|\mathrm{det} \nabla^{2}S_\ti\left(y\right)\right|}.
\end{equation}
Then we can compute \(\pi_{\CY|\CQ}(y)\):
\begin{eqnarray}\label{piTpq}
\pi_{\CY|\CQ}(y)&=&\pi_{\CP|\CQ}(p)\left|\mathrm{det} \CJ_T(p)\right|^{-1}\nonumber\\
&=&\frac{2^m(\pi\varepsilon)^{\frac{m}{2}}}{\CZ\pi_\CQ(q)}a_\ti(y)\exp\left(-\frac{|y-q|^2}{2\varepsilon}\right)\nonumber\\
&=&\left(\prod_{j=1}^m a_j\right)^{\frac{1}{4}}\frac{2^m\pi^{\frac{m}{4}}\varepsilon^{\frac{m}{2}}}{\CZ\pi_\CQ(q)}\exp\left[-\sum_{j=1}^m\frac{a_j}{2}(y_j-\tilde{x}_j)^2-\frac{1}{2\varepsilon}|y-q|^2\right].\nonumber\\&&
\end{eqnarray}
Note that the exponent term of the conditional density function \(\pi_{\CY|\CQ}(y)\) is quadratic of \(y\), therefore the conditional distribution of \(\CY\) with fixed \(\CQ\) is a normal distribution. For fixed \(\CQ=q\), the conditional distribution of \(\CY\) follows \(\CN\kl \mu_2,\Sigma_2\kr\), where the mean vector and covariance matrix is given by
\beq\label{mu2}
\mu_2=\kl\frac{\ep a_1\tilde{x}_1+q_1}{\ep a_1+1},\frac{\ep a_2\tilde{x}_2+q_2}{\ep a_2+1},\cdots,\frac{\ep a_n\tilde{x}_n+q_n}{\ep a_n+1}\kr^T,
\eeq
and
\beq\label{SM2}
\Sigma_2=\mathrm{diag}\kl\frac{\ep}{\ep a_1+1},\frac{\ep}{\ep a_2+1},\cdots,\frac{\ep}{\ep a_n+1}\kr.
\eeq
\eprf

Motivated by Proposition \ref{WKBpiPro}, one can sample the random variables \((q,p)\) obeying the probability density \(\pi\) defined in Equation \eqref{densityWKBsp} on the phase space by sampling two normally distributed random variable \(q\) and \(y\), where \(p=\nabla S_\ti(y)\). The initial sampling method bases on Proposition \ref{WKBpiPro}  can be summarized into the following algorithm:
\begin{algorithm}
	\caption{Initial sampling of \((q,p)\) on the phase space for WKB initial data \eqref{WKBInitial} with amplitude \eqref{a0WKB}} \label{alg:WKB}
	\begin{algorithmic}[1]
		\State \textbf{Initial position sampling:} Sample \(q\sim\CN\kl\tilde{x},\Sigma_1\kr\), where \(\tilde{x}\) and \(\Sigma_1\) is given by Equation \eqref{deftx} and \eqref{SM1}.
		\State \textbf{Initial momentum sampling:} For specified \(q\), sample \(y\sim\CN\kl\mu_2,\Sigma_2\kr\) and compute \(p=\nabla S_\ti(y)\). Here \(\mu_2\) and \(\Sigma_2\) are given in Equation \eqref{mu2} and \eqref{SM2}.
	\end{algorithmic}
\end{algorithm}

To conclude this subsection, we make two remarks on the initial sampling algorithm \ref{alg:WKB}. The first remark elaborate what happens when Assumption 1 is not satisfied. The second remark discuss the cases of WKB initial data \eqref{WKBInitial} with general amplitude \(a_\ti\) and generalize Proposition \ref{WKBpiPro} to these cases.

\brm\label{asv1}
Now We discuss the possible treatments when Assumption 1 is violated. If there exists more than one \(y\in\BR^n\) such that \(\nabla S_\ti(y)=p\) for some \(p\in\BR^n\), the asymptotic approximation \(A_\rsp(q,p)\) of the initial amplitude \(A(0,q,p)\)(see Equation \eqref{SPhaseA0}) can't be reduced to only one term as in Equation \eqref{WKBansA0}, i.e.
\ben\label{eqasv1}
A_{\rsp}(q,p)&=&2^m(\pi\varepsilon)^{\frac{m}{2}}\chi_{\CD(T)}(p)\sum_{\nabla S_{\text{in}}\kl y_k\kr=p}\left[ a_\text{in}\kl y_k\kr\exp\left(-\frac{|y-q|^2}{2 \varepsilon}\right)\right.\nonumber\\&&\left.\times\frac{\exp\left(\frac{\ii}{\varepsilon}(S_\text{in}\kl y_k\kr-p\cdot\kl y_k-q \kr)\right)+\frac{\ii \pi}{4} \operatorname{sgn}\left(\nabla^{2} \phi\left(y_{k}\right)\right)}{\left|\operatorname{det} \nabla^{2}S_\text{in}\kl y_k\kr\right|^{\frac{1}{2}}}\right].
\een
The sampling density associating with \(A_{\rsp}(q,p)\) in Equation \eqref{eqasv1} can be given through mimicking the density function \eqref{densityWKBsp}:
\beq
\pi_{\rsp}(q,p)=2^m(\pi\varepsilon)^{\frac{m}{2}}\chi_{\CD(T)}(p)\sum_{\nabla S_{\text{in}}\kl y_k\kr=p}a_\text{in}\kl y_k\kr\exp\left(-\frac{|y-q|^2}{2 \varepsilon}\right).
\eeq
Note that the density function \(\pi_{\rsp}\) determines a mixture distribution and each component of the summation in the density function \(\pi_{\rsp}\) can be sampled through algorithm \ref{alg:WKB}. Then it is natural to adopt the method of sampling the mixture model (see Section 4.2 in \cite{ELV}) to sample \(\pi_{\rsp}\). Therefore, the difficulties resulted from the invalidity of Assumption 1 is not essential.
\erm

\brm\label{no4.2}
We discuss the possible treatments when \(a_\ti\) is not Gaussian. Since \(a_\ti\) is of a general profile, Proposition \ref{WKBpiPro} no longer holds and there is no direct way to sample the probability density \(\pi\) given in Equation \eqref{densityWKBsp}. To overcome these problems, we may approximate \(a_\ti\) with a Gaussian amplitude \(a_\ti^S\) and sample the following probability with Algorithm \ref{alg:WKB}:
\beq
\pi^S(q,p)=\frac{1}{\CZ}2^m(\pi\varepsilon)^{\frac{m}{2}}\chi_{\CD(T)}(p)a^S_\text{in}(y)\frac{\exp\left(-\frac{|y-q|^2}{2\varepsilon}\right)}{\left|\operatorname{det} \nabla^{2}S_\text{in}\left(y\right)\right|},
\eeq
where we have replaced \(a_\ti\) in Equation \eqref{densityWKBsp} with the approximated Gaussian amplitude \(a_\ti^S\). Ideally, the density \(\pi^S\)  balances the need of variance reduction and convenience of sampling such that it serves as the sampling density of the FGS when \(a_\ti\) is not Gaussian. 
\erm

\section{Numerical examples}

In the previous sections, we have elaborated the FGS algorithm with Gaussian initial data and WKB initial data and studied the sampling error with respect to the sampling size \(M\), the scaled Planck parameter \(\ep\) and the dimension number \(m\). In this section, we present five typical examples using the FGS to numerically approximate the semiclassical Schr\"odinger equation \eqref{SemiSchrodinger} to illustrate the theoretical results on the sampling error and manifest the practical potential of the FGS in solving Equation \eqref{SemiSchrodinger}, especially in high-dimensional cases with small \(\ep\). The numerical examples include exploring the relationship between the sampling error and parameter \(\ep\), applying the FGS in challenging cases such as caustics and computing the physical observable with the FGS in high-dimensional cases. The rest of this section is organized as follows. In Section 5.1 and 5.2, we investigate the cases with Gaussian initial data and WKB initial data respectively. In section 5.3, we compute the physical observables for a high-dimensional problem with the FGS and compare the numerical results with the exact value obtained through the Hagedorn wave packet method. 

\subsection{Examples with Gaussian initial conditions }

\subsubsection*{Example 1: 1D Gaussian initial data examples } 

In this example, we aim to justify the theoretical results on the sampling error \(E_S\) in Section 3 through numerical examples. We implement the FGS to approximate Equation \eqref{SemiSchrodinger}  at time \(t=0.5\) and study the relationship among the sampling error \(E_S\), the scaled Planck parameter \(\ep\) and ensemble size \(M\).

In the Gaussian initial condition \eqref{InitialGaussian}, let \(m=1\), \(a_1=2\), \(\tilde{p}=-\frac{1}{2}\) and \(\tilde{q}=\frac{1}{2}\). Then we obtain the following Gaussian wave packet initial data:
\begin{equation}\label{Gaussianini1D}
u_{\ti}(x)=\left(\frac{\pi\varepsilon}{2}\right)^{\frac{1}{4}}\exp\left[-\frac{\ii}{2\varepsilon}\left(x-\frac{1}{2}\right)-\frac{1}{2\varepsilon}\left(x-\frac{1}{2}\right)^2\right].
\end{equation}
Here we consider two different potential functions respectively: the harmonic potential \(E_1(x)=\frac{x^2}{2}\) and the torsion potential \(E_2(x)=1-\cos x\). 

Recall that the \(L^2\) sampling error \(E_S\) is defined as follows:  
\beq\label{ES2}
E_S\left(M,\varepsilon\right)=\left\|\ufgs\left(t,x;M,\left\{z_0^{(j)}\right\}\right)-\ufga(t,x)\right\|_{L^2},
\eeq
where \(\ufga\) denotes the FGA ansatz \eqref{FGAansatz}, which deviates from the exact solution to Equation \eqref{SemiSchrodinger} by the asymptotic error of the FGA ansatz. Note that \(\ufga\) is numerically unavailable, we choose the FGS wave function with a sufficiently large ensemble size \(M_0\) as an accurate approximation of the FGA ansatz
\beq\label{FGScheck}
\ufga(t,x)\approx\ufgs\left(t,x;M_0,\left\{z_0^{(j)}\right\}_{j=1}^{M_0}\right)=\frac{1}{M_0}\sum_{j=1}^{M_0}\Lambda\kl t,x,z_0^{(j)}\kr,
\eeq
where \(M_0=1.5\times 10^5\) is chosen large enough. 

To compute the FGS wave function \(\ufgs\) and \(\ufga\), we use the fourth order Runge-Kutta method to solve the ODE system \eqref{PhaseQ}-\eqref{PhaseA} in the time interval \([0,0.5]\) with time step length \(\Delta t=0.01\) and reconstruct the wave function \(\ufgs\) at \(t=0.5\). We emphasize that \(\Delta t\) is chosen sufficiently small, hence the error due to the approximation of the ODE system is negligible. The wave function reconstruction is carried on the spatial interval \([-\pi,\pi]\). Also note that the FGS wave function \(\ufgs\) with a fixed ensemble size \(M\) may vary with different initial samples \(\left\{z_0^{(j)}\right\}_{j=1}^M\), the sampling error \(E_S\) is actually a random variable. Therefore, when approximating the FGS wave function \(\ufgs\), we average over numerical results obtained through 30 sets of initial samples for each simulation of \(E_S(M,\ep)\).

\begin{figure}
\centering
\includegraphics[width=6cm,height=5cm]{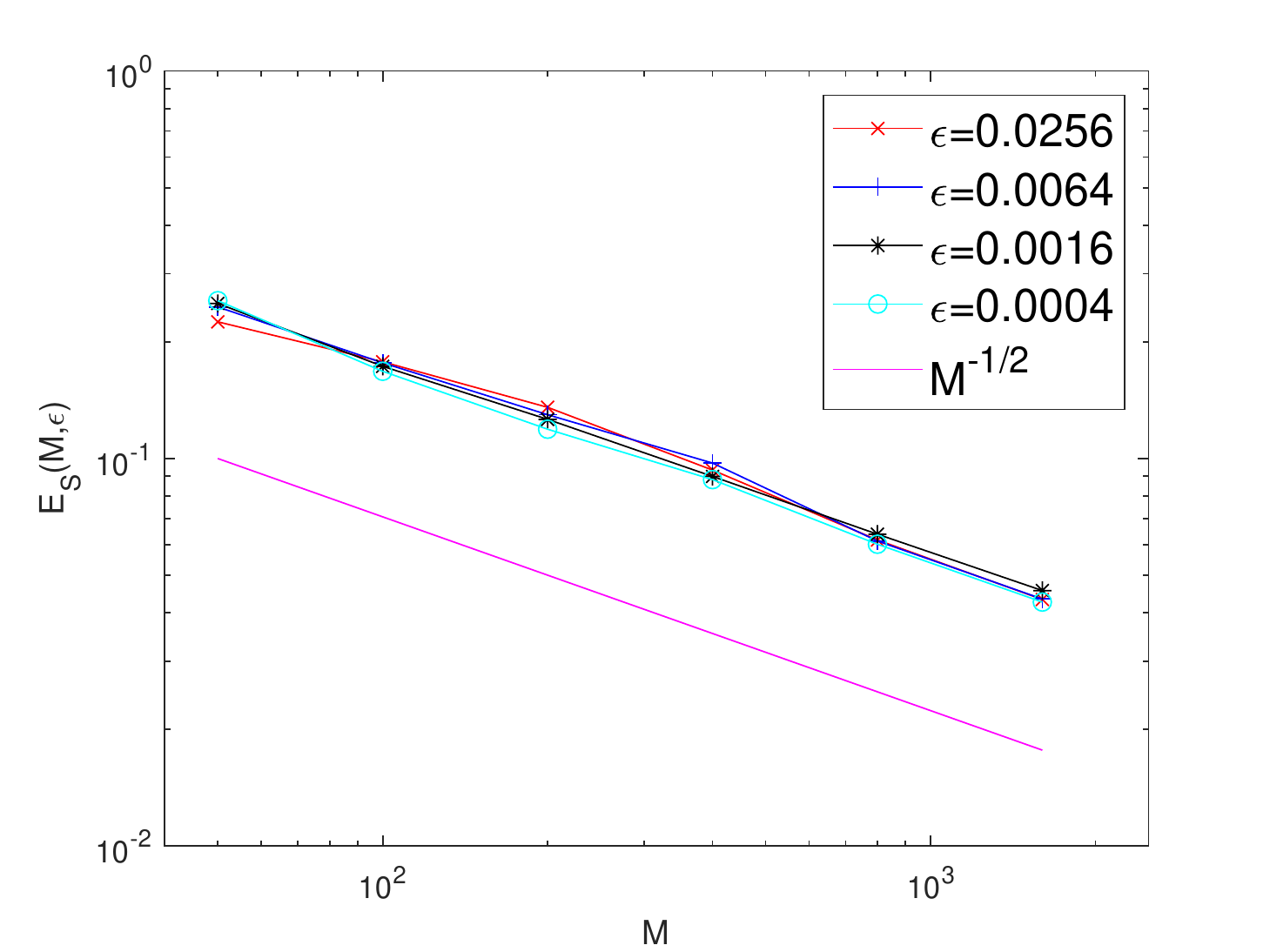}
\includegraphics[width=6cm,height=5cm]{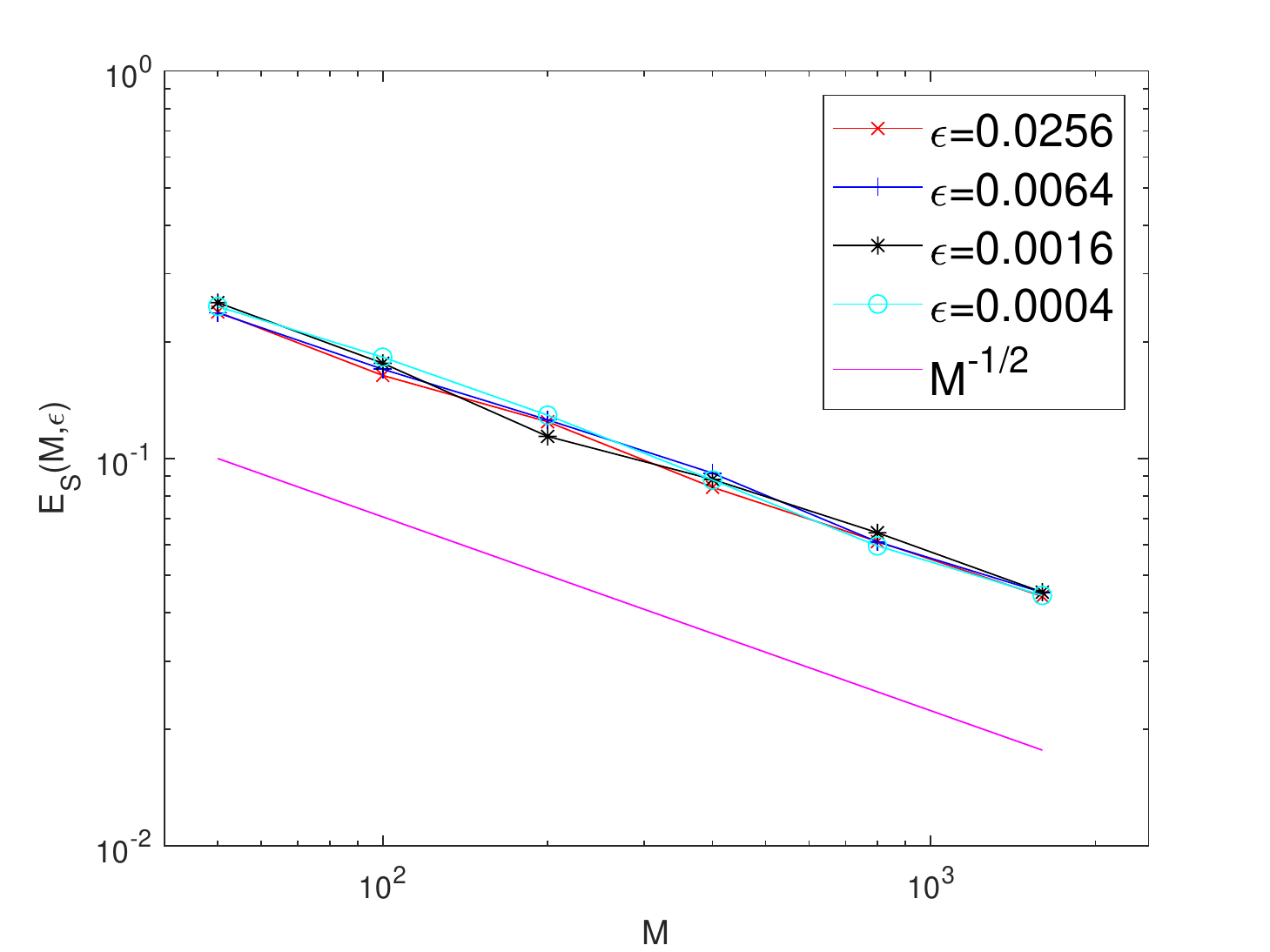}
\caption{\textbf{(Example 1) }The convergence of \(E_S\left(M,\varepsilon\right)\) with respect to the sampling size \(M\) with different scaling parameters \(\ep\) and potential functions. Left: quadratic potential \(E_1\); Right: torsion potential \(E_2\).}
\label{figlo}
\end{figure}

In Table \ref{TabGau} and Figure \ref{figlo}, we present the sampling error \(E_S(M,\ep)\) of the FGS with for different scaling Planck parameter \(\ep=0.0256,0.0064,0.0016,0.0004\) and two potential functions \(E_1\) and \(E_2\) respectively. One can easily notice that \(E_S\) is independent of \(\ep\) for both potential functions and the convergence rate of the sampling error \(E_S\) with respect to the ensemble size \(M\) is of order \(\CO(M^{-\frac{1}{2}})\), which is in accordance with the theoretical results in Theorem \ref{ThmIFGA} and Corollary \ref{CorFGA}. We also underline that it only takes 400 samples for the sampling error \(E_S\) to reach \(\CO(10^{-2})\). As a reference, we have stated in Section 3 that the computational cost of the mesh-based algorithms to reach certain error should be at least \(\CO(\ep^{-1})\), which indicates that the FGS is a practical and accurate algorithm to numerically approximate Equation \eqref{SemiSchrodinger} with 1D Gaussian initial data, especially for the cases when \(\ep\) is small.

\begin{center}
  \begin{tabular}{c|cccc}
   \toprule
   \textbf{Potential \(E_1\)} & \(\varepsilon=0.0256\) & \(\varepsilon=0.0064\) & \(\varepsilon=0.0016\) & \(\varepsilon=0.0004\)  \\
   \midrule
   \(M=50\) & 2.25e-01 & 2.44e-01 & 2.51e-01 & 2.55e-01  \\
   \(M=100\) & 1.78e-01 & 1.77e-01 & 1.73e-01 & 1.68e-01  \\
   \(M=200\) & 1.36e-01 & 1.30e-01 & 1.26e-01 & 1.19e-01  \\
   \(M=400\) & 9.34e-02 & 9.73e-02 & 8.99e-02 & 8.82e-02  \\
   \(M=800\) & 6.16e-02 & 6.12e-02 & 6.38e-02 & 6.01e-02  \\
   \(M=1600\) & 4.33e-02 & 4.34e-02 & 4.56e-02 & 4.26e-02   \\
   \midrule
   \textbf{Potential \(E_2\)} & \(\varepsilon=0.0256\) & \(\varepsilon=0.0064\) & \(\varepsilon=0.0016\) & \(\varepsilon=0.0004\)  \\
   \midrule
   \(M=50\) & 2.39e-01 & 2.38e-01 & 2.53e-01 & 2.47e-01  \\
   \(M=100\) & 1.64e-01 & 1.70e-01 & 1.76e-01 & 1.83e-01  \\
   \(M=200\) & 1.24e-01 & 1.26e-01 & 1.14e-01 & 1.29e-01  \\
   \(M=400\) & 8.43e-02 & 9.17e-02 & 8.83e-02 & 8.80e-02  \\
   \(M=800\) & 6.11e-02 & 6.09e-02 & 6.43e-02 & 5.95e-02  \\
   \(M=1600\) & 4.41e-02 & 4.51e-02 & 4.52e-02 & 4.43e-02   \\
   \bottomrule
  \end{tabular}
 \captionof{table}{\textbf{(Example 1) }The sampling error \(E_S\left(M,\varepsilon\right)\) of the FGS with different sampling size \(M\) and scaling parameters \(\ep\).}
 \label{TabGau}
 \end{center}

\subsubsection*{Example 2: 2D Gaussian initial data example}

As we have stated, the FGS wave functions is constructed by an ensemble of Gaussian wave packets \(\Lambda\kl z_0^{(j)}\kr\) (see Equation \eqref{Lambda} for definition) related to an ensemble of FGA trajectories with a certain probability distribution. In Example 1, the potential functions are so simple that the exact solutions to Equation \eqref{SemiSchrodinger} are still Gaussian wave packets at time \(t=0.5\), which lowers the difficulty of constructing the wave function from the ensemble of Gaussian wave packets \(\Lambda\kl z_0^{(j)}\kr\). In this example, we implement the FGS algorithm on 2D semiclassical Schr\"odinger equation with more complicated potential functions, where the exact solutions to Equation \eqref{SemiSchrodinger} are not Gaussian wave packets as time evolves. Though it is much more difficult to approximate a wave function that is not Gaussian from an ensembles of Gaussian wave packets, the FGS method is able to capture the correct evolution behavior with a reasonably small simple size.

\begin{figure}
\centering
\includegraphics[width=6cm,height=5cm]{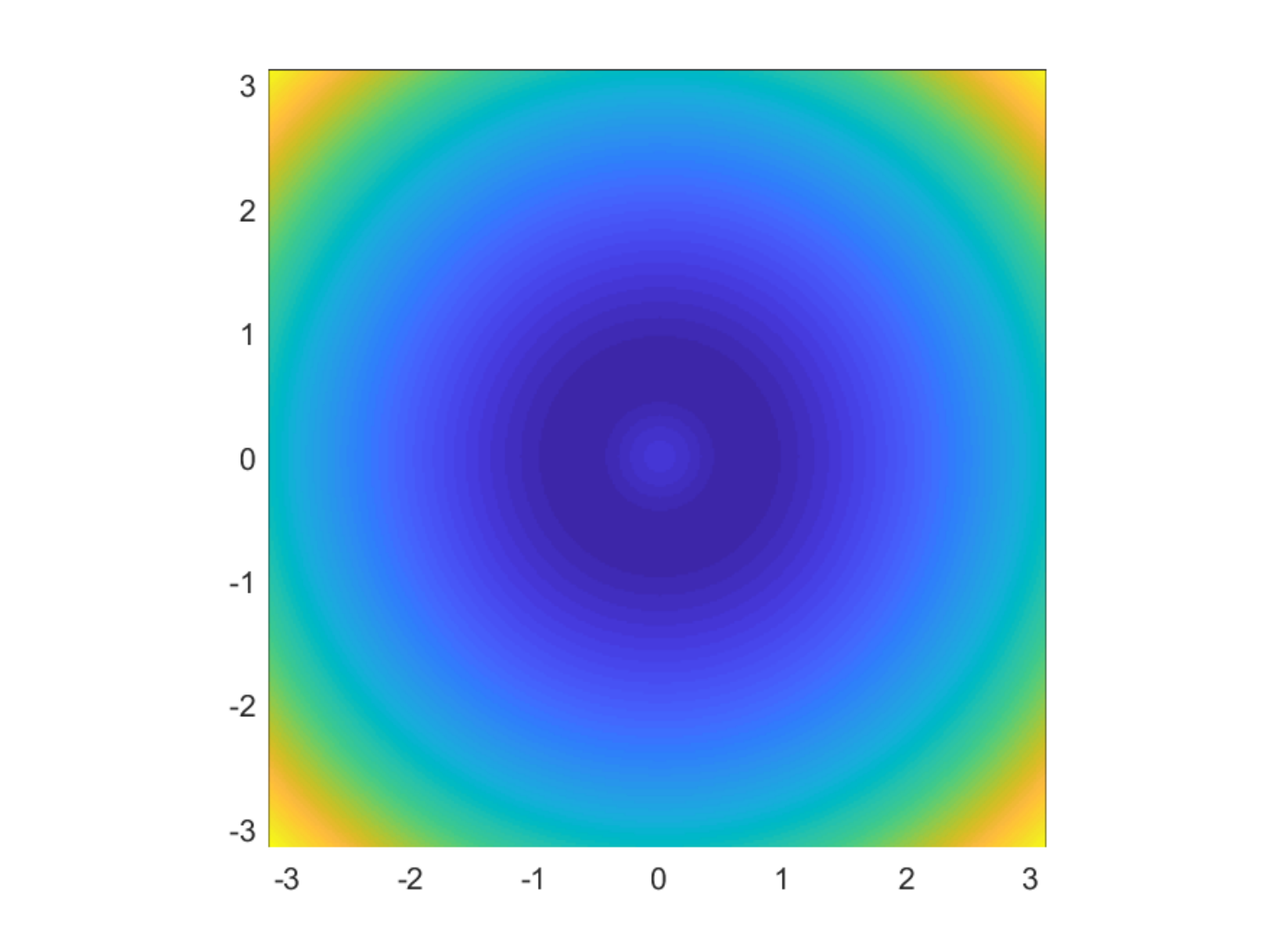}
\includegraphics[width=6cm,height=5cm]{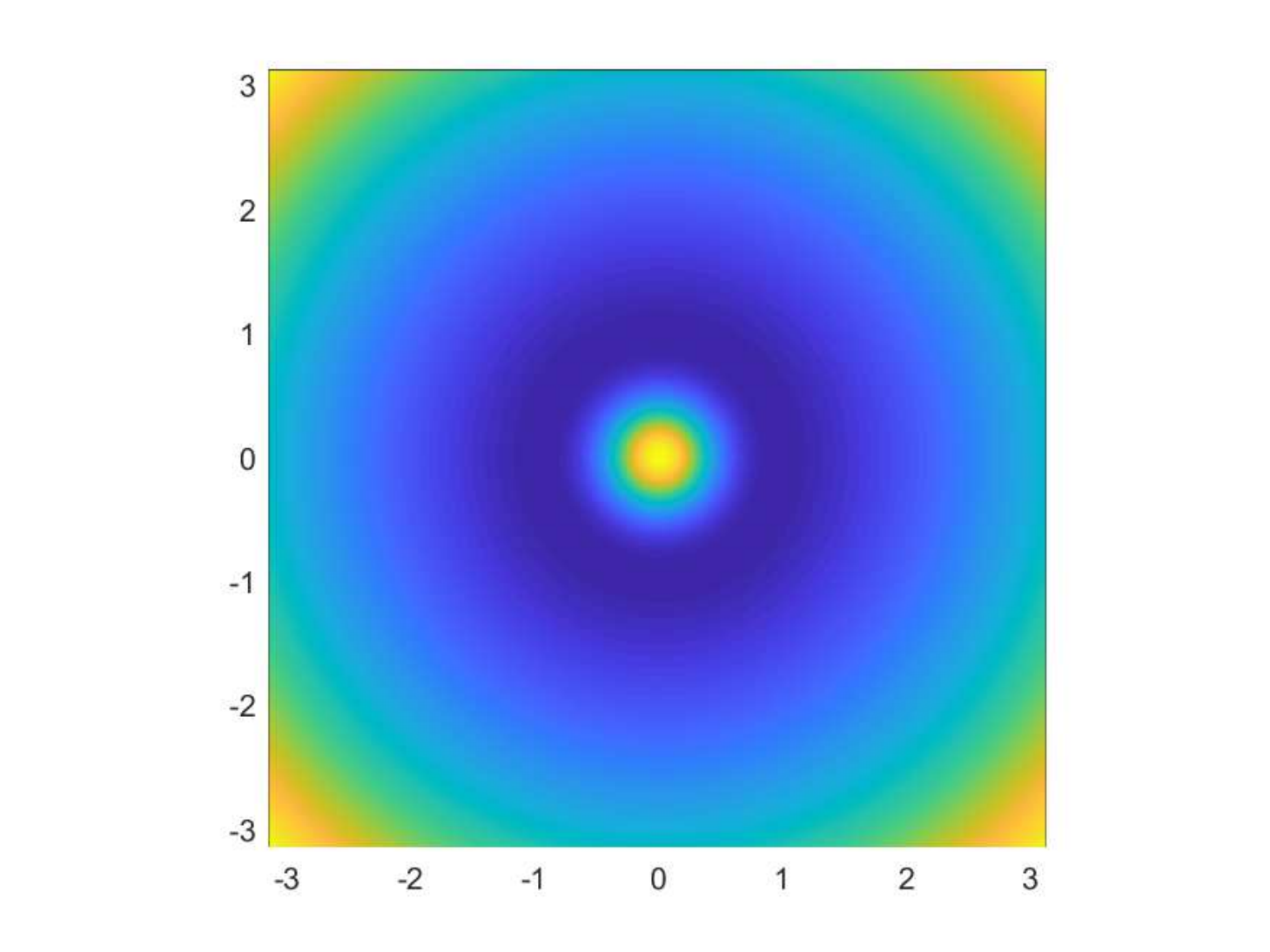}
\caption{\textbf{(Example 2) }An illustration of the 2D potential functions \(E_3\) and \(E_4\). Left: \(E_3\); Right: \(E_4\).}
\label{figsc}
\end{figure}

In the Gaussian initial condition \eqref{InitialGaussian}, let \(m=2\), \(a_1=a_2=2\), \(\tilde{p}=(-1,-1)\) and \(\tilde{q}=(1,1)\). The scaling parameter is fixed at \(\ep=\frac{1}{96}\). Then we obtain the following 2D Gaussian wave packet:
\beq\label{EG2wave}
u_\ti(x)=\sqrt{\pi\ep}\exp\kl-\frac{\ii}{\ep}(x_1-1)-\frac{\ii}{\ep}(x_2-1)-\frac{1}{2\ep}(x_1-1)^2-\frac{1}{2\ep}(x_2-1)^2\kr.
\eeq
Now we consider two complicated potential functions \(E_3(x)=\frac{|x|^2}{2}+\ee^{-5|x|^2}\) and \(E_4(x)=\frac{|x|^2}{2}+10\ee^{-5|x|^2}\). We illustrate of the potential \(E_3\) and \(E_4\) in Figure \ref{figsc}. Intuitively, due to the effect of the exponential term, the shapes of the potential functions \(E_3\) and \(E_4\) show peaks at the center. The main difference between the two potentials lies in that the peak of \(E_4\) is much higher than \(E_3\), which leads to different behaviors of the exact solutions to Equation \eqref{SemiSchrodinger}.

In this example, we apply the Strang Splitting Spectral Method (SP2) \cite{BJM2003} with sufficiently small step lengths to compute the exact solution to Equation \eqref{SemiSchrodinger}. We reconstruct the FGS wave function on spatial region \([-\pi,\pi]\times[-\pi,\pi]\) through Algorithm \ref{alg:FGA} and compare the resulting FGS density function
\[
\rho_\fgs(t,x)=\left|u_\fgs(t,x)\right|^2,
\]
with the exact density function \(\rho(t,x)\) approximated by the SP2. Similar to the previous example, we use the fourth order Runge-Kutta method to solve the ODE system \eqref{PhaseQ}-\eqref{PhaseA} in the time interval \([0,0.5]\) with time step length \(\Delta t=0.01\) and reconstruct the FGS wave function at different times. We choose the ensemble size \(M=2000\). Note that \(M=2000\) is of lower order than \(\CO(\ep^{-2})\), which is the degree of freedom of the mesh-based method SP2. 

\begin{figure}
\centering
\includegraphics[width=3.5cm,height=3cm]{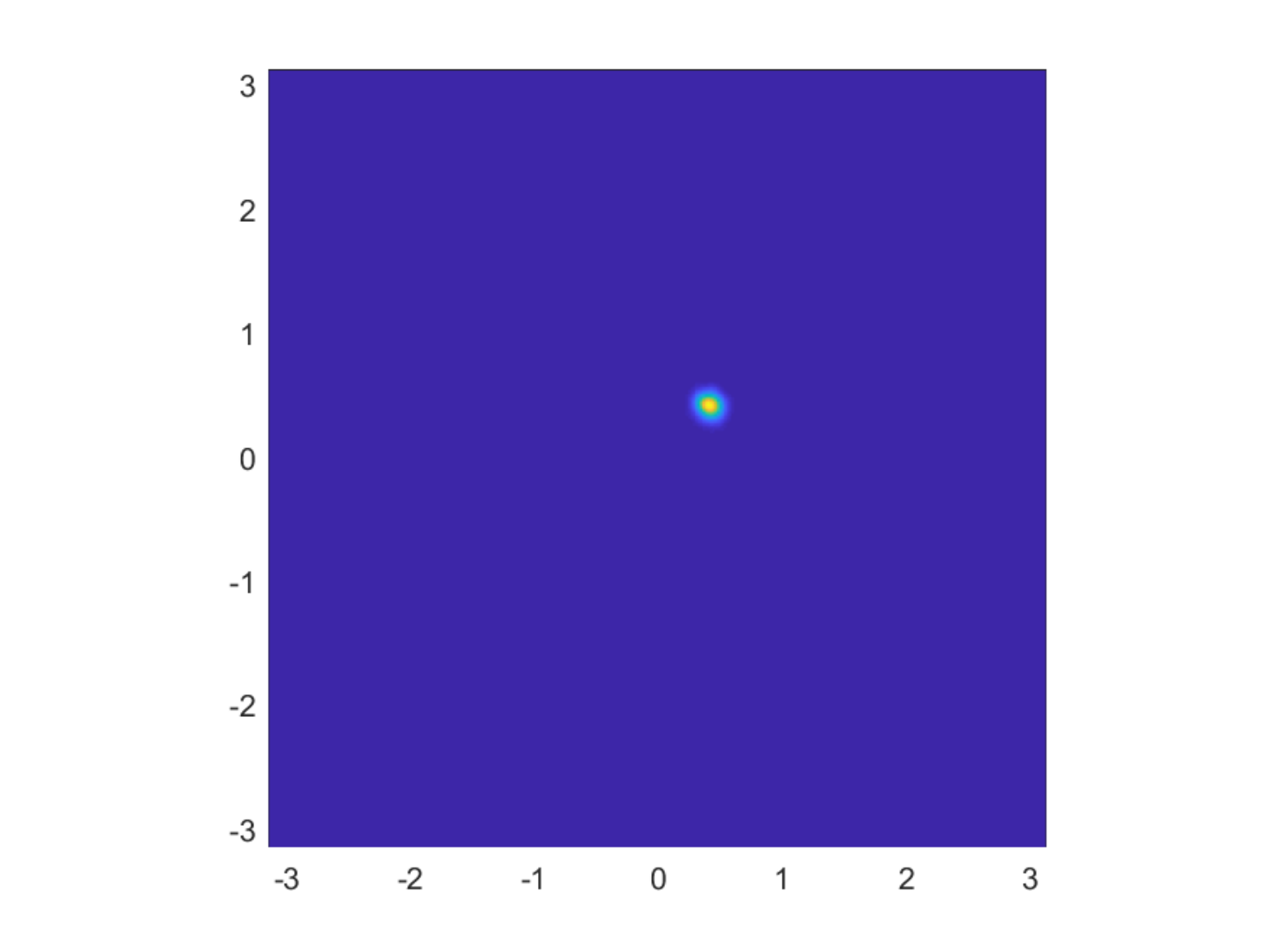}
\includegraphics[width=3.5cm,height=3cm]{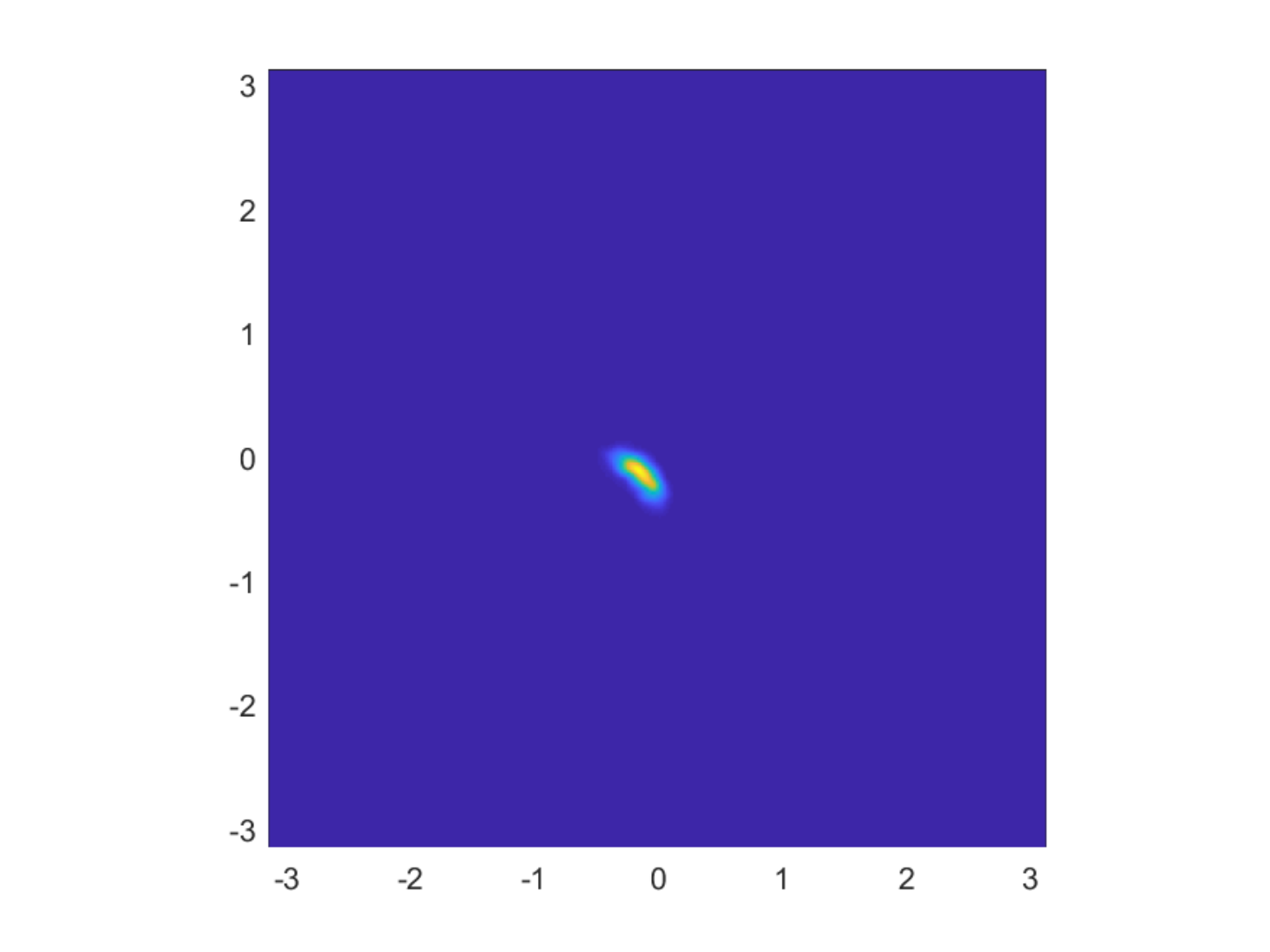}
\includegraphics[width=3.5cm,height=3cm]{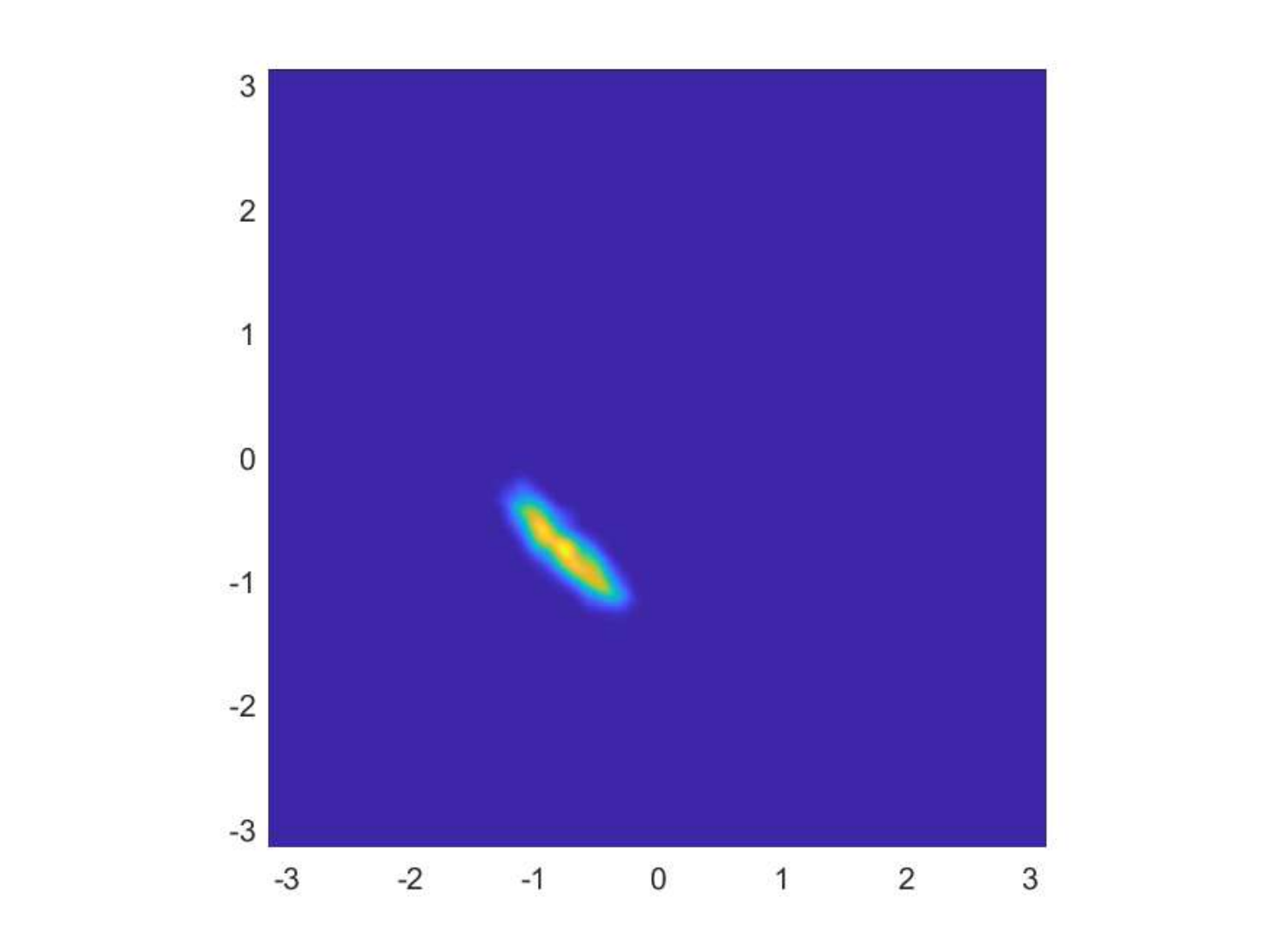}
\includegraphics[width=3.5cm,height=3cm]{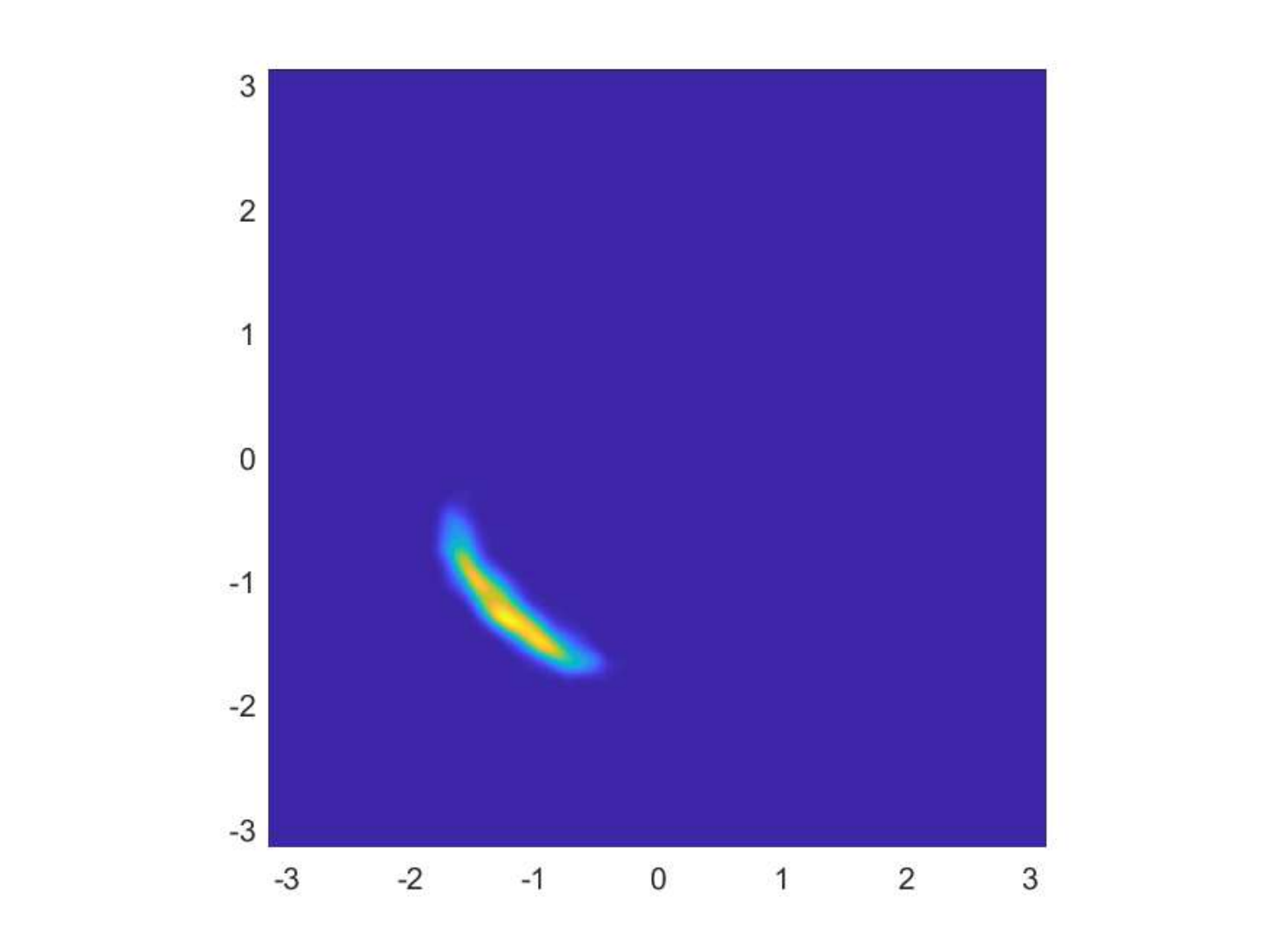}
\includegraphics[width=3.5cm,height=3cm]{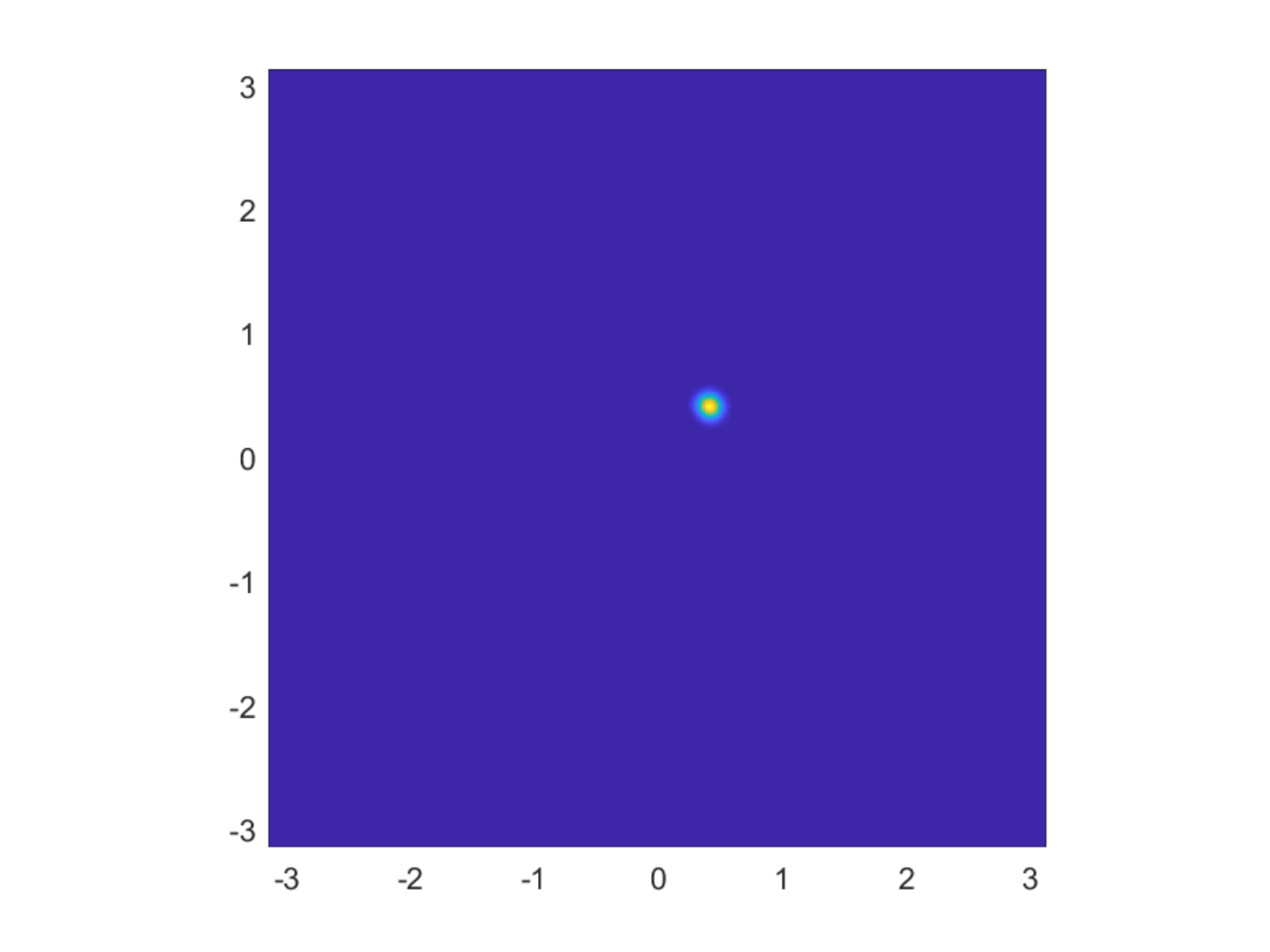}
\includegraphics[width=3.5cm,height=3cm]{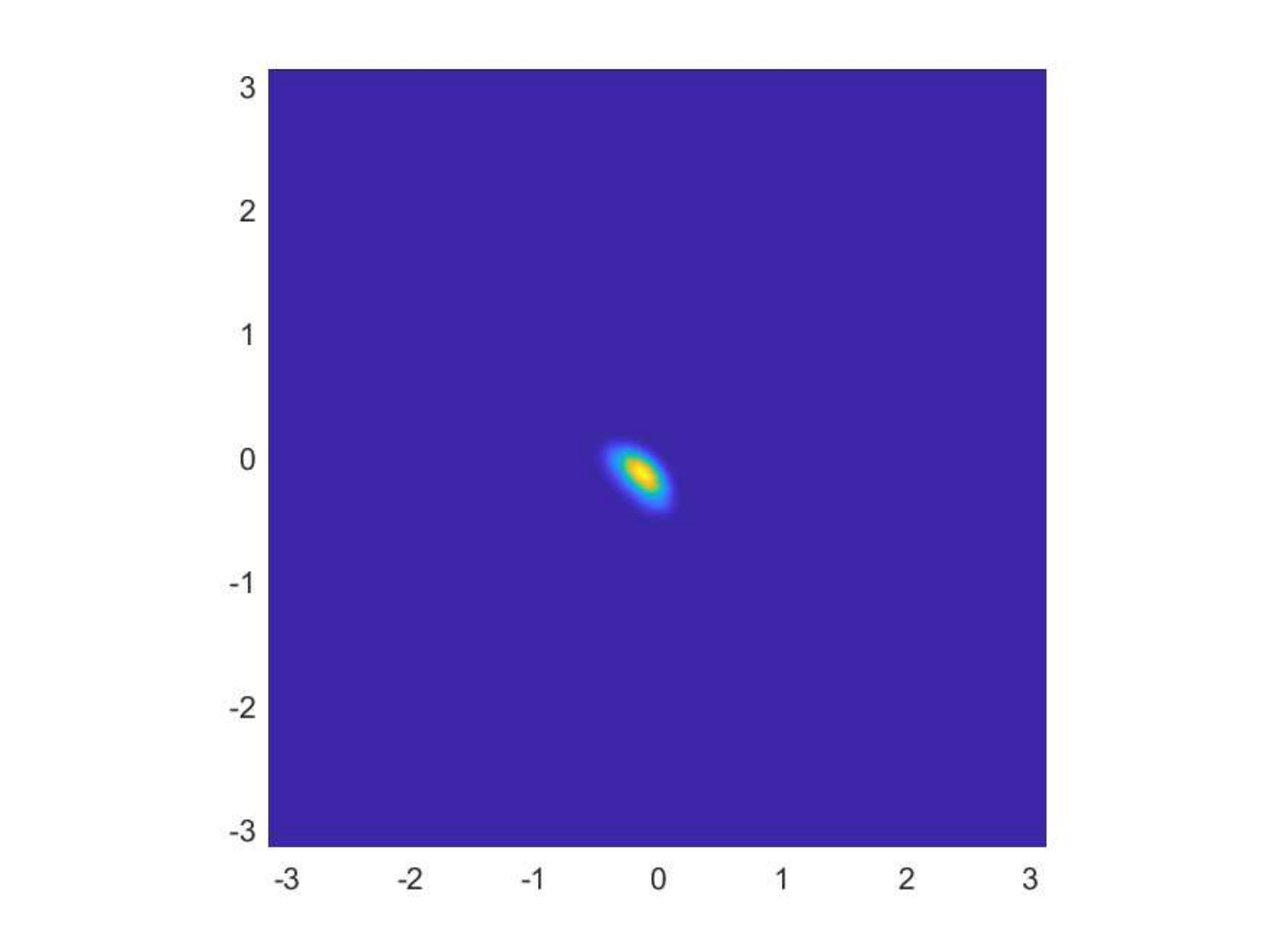}
\includegraphics[width=3.5cm,height=3cm]{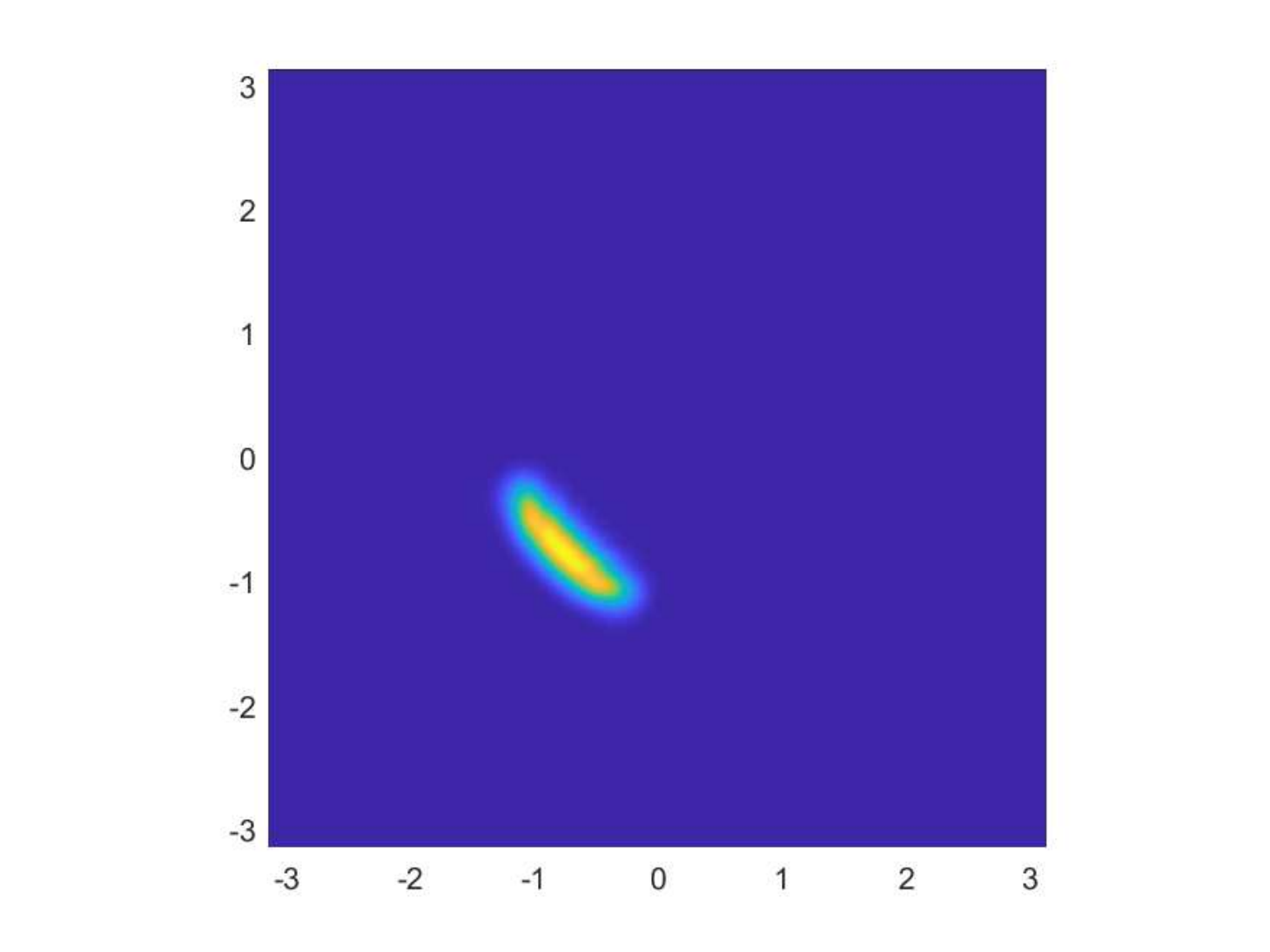}
\includegraphics[width=3.5cm,height=3cm]{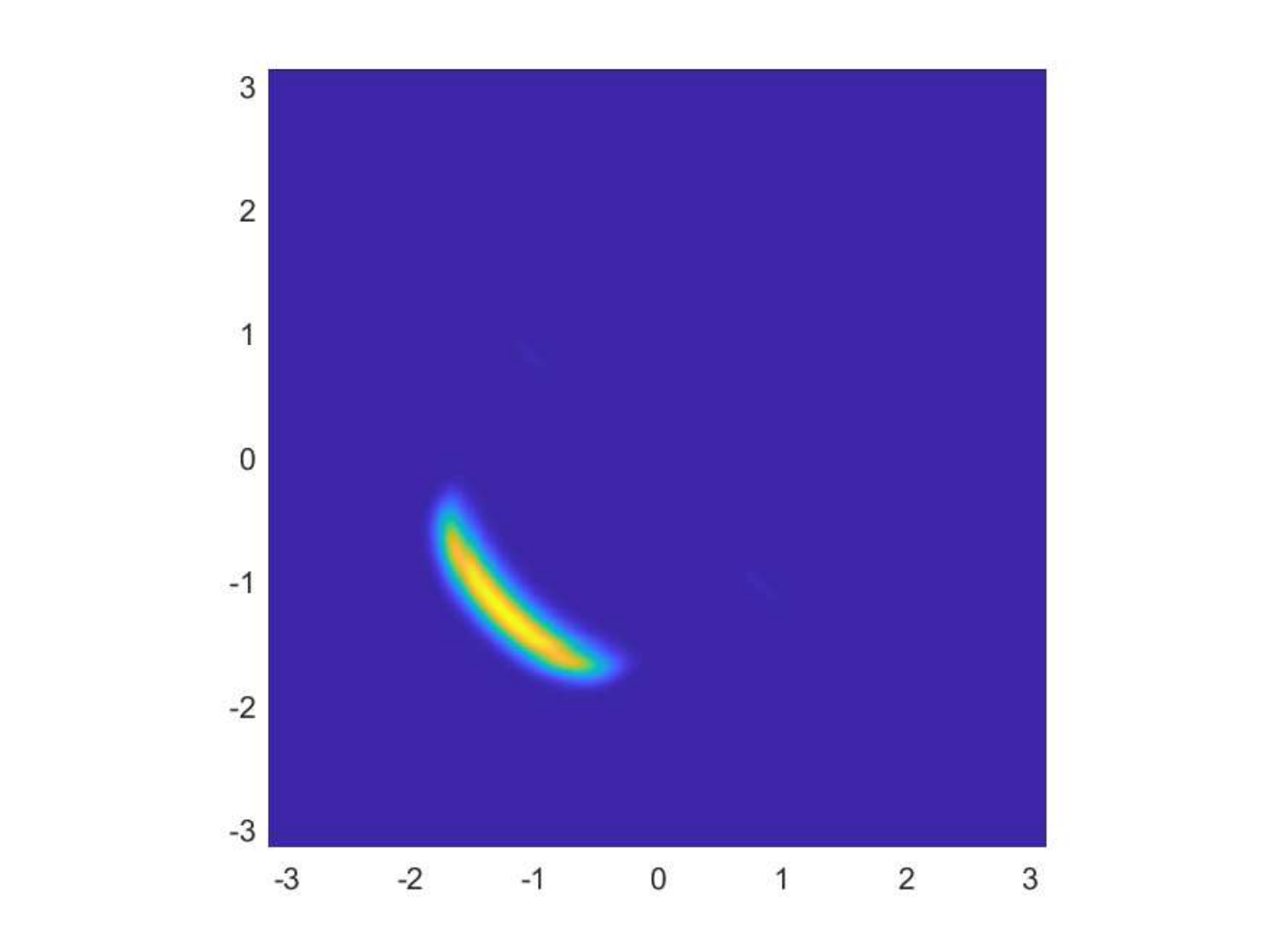}
\caption{\textbf{(Example 2) }We choose the potential function \(E_3\) and the figures show the evolution of the FGS density functions \(\rho_\fgs(t,x)\) compared with the exact density functions computed via the SP2 method. The top row: the FGS density function; The bottom row: the exact density function. From the left column to the right column, the figures show the density function at time \(t=0.5,1,1.5,2\).}
\label{figSC1}
\end{figure}

\begin{figure}
\centering
\includegraphics[width=3.5cm,height=3cm]{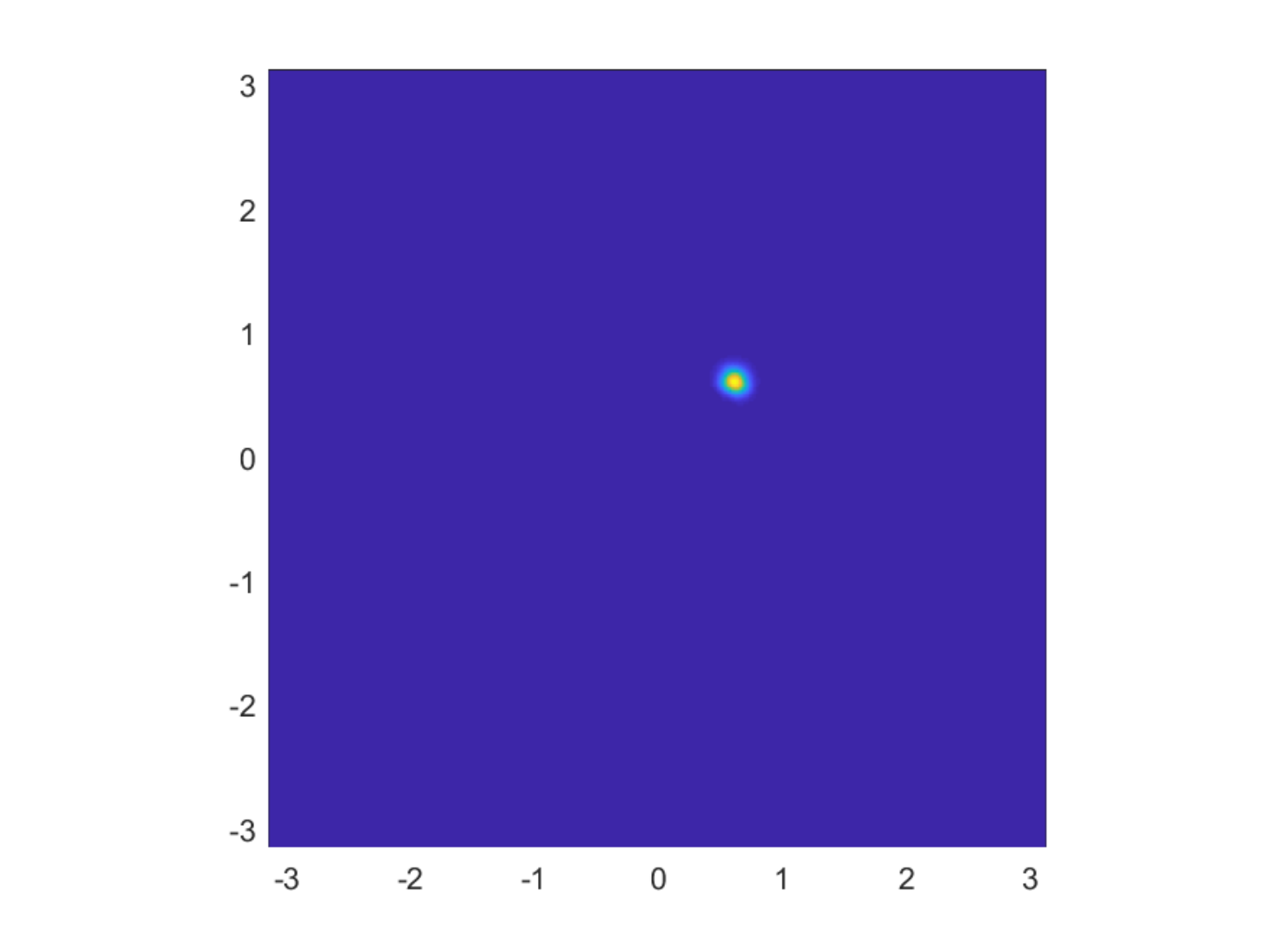}
\includegraphics[width=3.5cm,height=3cm]{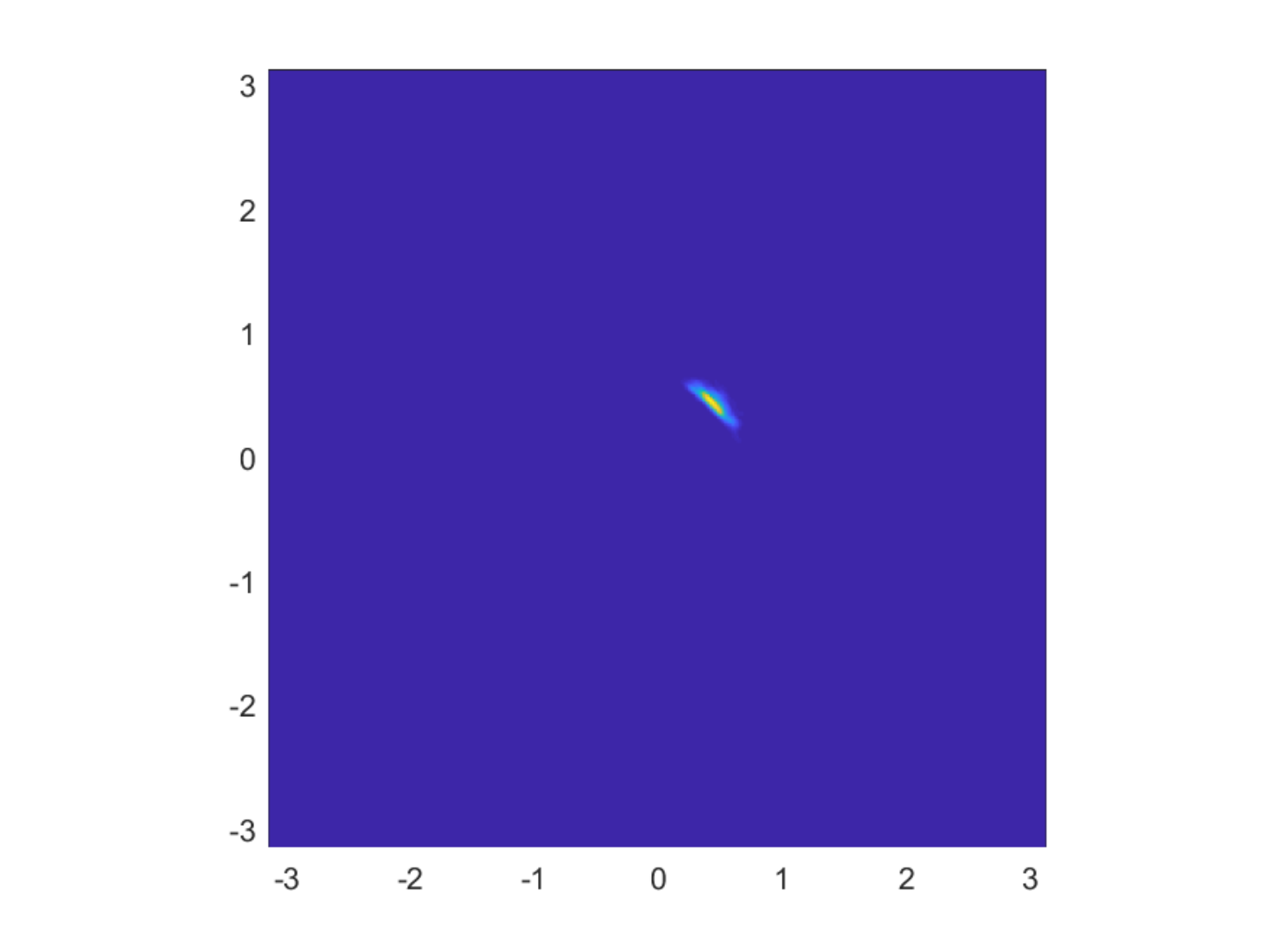}
\includegraphics[width=3.5cm,height=3cm]{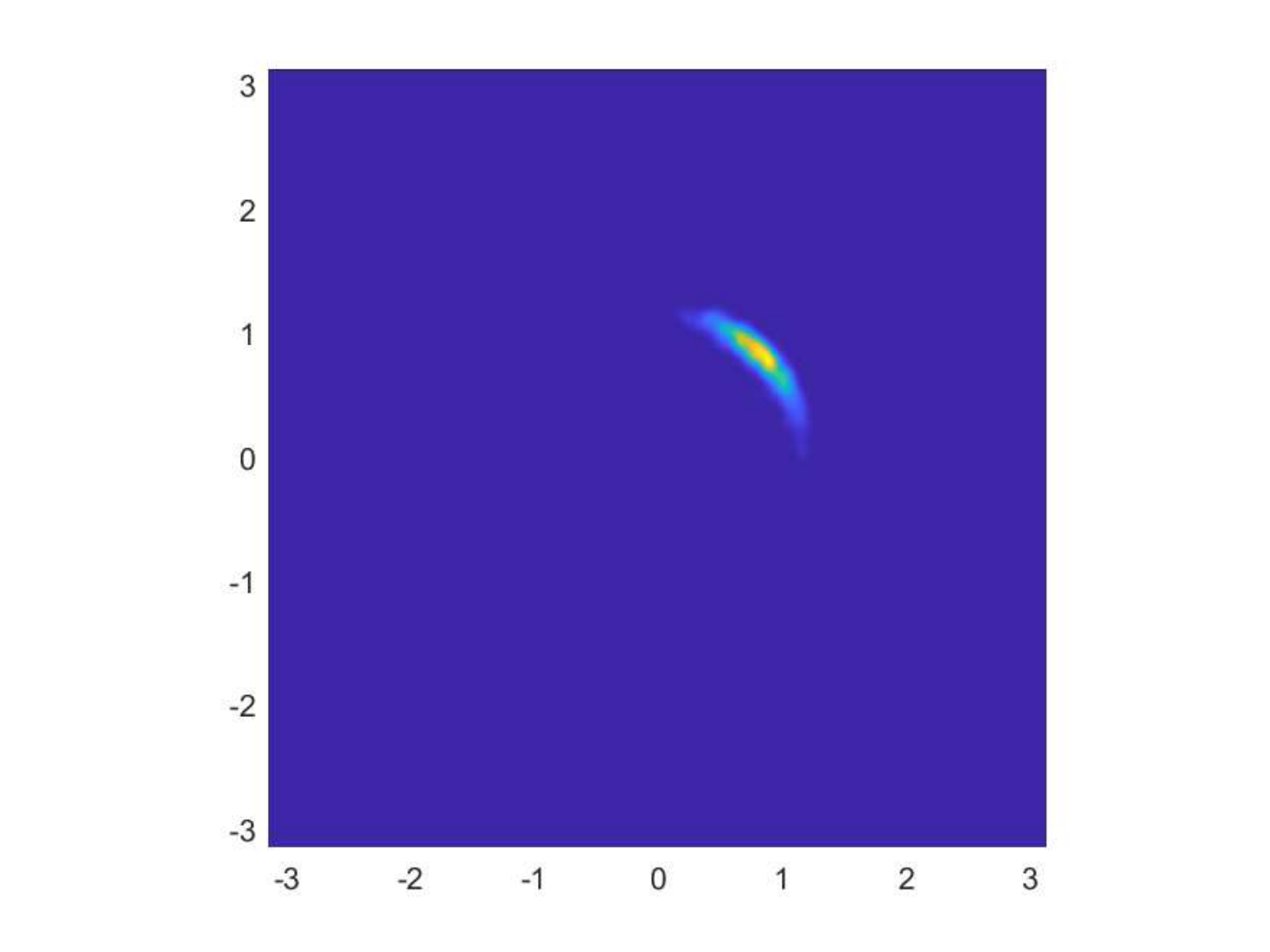}
\includegraphics[width=3.5cm,height=3cm]{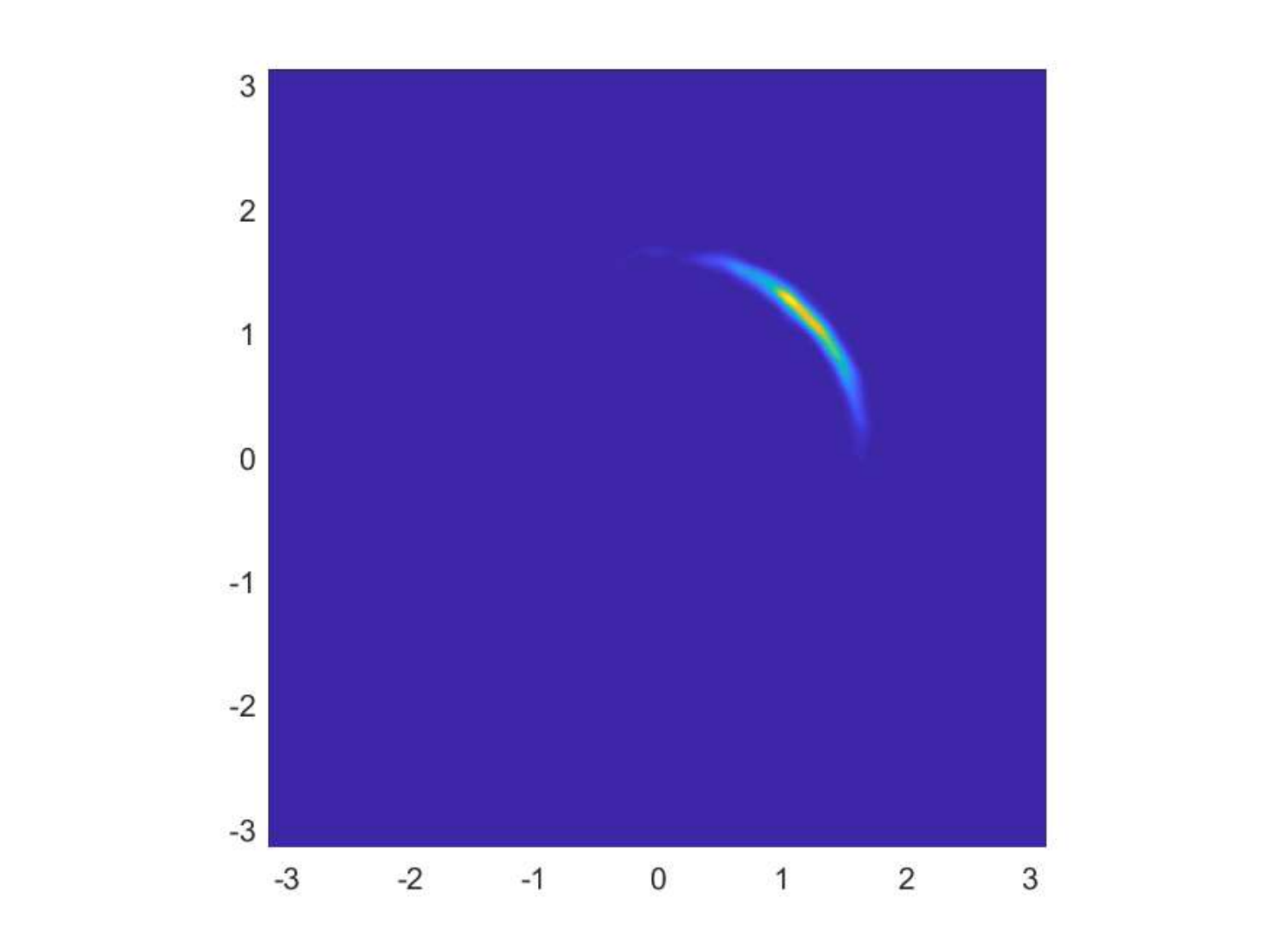}
\includegraphics[width=3.5cm,height=3cm]{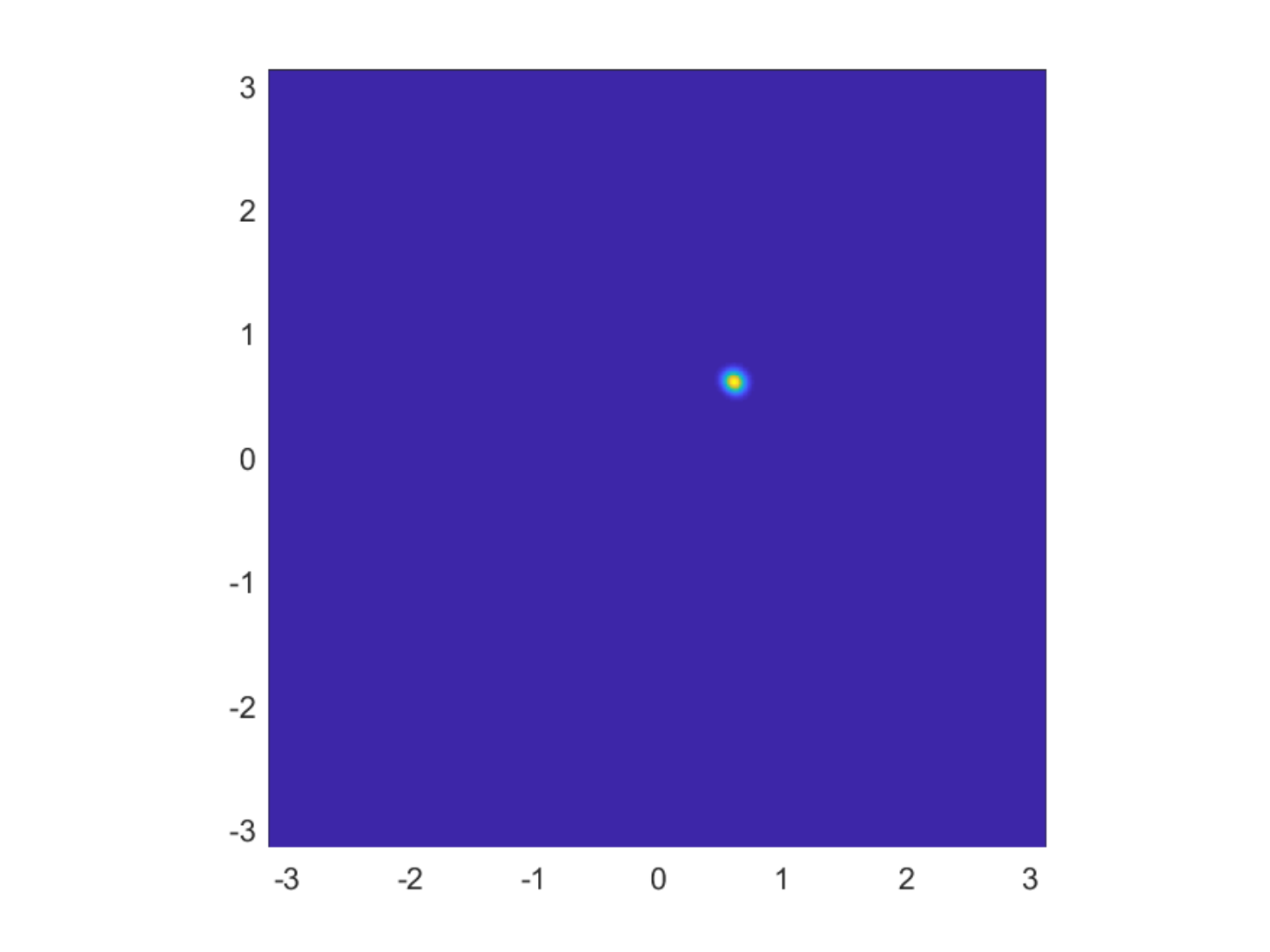}
\includegraphics[width=3.5cm,height=3cm]{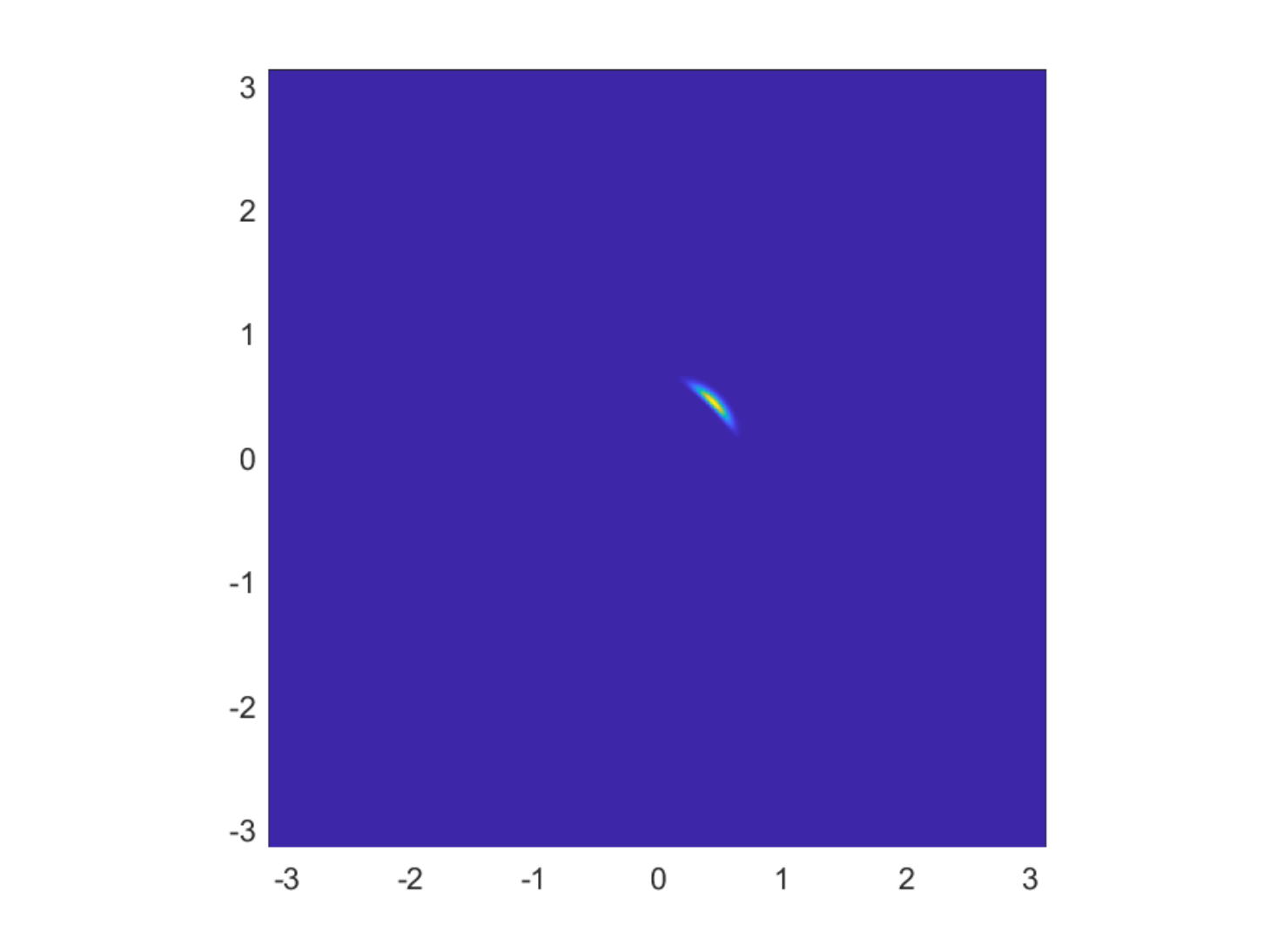}
\includegraphics[width=3.5cm,height=3cm]{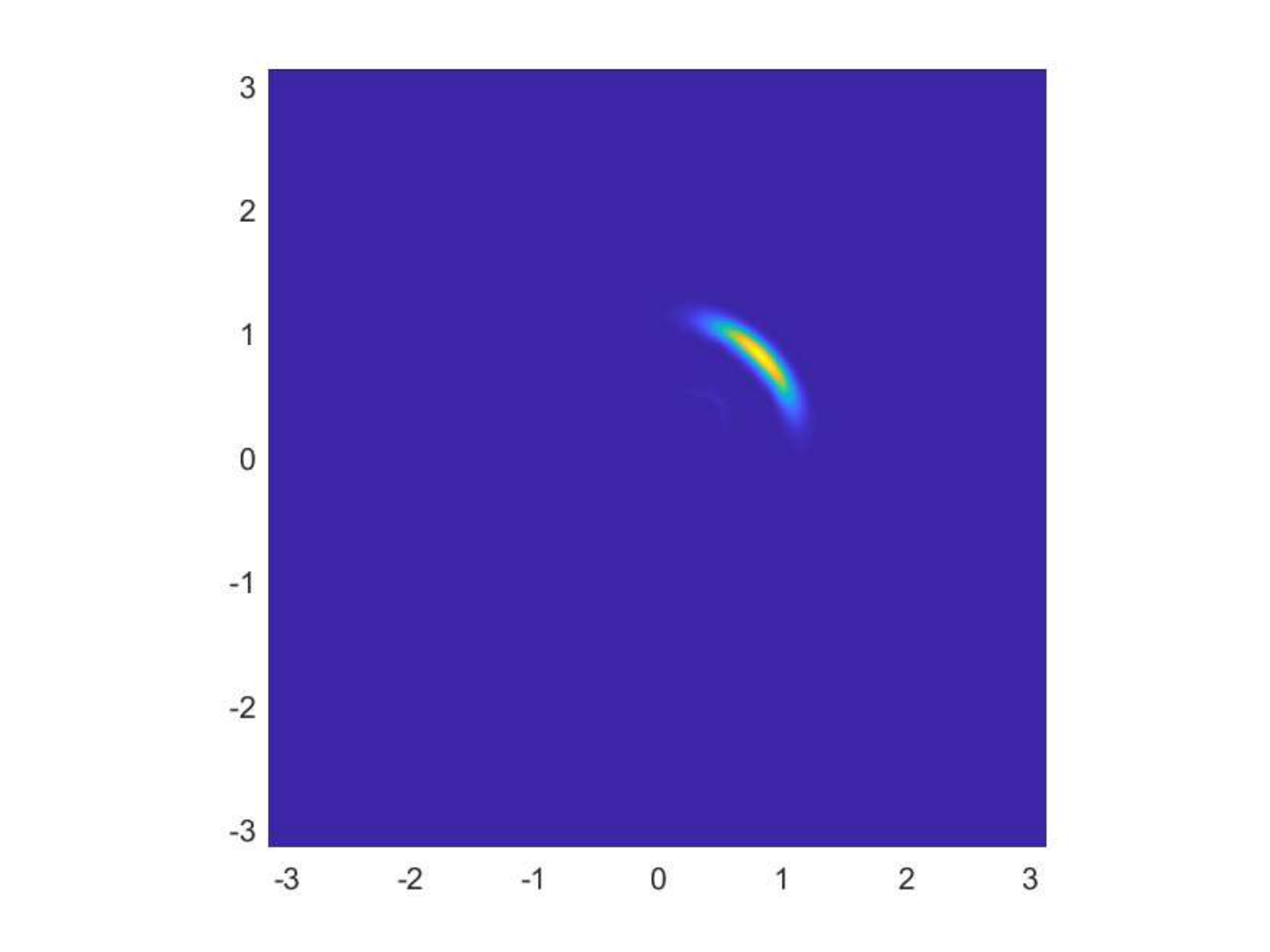}
\includegraphics[width=3.5cm,height=3cm]{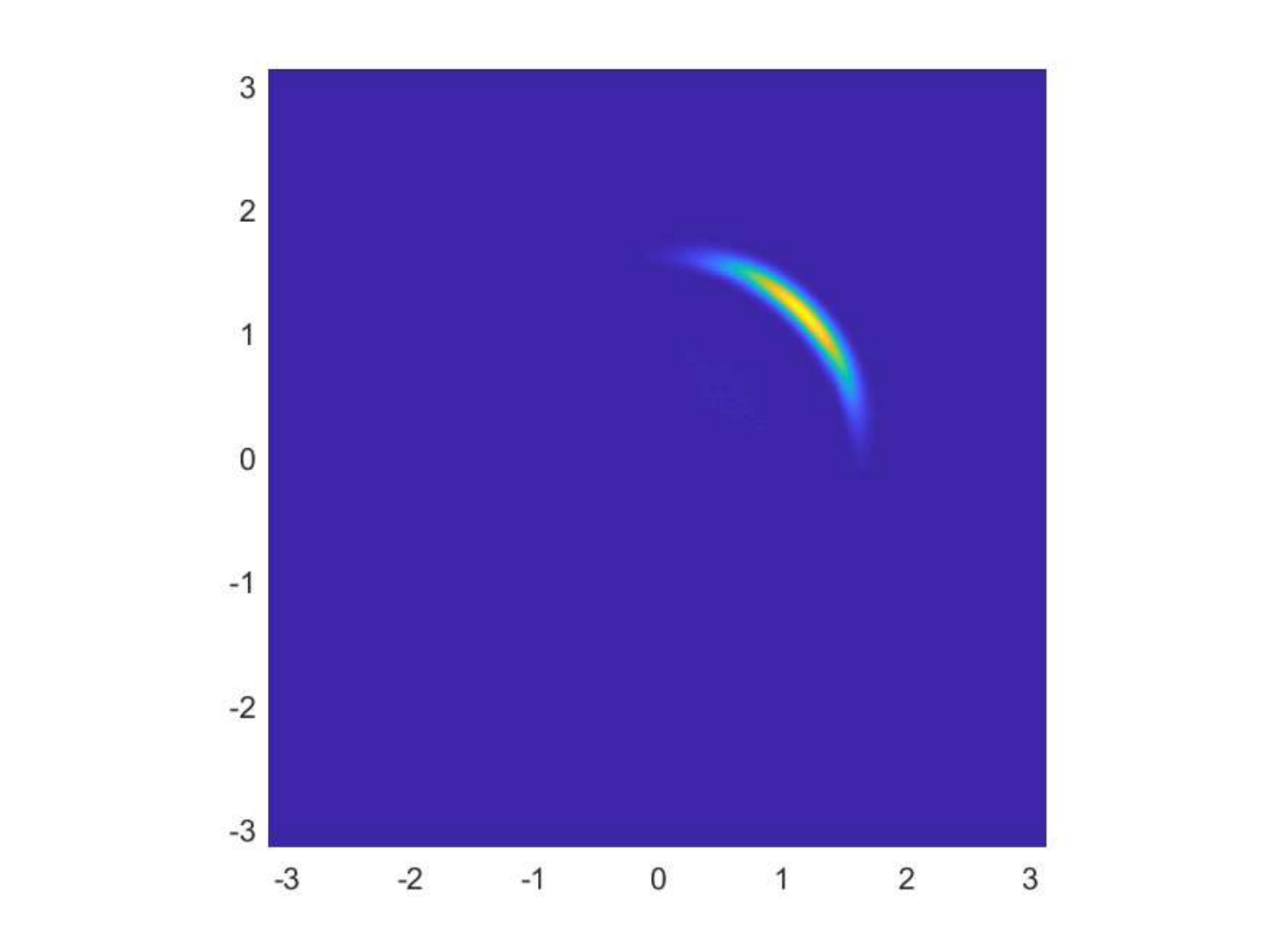}
\caption{\textbf{(Example 2) }We choose the potential function \(E_4\) and the figures show the evolution of the FGS density functions \(\rho_\fgs(t,x)\) compared with the exact density functions computed via the SP2 method. The top row: the FGS density function; The bottom row: the exact density function. From the left column to the right column, the figures show the density function at time \(t=0.35,0.7,1.05,1.4\).}
\label{figSC2}
\end{figure}

In Figure \ref{figSC1} and \ref{figSC2}, we show the evolution of the FGS density functions with two different potential functions \(E_3\) and \(E_4\) and compare them with the density functions computed via the SP2 algorithm. The density functions under both potential functions are no longer Gaussian after certain time. For \(E_3\) in Figure \ref{figSC1}, the wave packets pass through the potential peak at the center and disperse into different directions, which makes the wave packet crescent. In contrast, the potential peak of \(E_4\) is so high and steep that the wave packets in Figure \ref{figSC2} rebound against the potential peak and disperse into the opposite directions against the wave functions in Figure \ref{figSC1}. For both cases, the profile of the FGS density functions conform with the exact solutions when the wave packets are no longer Gaussian, which ensure the accuracy of the FGS method in computing non-Gaussian wave packets. In future work, we are ready to apply the FGS to approximating more complicated wave packets arising in scientific problems.

\subsection{Examples with WKB initial condition }

In this subsection, we apply the FGS with the initial sampling algorithm put forward in Section 4 to two examples of Equation \eqref{SemiSchrodinger} with WKB initial data. Caustics forms in the former example and does not form in the latter, which makes the former example more challenging for the FGS. To avoid confusion, we would like to point out that caustics become relevant to  Equation \eqref{SemiSchrodinger} only in the semiclassical limit, $\varepsilon \rightarrow 0$, but it also causes numerical challenges for moderately small $\varepsilon$.

\subsubsection*{Example 3: WKB initial data example with caustics.} 

In the WKB initial condition \eqref{WKBInitial}, let
\[
a_{\ti}(x)=\left(\frac{50}{\pi}\right)^{\frac{1}{4}}\ee^{-25(x-0.5)^{2}}, \quad S_{\ti}(x)=-\frac{1}{5} \ln \left(\ee^{5(x-0.5)}+\ee^{-5(x-0.5)}\right), \quad x \in \mathbb{R},
\]
we obtain the following WKB initial data:
\begin{equation}\label{WKBini1DC}
u_{\ti}(x)=\left(\frac{50}{\pi}\right)^\frac{1}{4}\exp\left[-25\left(x-\frac{1}{2}\right)^2-\frac{\ii}{5\varepsilon}\ln \left(e^{5(x-0.5)}+e^{-5(x-0.5)}\right)\right].
\end{equation}
We choose the constant potential \(E_5(x)=10\). Due to the compressive initial velocity \(\frac{\dd}{\dd x}S_\ti\), caustics will form. This example was used in \cite{MPP1999,BJM2003,JWY2008} as a reference to test whether certain algorithms can capture caustics.

In this example, we plot the FGS wave function and the resulting FGS current density function
\[
J_{\fgs}=\frac{\varepsilon}{2 \ii} \left(\ufgs  \overline{\nabla\ufgs}-\overline{\ufgs} \nabla \ufgs\right),
\]
and compared them with the reference solutions computed via the SP2 method with sufficiently small step lengths. The computation is carried out on spatial interval \([-\pi,\pi]\). Similarly, we use the fourth order Runge-Kutta method to solve the ODE system \eqref{PhaseQ}-\eqref{PhaseA} in the time interval \([0,0.54]\) with time step length \(\Delta t=0.01\)  and construct the wave function and current density at \(t=0.54\).

\begin{figure}
\centering
\includegraphics[width=6cm,height=5cm]{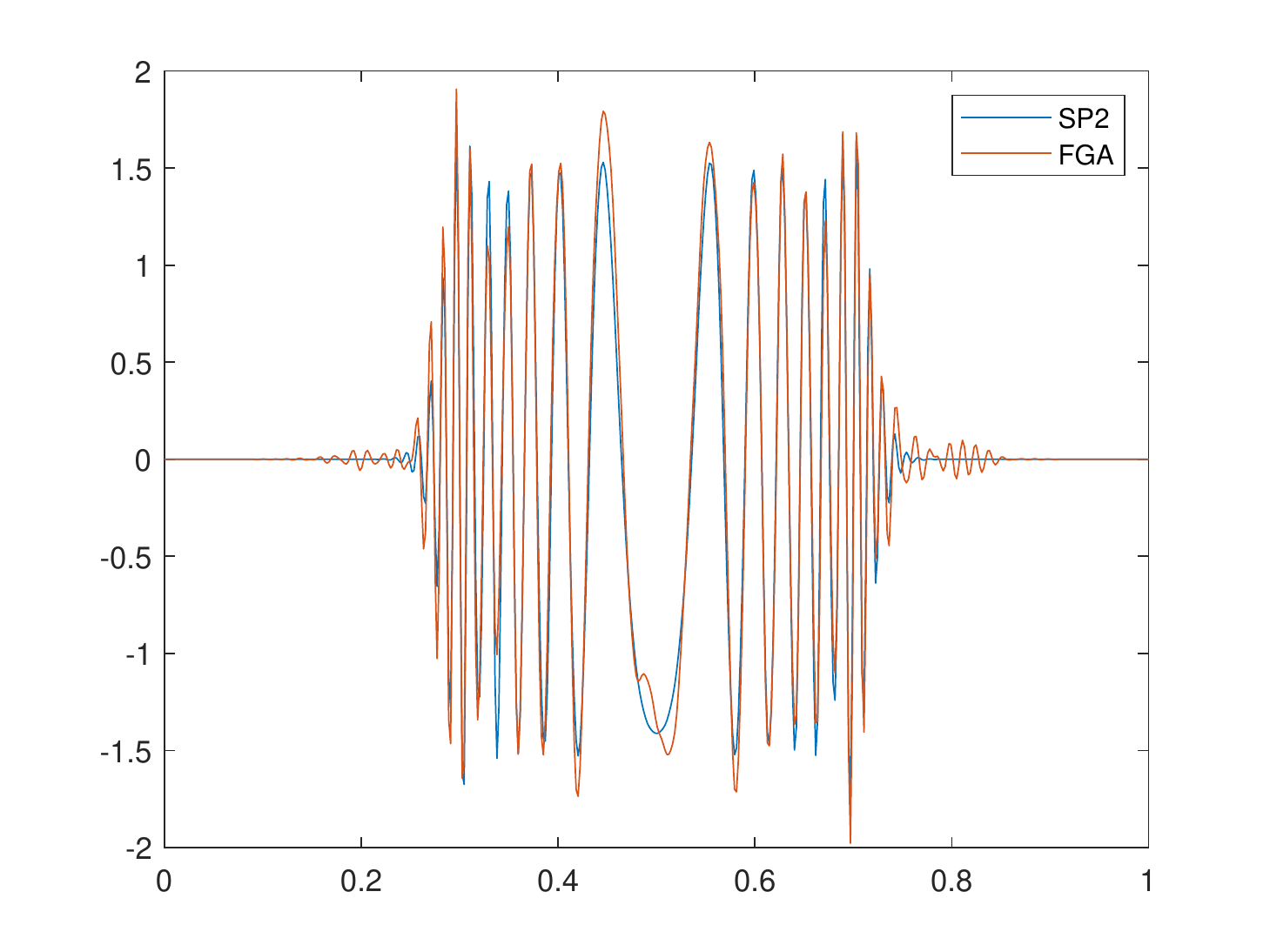}
\includegraphics[width=6cm,height=5cm]{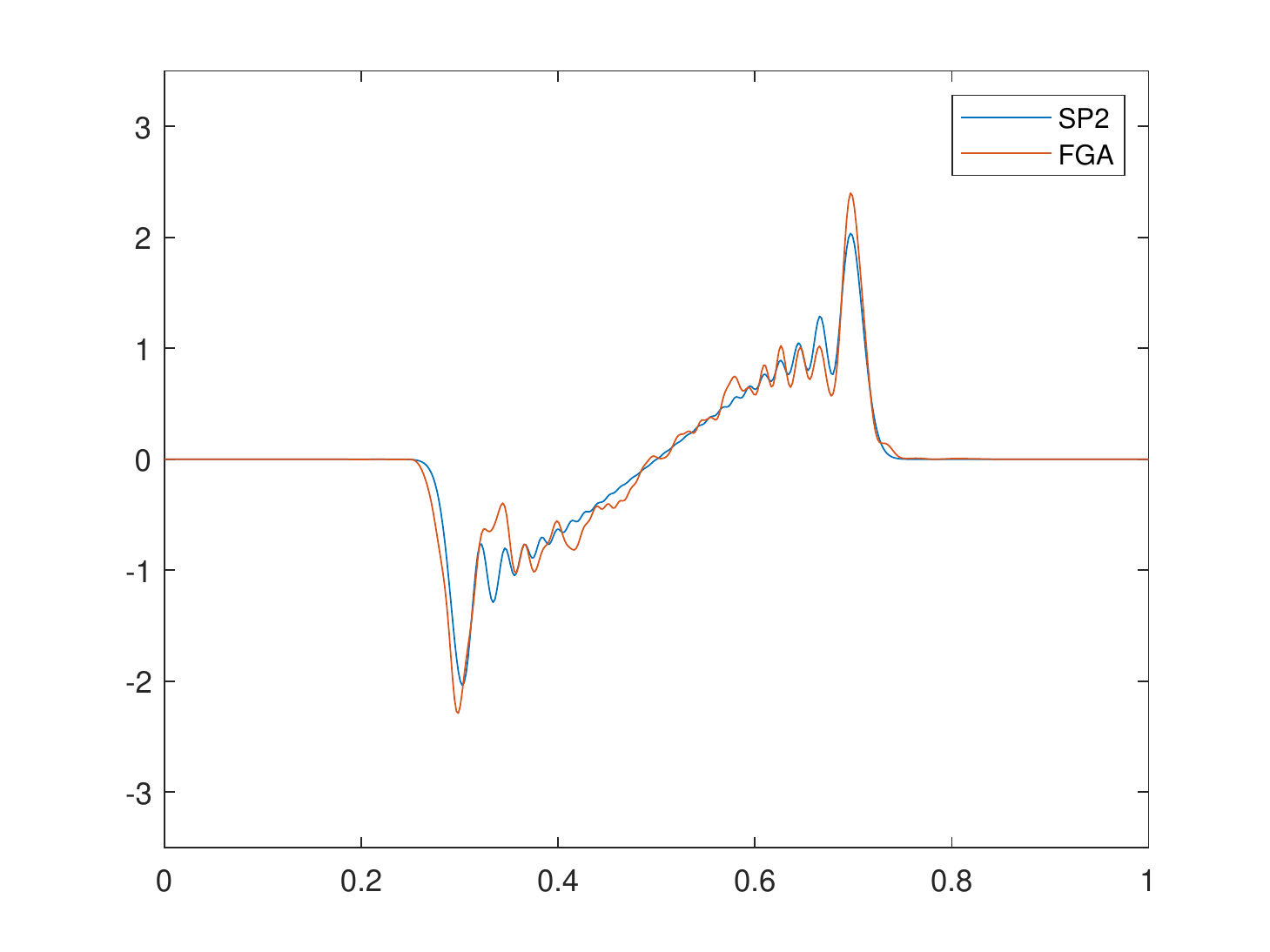}
\includegraphics[width=6cm,height=5cm]{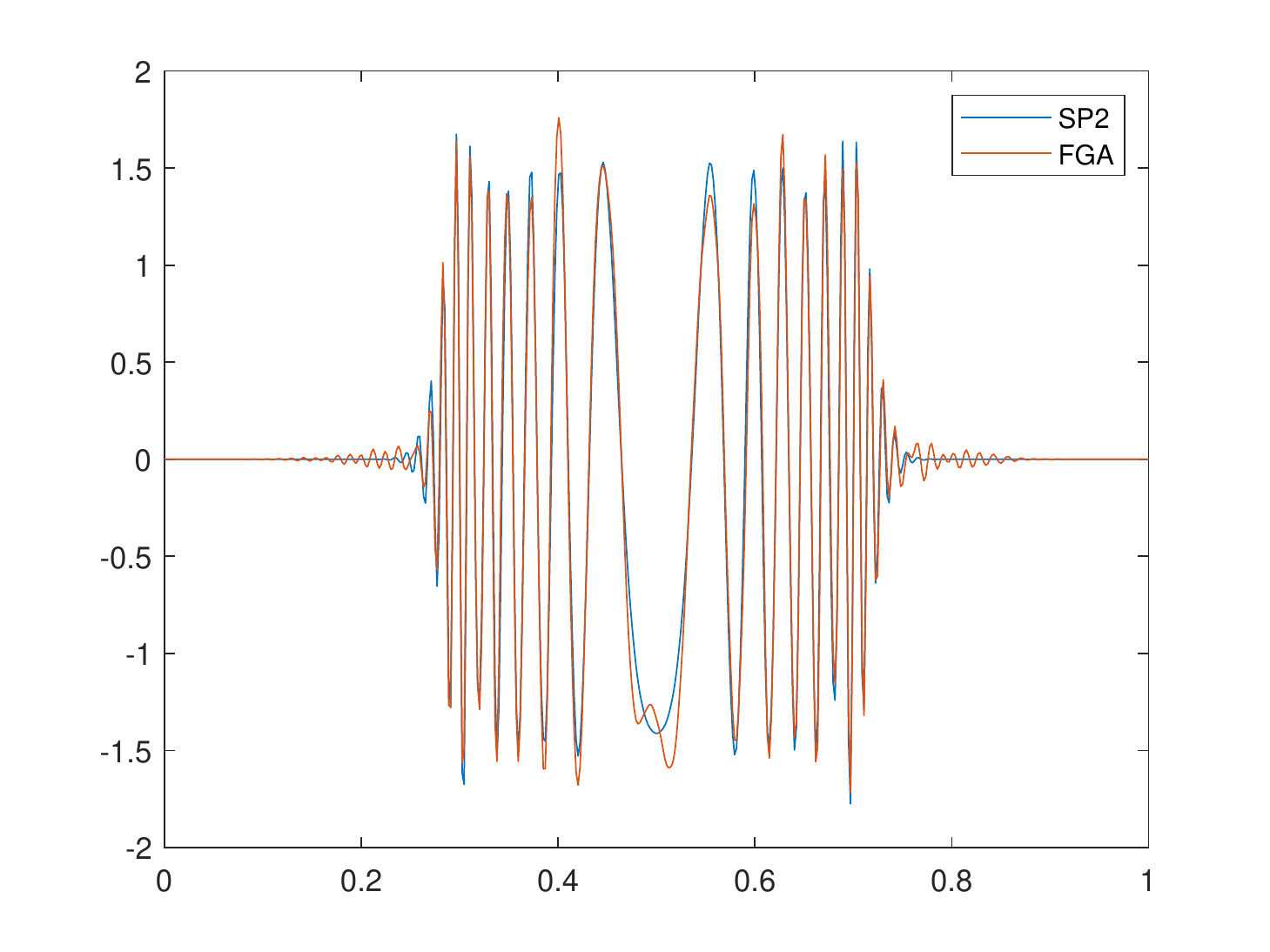}
\includegraphics[width=6cm,height=5cm]{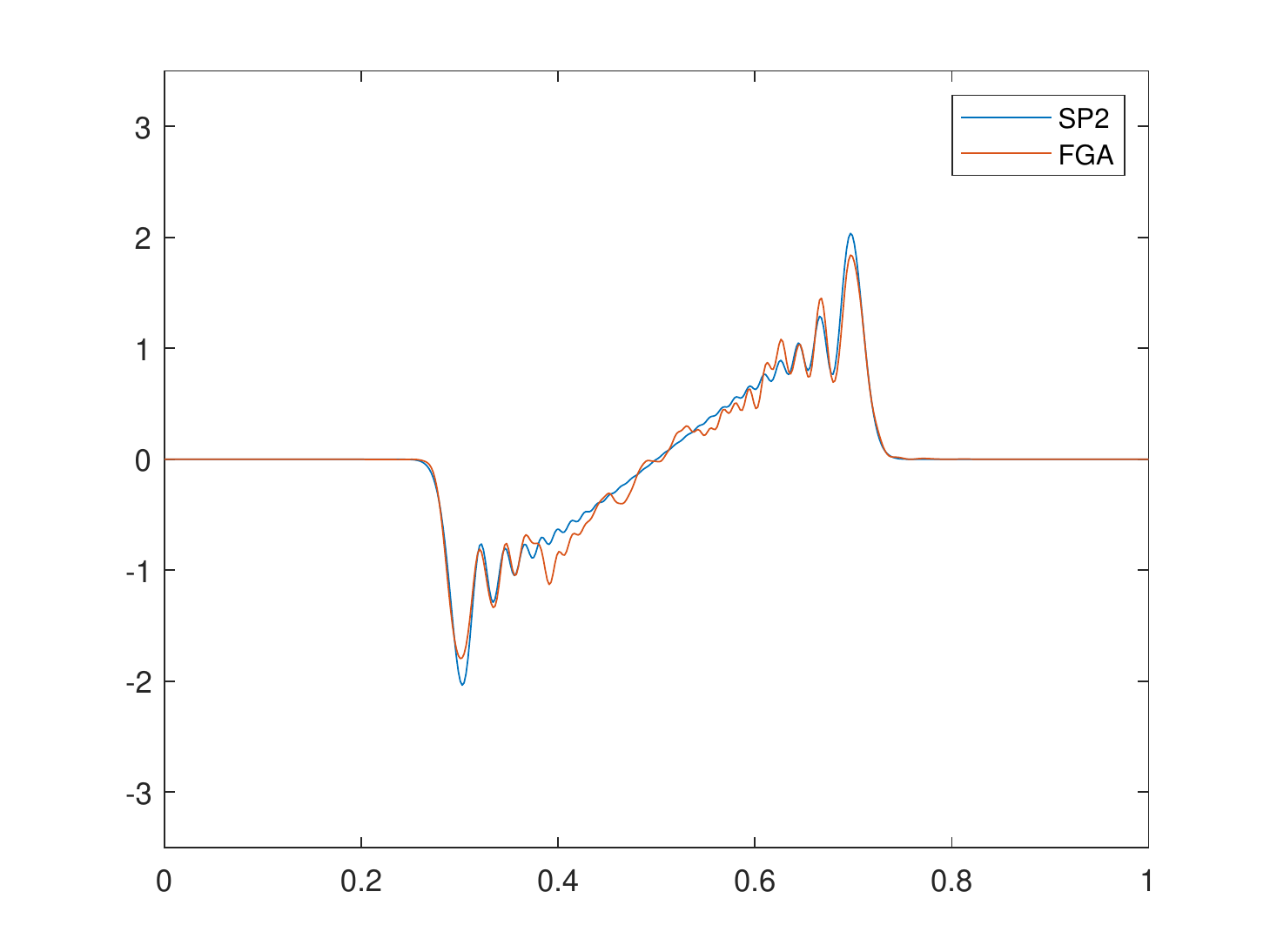}
\includegraphics[width=6cm,height=5cm]{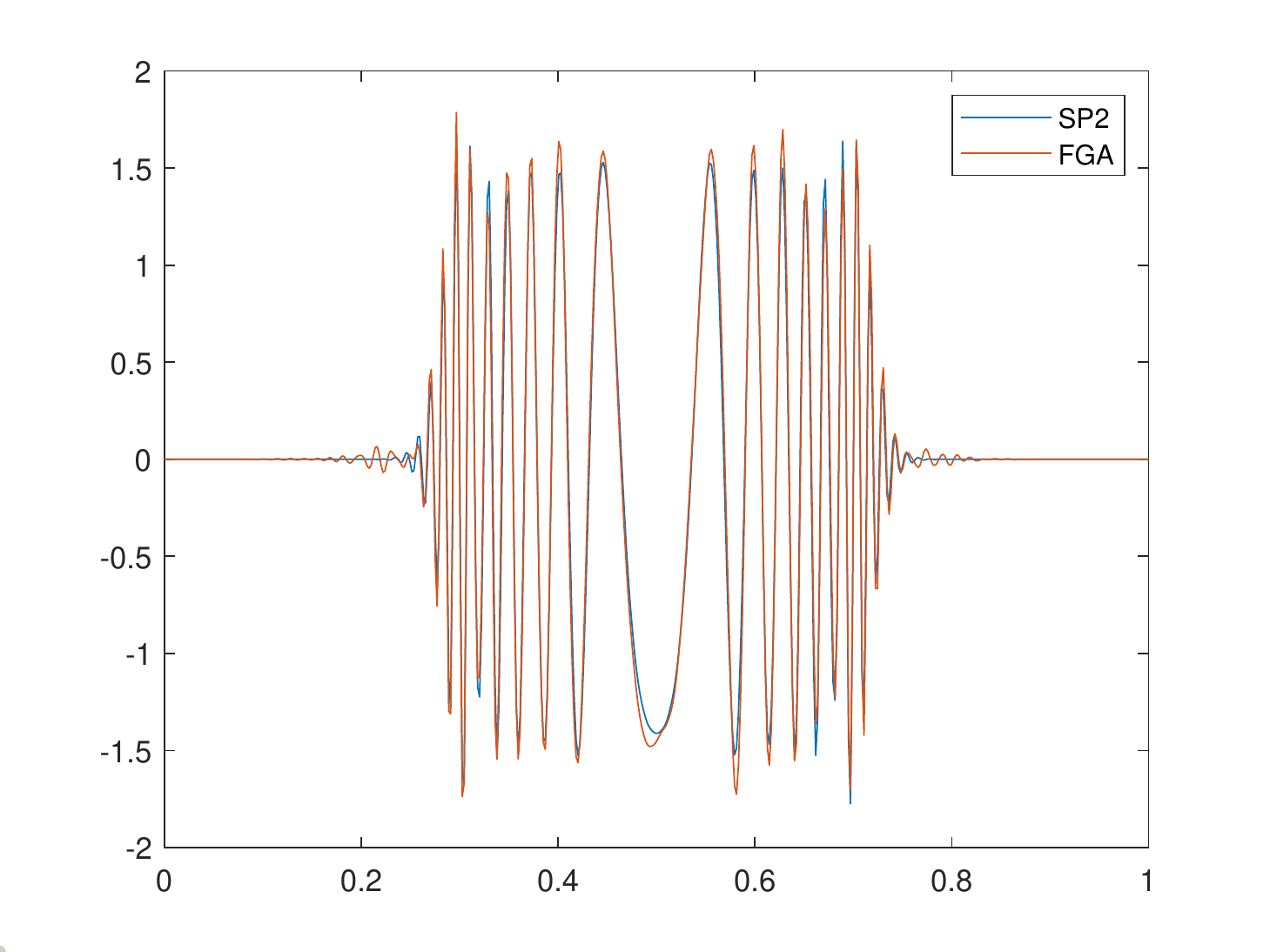}
\includegraphics[width=6cm,height=5cm]{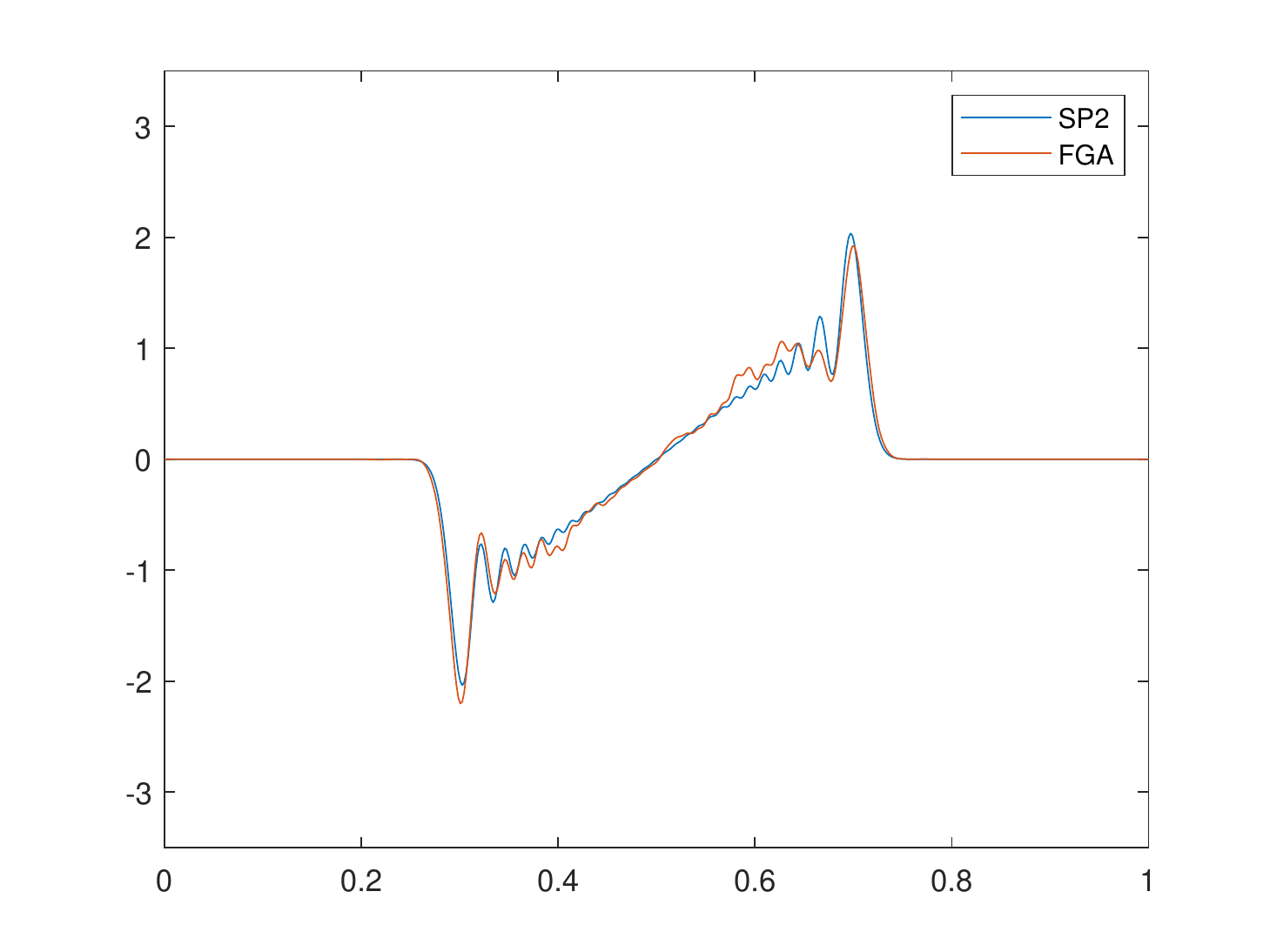}
\caption{\textbf{(Example 3)} We plot the wave functions and resulting current density functions of Equation \eqref{SemiSchrodinger} computed via the FGS versus the SP2 with different sampling size \(M\). The scaling parameter is fixed to \(\ep=0.0016\). The wave functions are presented in the left column and the current density are presented in the right column. From the top row to the bottom row, the FGS ensemble size \(M=400,800,1600\).}
\label{figWKB:C2}
\end{figure}

We fix the scaling Planck parameter as \(\ep=0.0016\) and plot the FGS wave functions and current density functions with ensemble sizes \(M=400,800,1600\) in Figure \ref{figWKB:C2}. The wave functions and current density functions computed via the FGS can reproduce the same profile as the reference solutions, which ensures the sampling strategy proposed in Section 4 is reliable and of low computational cost to approximate Equation \eqref{SemiSchrodinger} with the WKB initial conditions even in the presence of caustics.

\subsubsection*{Example 4: WKB initial data example without caustics} 

Now that we have shown the reliable results of the FGS in computing complicated WKB initial data problems with caustics, we further quantitatively study the sampling error of the FGS with the new initial sampling method in Section 4 as we have done in Example 1. In the WKB initial condition \eqref{WKBInitial}, let 
\[
a_{\ti}(x)=\left(\frac{50}{\pi}\right)^{\frac{1}{4}}\ee^{-25(x-0.5)^{2}}, \quad S_{\ti}(x)=-x^2, \quad x \in \mathbb{R},
\]
we obtain the following WKB initial data:
\begin{equation}\label{WKBini1DC}
u_{\ti}(x)=\left(\frac{50}{\pi}\right)^\frac{1}{4}\exp\left[-25\left(x-\frac{1}{2}\right)^2+\frac{\ii x^2}{\varepsilon}\right].
\end{equation}
We consider the null potential function \(E_6(x)=0\) and the quadratic potential function \(E_7(x)=\frac{x^2}{2}\). These examples were used in \cite{JWY2008} and caustics does not form.

Similar to Example 1 in Section 5.1, we aim to compute the sampling error \(E_S\) defined in Equation \eqref{ES2}. The FGA ansatz \(u_\fga\) in Equation \eqref{ES2} is approximated by a FGS wave function with sufficiently large ensemble size \(M_0=1.5\times 10^5\) as shown in Equation \eqref{FGScheck}. We use the fourth order Runge-Kutta method to solve the ODE system \eqref{PhaseQ}-\eqref{PhaseA} in the time interval \([0,0.5]\) and reconstruct the wave function \(\ufgs\) at \(t=0.5\). The computation is carried out on spatial interval \([-\pi,\pi]\).

\begin{figure}
\centering
\includegraphics[width=6cm,height=5cm]{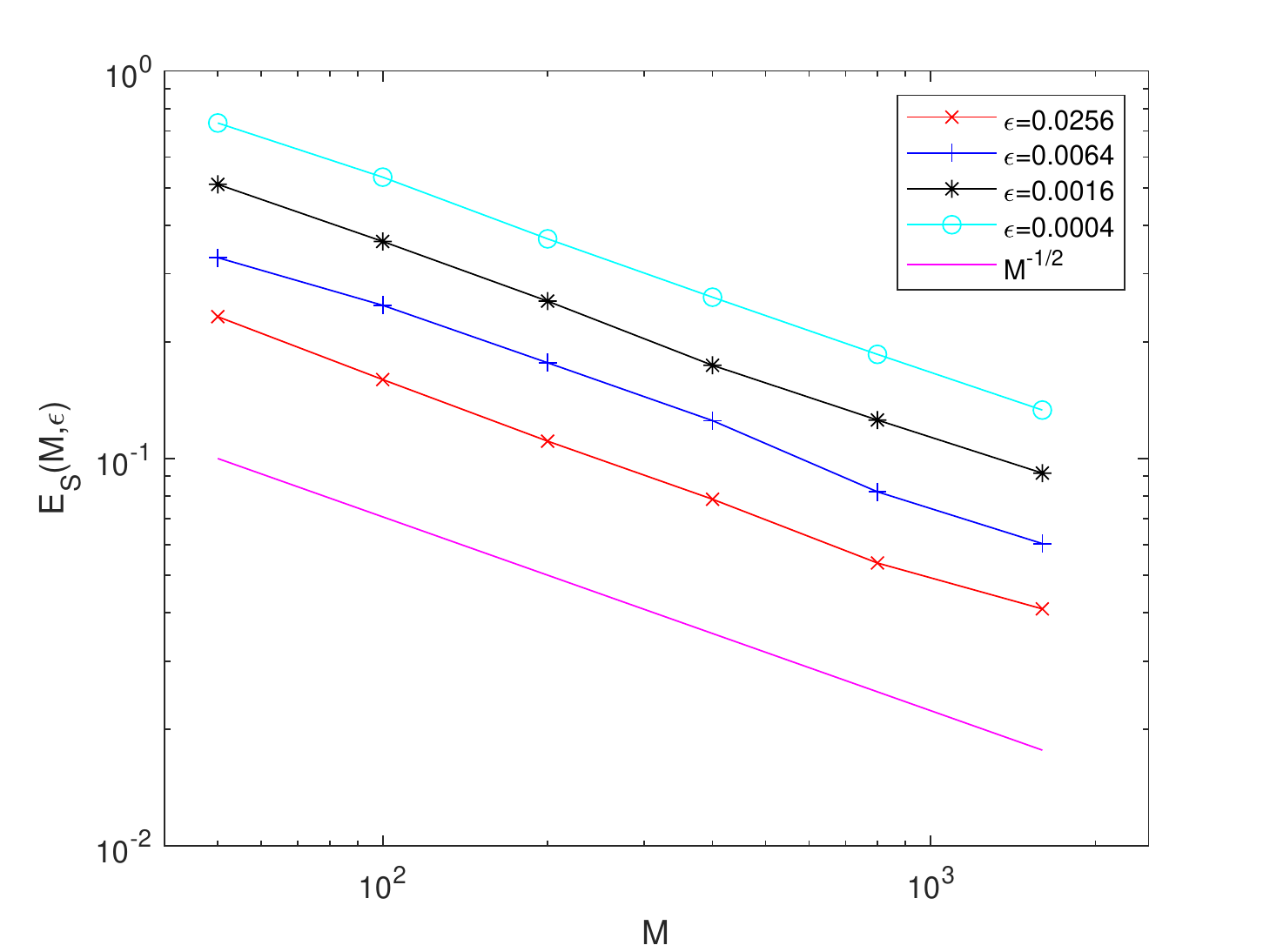}
\includegraphics[width=6cm,height=5cm]{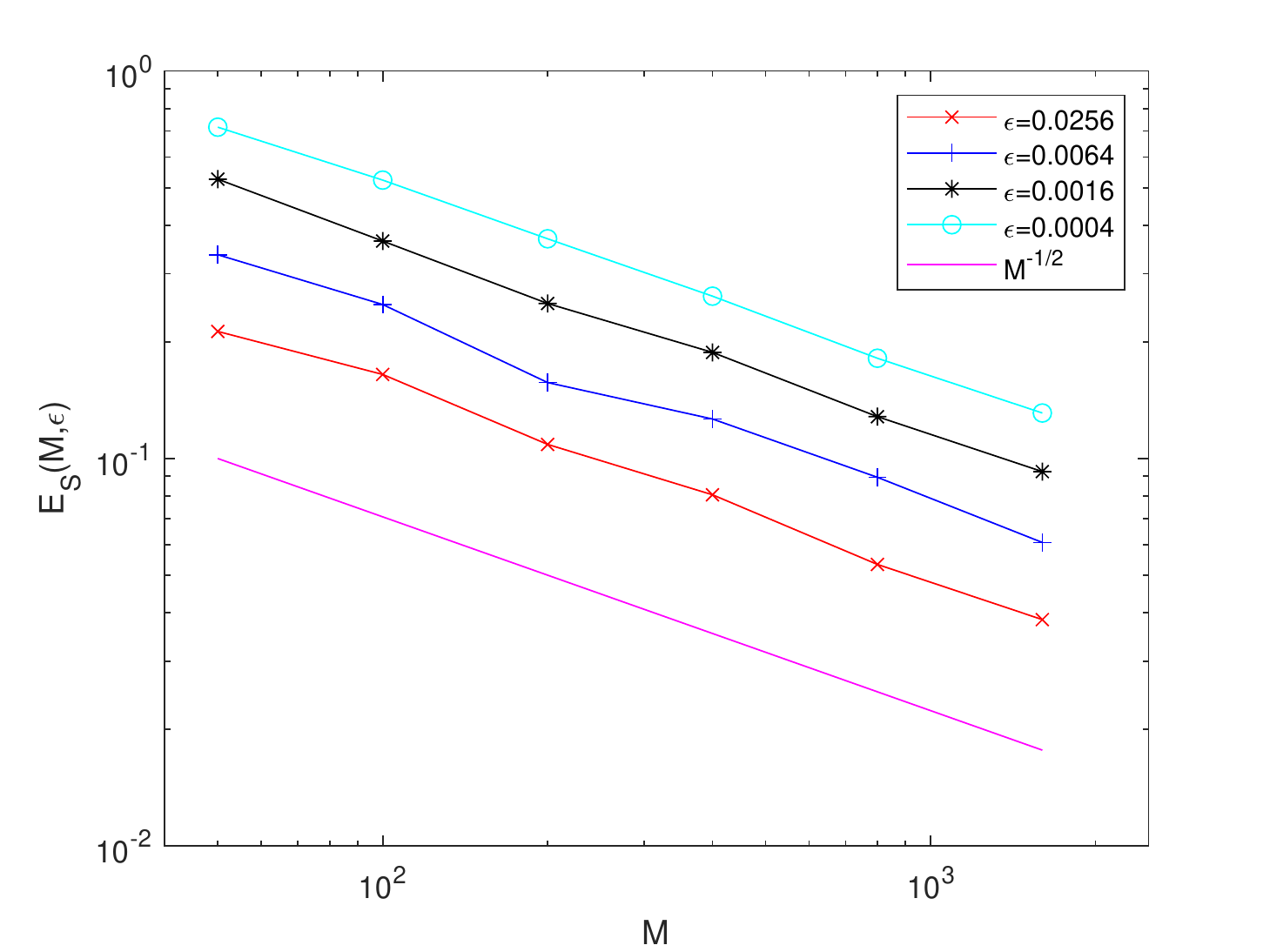}
\caption{\textbf{(Example 4) }The convergence of \(E_S\left(M,\varepsilon\right)\) with respect to the sampling size \(M\) with different scaling parameters \(\ep\) and potential functions. Left: null potential \(E_6\); Right: quadratic potential \(E_7\).}
\label{figlog2}
\end{figure}

In Table \ref{TabWKB} and Figure \ref{figlog2}, we present the sampling errors \(E_S(M,\ep)\) of the FGS with for different scaling Planck parameter \(\ep=0.0256,0.0064,0.0016,0.0004\) and two potential functions \(E_6\) and \(E_7\) respectively. Similar to Example 1, the convergence of the sampling error \(E_S\) is of order \(\CO(M^{-\frac{1}{2}})\) as \(M\to+\infty\). Moreover, we observe that the FGS sample errors of the WKB initial data cases mildly increase as the scaling parameter \(\ep\to 0\) for both potential functions. Nevertheless, we underline that the increase of \(E_S\) is insensitive to \(\ep\), which is manifested by the results in Table \ref{TabWKB}: when \(\ep\) becomes \(\frac{1}{64}\) times smaller, the errors grow at most a couple times larger. Compared to the mesh-based algorithm whose computational cost is of order \(\CO(\ep^{-1})\), the increase of the sampling error of the FGS is rather acceptable.

\begin{center}
  \begin{tabular}{c|cccc}
   \toprule
   \textbf{Potential \(E_6\)} & \(\varepsilon=0.0256\) & \(\varepsilon=0.0064\) & \(\varepsilon=0.0016\) & \(\varepsilon=0.0004\)  \\
   \midrule
   \(M=50\) & 2.32e-01 & 3.30e-01 & 5.10e-01 & 7.34e-01  \\
   \(M=100\) & 1.60e-01 & 2.49e-01 & 3.63e-01 & 5.32e-01  \\
   \(M=200\) & 1.11e-01 & 1.77e-01 & 2.55e-01 & 3.69e-01  \\
   \(M=400\) & 7.85e-02 & 1.25e-01 & 1.74e-01 & 2.61e-01  \\
   \(M=800\) & 5.37e-02 & 8.21e-02 & 1.26e-01 & 1.86e-01  \\
   \(M=1600\) & 4.09e-02 & 6.03e-02 & 9.17e-02 & 1.33e-01   \\
   \midrule
   \textbf{Potential \(E_7\)} & \(\varepsilon=0.0256\) & \(\varepsilon=0.0064\) & \(\varepsilon=0.0016\) & \(\varepsilon=0.0004\)  \\
   \midrule
   \(M=50\) & 2.13e-01 & 3.36e-01 & 5.26e-01 & 7.16e-01  \\
   \(M=100\) & 1.65e-01 & 2.50e-01 & 3.64e-01 & 5.23e-01  \\
   \(M=200\) & 1.08e-01 & 1.57e-01 & 2.51e-01 & 3.69e-01  \\
   \(M=400\) & 8.06e-02 & 1.26e-01 & 1.88e-01 & 2.62e-01  \\
   \(M=800\) & 5.33e-02 & 8.94e-02 & 1.28e-01 & 1.81e-01  \\
   \(M=1600\) & 3.84e-02 & 6.07e-02 & 9.25e-02 & 1.31e-01   \\
   \bottomrule
  \end{tabular}
 \captionof{table}{\textbf{(Example 4) }The sampling error \(E_S\left(M,\varepsilon\right)\) of the FGS with different sampling size \(M\) and scaling parameters \(\ep\).}
 \label{TabWKB}
 \end{center}

\subsection{Computing physical observables in high dimensions}

\subsubsection*{Example 5: High dimensional example}

In this subsection, we aim to compute the position observable and momentum observable as follows:
\begin{eqnarray}
O^q&=&\langle u|\hat{O}^q|u\rangle,\quad \hat{O}^q=x,\\
O^p&=&\langle u|\hat{O}^p|u\rangle,\quad \hat{O}^p=-\ii\ep\nabla_x.
\end{eqnarray}
where \(u\) denotes the exact solution to semiclassical Schr\"odinger equation \eqref{SemiSchrodinger}. According to Section 2.3, the FGS provides a mesh-free algorithm to compute the position observable and momentum observable without reconstructing the FGS wave function \(\ufgs\). The physical observables computed via the FGS are denoted as follows:
\begin{eqnarray}
\label{OQF}O^q_\fga&=&\langle \ufgs|\hat{O}^q|\ufgs\rangle,\quad \hat{O}^q=x,\\
\label{OPF}O^p_\fga&=&\langle \ufgs|\hat{O}^p|\ufgs\rangle,\quad \hat{O}^p=-\ii\ep\nabla_x,
\end{eqnarray}
where we refer to Equation \eqref{EqDoubleS} for  the computing expressions of \(O^q_\fga\) and \(O^p_\fga\). We compute the mean-square error of the quantities of the physical observables: 
\begin{eqnarray}
\label{Oqdef}\CE^q(M,m)&=&\BE\left|O^q-O^q_\fgs\right|^2,\\
\label{Opdef}\CE^p(M,m)&=&\BE\left|O^p-O^p_\fgs\right|^2,
\end{eqnarray}
where we view the observables \(O^q_\fgs\) and \(O^p_\fgs\) computed via the FGS as a random variables, whose probabilistic distribution is determined by the FGS random samples \(\left\{z_0^{(j)}\right\}_{j=1}^m\) on the phase space. We aim to discuss the mean-square errors \(\CE^q\) and \(\CE^p\) with respect to the semiclassical parameter \(\ep\) and the ensemble size \(M\).

To approximate the mean-square errors, we compute
\begin{eqnarray}
\label{OqdefC}\CE^q(M,m)&\approx&\frac{1}{\CM}\sum_{j=1}^\CM\left|O^q-O^{q(j)}_\fgs\right|^2,\\
\label{OpdefC}\CE^p(M,m)&\approx&\frac{1}{\CM}\sum_{j=1}^\CM\left|O^p-O^{p(j)}_\fgs\right|^2,
\end{eqnarray}
where \(O^{q(j)}_\fgs\) and \(O^{p(j)}_\fgs\) are sampled by independent phase space samples. We choose the ensemble size \(\CM=60\) in our tests. To compute  the exact quantities \(O^q\) and \(O^p\) of the position and momentum observables which are analytically unavailable, we consider a special case when the quantities of \(O^q\) and \(O^p\) can be given explicitly by a semiclassical Gaussian wave packet named the Hagedorn wave packet \cite{H1980,H1998,Z2014,LL2020}, which can be approximated numerically through evolving ODEs. For this special case, the position and momentum observables involved to the wave packet can be derived explicitly through the wave packet parameters. We directly quote Proposition 3.18 in \cite{LL2020} as the theoretical foundation of the Hagedorn wave packet method:

\begin{theorem}
Consider the semiclassical Schr\"odinger equation \eqref{SemiSchrodinger} with a quadratic potential \(E(x)=\frac{|x|^2}{2}\). Let vectors \(q_H(t),p_H(t)\in\BR^m\) and matrices \(Q_H(t),P_H(t)\in\BC^{m\times m}\) be the solution to the following ODE system:
\begin{eqnarray}
\label{qH}\dot{q_H}&=&p_H,\\
\label{pH}\dot{p_H}&=&-\nabla E(q_H),\\
\label{QH}\dot{Q_H}&=&P_H,\\
\label{PH}\dot{P_H}&=&-\nabla^{2} E(q_H) Q_H.
\end{eqnarray}
Let \(S_H(t)=\int_{0}^{t}\left(\frac{1}{2}|p_H(s)|^{2}-E(q_H(s))\right) \dd s\) be the corresponding action function. Then the following complex Gaussian wave packet
\begin{eqnarray}\label{HagedornWave}
u_H(t,x)&=& \ee^{\frac{\ii S_H(t)}{\ep}}(\pi \varepsilon)^{-\frac{m}{4}}(\mathrm{det} Q_H(t))^{-\frac{1}{2}}\times\nonumber\\&&\exp \left(\frac{\mathrm{i}}{2 \varepsilon}(x-q_H(t))^{T} P_H(t) Q_H(t)^{-1}(x-q_H(t))+\frac{\ii}{\varepsilon} p_H(t)^{T}(x-q_H(t))\right),\nonumber\\&&
\end{eqnarray}
is an exact solution to equation \eqref{SemiSchrodinger} of unit \(L^2\) norm for any \(t\in \BR^+\) provided that:
\begin{enumerate}[itemindent=0em]
\item the initial matrices \(P(0)\) and \(Q(0)\) satisfies the symplecticity relationship, i.e.
\begin{eqnarray}
\label{symp1}Q_H(0)^{T} P_H(0)-P_H(0)^{T} Q_H(0)&=&0, \\
\label{symp2}Q_H(0)^{*} P_H(0)-P_H(0)^{*} Q_H(0)&=&2 \ii\mathrm{I}_m,
\end{eqnarray}
where \(\mathrm{I}_m\) denotes the \(m\)-th order identity matrix.
\item The initial condition of \(u_H(t,x)\) is  exact, i.e. \(u_H(0,\cdot)=u_\ti\).
\end{enumerate}
\end{theorem}

To apply the Hagedorn wave packet method to compute the exact values of the observables, we only consider the semiclassical Gaussian wave packet initial data \eqref{InitialGaussian} and fix the potential function as \(E(x)=\frac{|x|^2}{2}\) in this example. Moreover, in order to ensure \(\kl q_H(0),p_H(0),Q_H(0),P_H(0)\kr\) reproduce the Gaussian initial condition \eqref{InitialGaussian} exactly and satisfies the symplecticity relationship \eqref{symp1} and \eqref{symp2}, the initial condition is chosen to be 
\begin{equation}
q_H(0)=\tilde{q},\quad p_H(0)=\tilde{p},\quad Q_H(0)=\ii A_m^{-\frac{1}{2}}U,\quad P_H(0)=A_m^{\frac{1}{2}}U,
\end{equation} 
where \(A_m=\mathrm{diag}\{a_1,a_2,\cdots,a_m\}\) and \(U\) denotes any unitary matrix. With the preparations above, the exact solution to Equation \eqref{SemiSchrodinger} is of the form \(u_H(t,x)\) given in Equation \eqref{HagedornWave}, whose parameters are numerically available through evolving the ODE system \eqref{qH}-\eqref{PH}. Now recall that the Hagedorn wave packet \(u_H(x,t)\) is a semiclassical Gaussian wave packet, hence the quantities of the position and momentum observables are analytically available:
\ben
\label{eoq}O^q_H&=&\langle u_H|\hat{O}^q|u_H\rangle=\frac{\left[\det\kl\mathrm{Im}\kl P_HQ_H^{-1}\kr\kr\right]^{-\frac{1}{2}}}{\det\kl Q_H\kr}q_H,\\
\label{eop}O^p_H&=&\langle u_H|\hat{O}^p|u_H\rangle=\frac{\left[\det\kl\mathrm{Im}\kl P_HQ_H^{-1}\kr\kr\right]^{-\frac{1}{2}}}{\det\kl Q_H\kr}p_H,
\een
where \(\mathrm{Im}\) denotes the imaginary part of the matrix. As a result, the observables \(O^q_H\) and \(O^p_H\) computed via the Hagedorn wave packet method serves as an accurate and efficient reference to the exact quantities of the physical observables in high-dimensional cases.

Now we move to design our numerical experiments. In the Gaussian initial condition \eqref{InitialGaussian}, let the initial position and momentum to be
\beq
\tilde{q}=\frac{1}{\sqrt{m}}\kl 1,1,\cdots,1\kr\in\BR^m,
\eeq
and
\beq
\tilde{p}=-\frac{1}{\sqrt{m}}\kl 1,1,\cdots,1\kr\in\BR^m.
\eeq
Moreover, the initial amplitude  \(\kl a_1,a_2,\cdots,a_m\kr\in\BR^m\) is given as follows
\beq
a_j=1+\{0.2j\},\quad j=1,2,\cdots,m,
\eeq
where \(\{x\}\) denotes the fractional part of a rel number \(x\). 

\brm 
As shown in Equation \eqref{CES}, the sampling error of the FGS is relevant to the term \(\prod_{j=1}^m\kl\frac{1+a_j}{\sqrt{a_j}}\kr^{\frac{1}{2}}\), which is related to the initial amplitude. To avoid the influence resulting from the initial conditions, we choose \(a_j\) to be restricted in the interval \([1,2)\).
\erm

We aim to compute the mean-square error \(\CE^q\) and \(\CE^p\) given by Equation \eqref{OqdefC} \eqref{OpdefC} at \(t=0.5\), where the exact quantities of the physical observables \(O^q,O^p\) are approximate via Equation \eqref{eoq} and \eqref{eop} associated with  the Hagedorn wave packets. We use the fourth order Runge-Kutta method to solve Equation \eqref{qH}-\eqref{PH} associated with the Hagedorn wave packets and Equation \eqref{PhaseQ}-\eqref{PhaseA} associated with the FGS with a sufficiently small step length \(\Delta t=0.01\) on time interval \([0,0.5]\). The computation is carried out on spatial interval \([-\pi,\pi]\) and the scaling parameter is fixed to \(\ep=0.0256\).

\begin{figure}
\centering
\includegraphics[width=6cm,height=5cm]{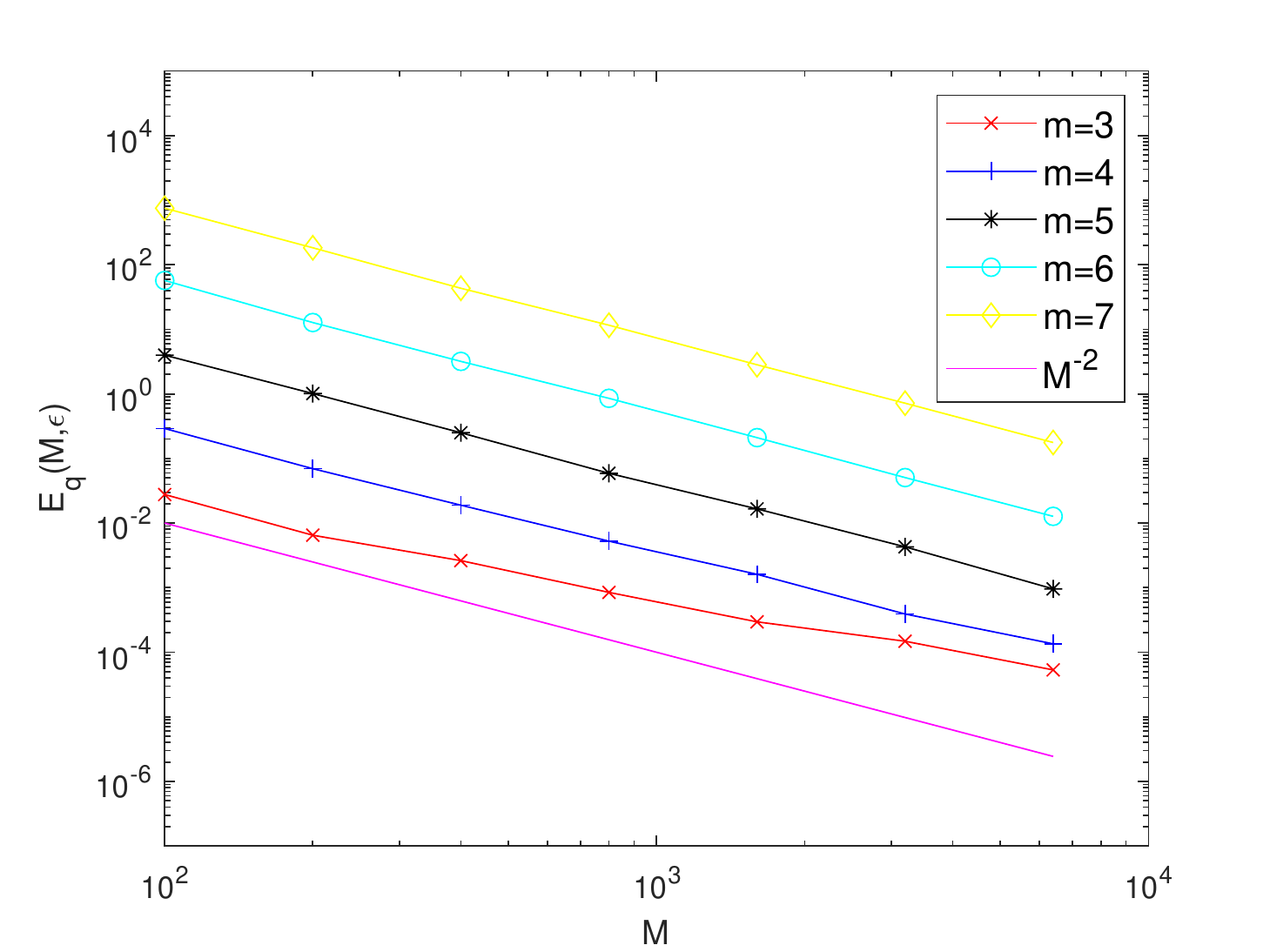}
\includegraphics[width=6cm,height=5cm]{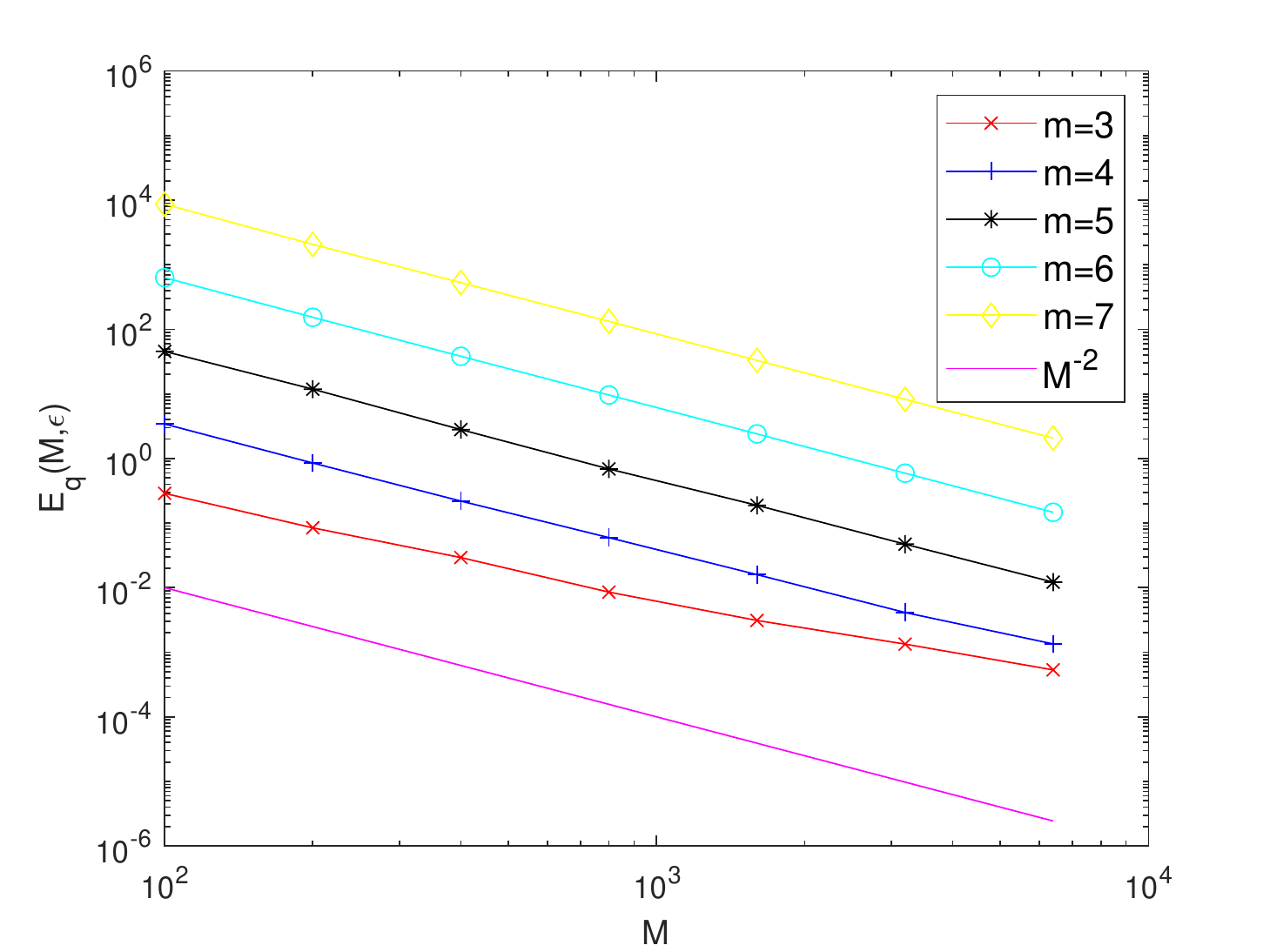}
\caption{\textbf{(Example 5) }We choose \(\ep=0.0256\). We plot the observable errors \(\CE^q\) and \(\CE^p\) (defined in Equation \eqref{Oqdef} and \eqref{Opdef}) with different sampling sizes \(M\) and dimension numbers \(m\). Left: position observable error \(\CE^q\); Right: momentum observable error \(\CE^p\).}
\label{fighighD}
\end{figure}

In figure \ref{fighighD}, we plot the position and momentum observable errors \(\CE^q\) and \(\CE^p\) with different ensemble sizes \(M\) and dimension number \(m=3,4,\cdots,7\). For fixed \(m\), the convergence rate of \(\CE^q\) and \(\CE^p\) is roughly \(\CO(M^{-2})\), whereas the convergence slows down when \(M\) is large for \(m=3\). This is due to the fact that, when \(M\) is sufficiently large, the Monte Carlo sampling errors becomes much smaller than the asymptotic errors of the FGA ansatz that dominate the observable errors \(\CE^q\) and \(\CE^p\). Moreover, we observe that, with fixed ensemble size \(M\), the errors of the physical observables mildly increase as the dimension number \(m\) increases. As is manifested in \label{Tabobs}, in order to reached a uniformly bounded mean-square error, the sample size should quadruple as the dimension number increases by one. Therefore, it is recommended that the sample size to compute physical observable in a \(m\)-dimensional problem should be of the order \(\CO(4^m)\). 

Furthermore, we record in Table \ref{Tabobs2} the system time to implement the FGS to compute an \(m-\)dimensional physical observables \(O^q_{\fga}\) or \(O^p_{\fga}\). Since the complexity to compute the position and momentum observables are almost the same, we only present the system time to compute \(O^q_{\fga}\) in Table \ref{Tabobs2}. To reached a uniformly bounded mean-square error, the system time increases at five times as  \(m\) increases by one.

\begin{center}
  \begin{tabular}{c|c|c|c|c}
   \toprule
   \multirow{8}{*}{\(\CE^q\)}&
   \multirow{2}{*}{m=6} & \(M=12800\) & \(M=6400\) & \(M=3200\)  \\\cline{3-5}
   & &3.23e-03  &1.27e-02  &5.07e-02  \\\cline{2-5}
   &\multirow{2}{*}{m=5} & \(M=3200\) & \(M=1600\) & \(M=800\) \\\cline{3-5}
   & &4.29e-03  &1.66e-02  &5.90e-02  \\\cline{2-5}
   &\multirow{2}{*}{m=4} & \(M=800\) & \(M=400\) & \(M=200\) \\\cline{3-5}
   & &5.24e-03 & 1.89e-02 &6.96e-02   \\\cline{2-5}
   &\multirow{2}{*}{m=3} & \(M=200\) & \(M=100\) & \(M=50\) \\\cline{3-5}
   & &6.49e-03  &2.78e-02  &7.28e-02  \\
   \midrule
   \multirow{8}{*}{\(\CE^q\)}&
   \multirow{2}{*}{m=6} & \(M=12800\) & \(M=6400\) & \(M=3200\)  \\\cline{3-5}
   & &3.76e-02  &1.47e-01  &5.93e-01  \\\cline{2-5}
   &\multirow{2}{*}{m=5} & \(M=3200\) & \(M=1600\) & \(M=800\) \\\cline{3-5}
   & &4.72e-02  &1.89e-01  &6.89e-01  \\\cline{2-5}
   &\multirow{2}{*}{m=4} & \(M=800\) & \(M=400\) & \(M=200\) \\\cline{3-5}
   & &5.96e-02 &2.19e-01  &8.52e-01   \\\cline{2-5}
   &\multirow{2}{*}{m=3} & \(M=200\) & \(M=100\) & \(M=50\) \\\cline{3-5}
   & &8.45e-02  &2.89e-01  &9.08e-01  \\
   \bottomrule
  \end{tabular}
  
 \captionof{table}{\textbf{(Example 5) }The table presents the observable errors \(\CE^q\) and \(\CE^p\) (defined in Equation \eqref{Oqdef} and \eqref{Opdef}) with different sampling sizes \(M\) and dimension numbers \(m\). }
 \label{Tabobs}
 \end{center}
 
 \begin{center}
  \begin{tabular}{c|c|c|c|c}
   \toprule
   \multirow{8}{*}{system time}&
   \multirow{2}{*}{m=6} & \(M=12800\) & \(M=6400\) & \(M=3200\)  \\\cline{3-5}
   & &251.2  &119.6  &58.8  \\\cline{2-5}
   &\multirow{2}{*}{m=5} & \(M=3200\) & \(M=1600\) & \(M=800\) \\\cline{3-5}
   & &50.8  &25.3  &12.0  \\\cline{2-5}
   &\multirow{2}{*}{m=4} & \(M=800\) & \(M=400\) & \(M=200\) \\\cline{3-5}
   & &10.9 & 5.5 &2.8   \\\cline{2-5}
   &\multirow{2}{*}{m=3} & \(M=200\) & \(M=100\) & \(M=50\) \\\cline{3-5}
   & &2.8  &1.4  &0.8  \\
   \bottomrule
  \end{tabular}
  
 \captionof{table}{\textbf{(Example 5) }The table presents the system times to compute the physical observable \(O^q_{\fga}\) or \(O^p_{\fga}\) (defined in Equation \eqref{OQF} and \eqref{OPF}) via the FGS with different sampling sizes \(M\) and dimension numbers \(m\).}
 \label{Tabobs2}
 \end{center}

\section{Conclusion and discussion}

In this work, we have developed a novel stochastic simulation methodology for the semiclassical  Schr\"odinger equation that shows great promises for high dimensional simulations. Its numerical cost is insensitive to the  semiclassical parameter, grows gently with respect to the dimension number and it provides an efficient way to compute the physical observables without reconstructing the wave function. The core component of the algorithm is to properly choose a density function, that achieves a balance between variance reduction and sampling convenience, and should be specified for different types of initial wave functions. Although we have already obtained satisfactory results for the Gaussian initial conditions, and the WKB initial conditions, more advanced sampling techniques, in theory, can be utilized in the framework of the FGS methodology, and hence further improve its performance. Besides, the idea of turning computing highly oscillatory waves into sampling and propagating wave packets may also be applied to other wave equations in the high frequency regime.  Thus, this work initiates many directions that we shall explore in the future.  

\section*{\textbf{Appendix A: }The nondimensionalisation of the Schr\"odinger equation}


When considering equations with dimensional physical quantities, one can scale all the variables by dimensional constant to yield dimensionless equations, named the process of nondimensionalisation, effectively leaving a "clean" mathematical system with dimensionless magnitudes and dimensionless parameters.

In the case of Schr\"odinger equation \eqref{physch}, we choose the spatial and temporary characteristic scales \(x_c\) and \(t_c\) in order to make the scaled variable \(\bar{x}=\frac{x}{x_c}\) and \(\bar{t}=\frac{t}{t_c}\) dimensionless. By taking \(\bar{u}(\bar{t},\bar{x})=u(t,x)\), Equation \eqref{physch} can be rewritten as
\beq\label{apd1}
\ii \hbar \frac{1}{t_{c}} \partial_{\bar{t}} \bar{u}=-\frac{\hbar^{2}}{2 m_o x_{c}^{2}} \partial_{\bar{x} \bar{x}} \bar{u}+E(x) \bar{u}.
\eeq
Then we multiply \(\frac{\kl t_c\kr^2}{m_o\kl x_c\kr^2 }\) on both sides of Equation \eqref{apd1}:
\beq
\ii\left(\frac{\hbar t_{c}}{m_o x_{c}^{2}}\right) \partial_{\bar{t}} \bar{u}=-\frac{1}{2}\left(\frac{\hbar t_{c}}{m x_{c}^{2}}\right)^{2} \partial_{\bar{x} \bar{x}} \bar{u}+\frac{t_{c}^{2}}{m x_{c}^{2}} E\left(\bar{x} x_{c}\right) \bar{u}.
\eeq
Note that the parameter \(\ep= \frac{\hbar t_{c}}{m x_{c}^{2}}\) and potential \(\bar{E}(\bar{x})=\frac{t_{c}^{2}}{m x_{c}^{2}} E\left(\bar{x} x_{c}\kr\) are dimensionless, we obtain the dimensionless Schr\"odinger equation:
\beq\label{apd3}
\ii\ep\partial_{\bar{t}} \bar{u}=-\frac{\ep^2}{2}\bar{u}+\bar{E}\kl\bar{x}\kr\bar{u}.
\eeq
Finally, Equation \eqref{apd3} is actually identical with the semiclassical Schr\"odinger equation \eqref{SemiSchrodinger} if we remove the bars in Equation \eqref{apd3}.

\section*{\textbf{Appendix B: } The WKB approximation and the caustics}

The WKB approximation is an asymptotic description of  the exact solution to the semiclassical Schr\"odinger equation \eqref{SemiSchrodinger}. To match the \(\CO(\ep)\) frequency oscillation structure of the wave function, we give the following WKB ansatz:
\beq\label{WKBA}
u(t,x)\sim u_W(t,x)= a_W(t,x)\exp\kl\frac{\ii}{\ep} S_W\kr,\quad\ep\to 0,
\eeq
where \(a_W(t,x)\) and \(S_W(t,x)\) are two real-valued functions. Through asymptotic matching at order \(\CO(1)\) and \(\CO(\ep)\), we derive the eikonal equation:
\beq
\label{WKB2}\partial_{t} S_W+\frac{1}{2}\left|\nabla S_W \right|^{2}+E(x)=0, \quad S_W(0, x)=S_{\ti }(x),
\eeq
and the transport equation:
\beq
\label{WKB1}\partial_{t} a_W+\nabla S_W \cdot \nabla a_W+\frac{a_W}{2} \Delta S_W=0, \quad a_W(0, x)=a_{\ti}(x),
\eeq
which solve the amplitude function \(a_w(t,x)\) and the phase function \(S_W(t,x)\).

Note that the eikonal equation \eqref{WKB2} is a nonlinear Hamiltonian-Jacobi equation and can be solved through the method of characteristic. The characteristic flow is given by the following Hamiltonian system of ODEs:
\beq\label{ODEsys}
\left\{\begin{array}{l}
\dot{x_W}(t, y)=\xi_W(t, y), \quad x_W(0, x)=y, \\
\dot{\xi_W}(t, y)=-\nabla_{x} E(x_W(t, y)), \quad \xi_W(0, y)=\nabla S_{\ti}(y).
\end{array}\right.
\eeq
With the characteristic flow \(X_t:y\mapsto x_W(r,y)\), the ODE system \eqref{ODEsys} yields a phase function as follows:
\beq\label{SW}
S_W(t, x)=S_\ti(x)-\int_{0}^{t}\left[ \frac{1}{2}\left|\nabla S_W\kl\tau, X_\tau^{-1}(x)\kr\right|^{2}+E\kl X_\tau^{-1}(x)\kr\right] \dd \tau.
\eeq
Note that the characteristic flow \(X_t:y\mapsto x_W(r,y)\) may not be one-by-one, hence the phase function \(S_W\) given in Equation \eqref{SW} may be smooth only in a small time interval when characteristic flow \(X_t:y\mapsto x_W(r,y)\) doesn't cross. The cases when \(X_t\) ceases to be a diffeomorphism are called caustics phenomena. See Figure 2 in \cite{GJL2003} for an illustration of the caustics characteristic.

Before caustics onset, the WKB ansatz in Equation \eqref{WKBA} is an asymptotic approximation to the exact solution to Equation \eqref{SemiSchrodinger} of order \(\CO(\ep\)), provided that the initial considition is exact and the amplitude \(a_\ti\) decays rapidly (see Theorem 2.3 in \cite{JMS2011}). After caustics forms, the exact solution may no longer be approximated by the WKB ansatz \eqref{WKBA}, hence the WKB approximation breaks down. The readers may refer to \cite{C2008} for detailed discussion on the WKB approximation and analytic techniques after caustics onset.

\section*{Acknowledgment}
Z. Zhou is supported by the National Key R\(\&\)D Program of China, Project Number 2020YFA0712000, 2021YFA1001200  and NSFC grant Number 12031013, 12171013. Z. Zhou thanks Shi Jin for helpful discussions.


\end{document}